\documentclass[10pt,draftclsnofoot,letterpaper,onecolumn]{IEEEtran}
\usepackage{amsmath,amsfonts,amssymb}%,amsbsy,mathbbol}
\usepackage[multiple]{footmisc}

\newtheorem{assumption}{Assumption}

\newtheorem{theorem}{Theorem}

\newtheorem{lemma}{Lemma}

\newtheorem{proposition}{Proposition}

\newtheorem{corollary}{Corollary}

\renewcommand{\theenumi}{\roman{enumi}}

\newtheorem{remark}{}

\begin{document}

\title{Analyticity, Convergence and Convergence Rate of 
Recursive Maximum Likelihood Estimation in 
Hidden Markov Models} 

\author{Vladislav B. Tadi\'{c}$^{1}$\thanks{$^{1}$Department of Mathematics, University of Bristol,
University Walk, Bristol BS8 1TW, United Kingdom. 
(v.b.tadic@bristol.ac.uk). }}

\maketitle

\begin{abstract} 
This paper considers the asymptotic behavior of the recursive maximum likelihood 
estimation in hidden Markov models. 
The paper is focused on the analytic properties 
of the asymptotic log-likelihood and on the 
point-convergence and convergence rate of
the recursive maximum likelihood estimator. 
Using the principle of analytical continuation, 
the analyticity of the asymptotic log-likelihood is 
shown for analytically parameterized hidden Markov models. 
Relying on this fact and some results from differential geometry
(Lojasiewicz inequality), 
the almost sure point-convergence of 
the recursive maximum likelihood algorithm is demonstrated, 
and relatively tight bounds on the convergence rate are derived. 
As opposed to the existing result on the asymptotic behavior of 
maximum likelihood estimation in hidden Markov models, 
the results of this paper are obtained without assuming that 
the log-likelihood function has an isolated maximum at which the Hessian 
is strictly negative definite. 
\end{abstract}

\begin{IEEEkeywords} 
Hidden Markov models, maximum likelihood estimation, recursive identification, 
analyticity, Lojasiewicz inequality, point-convergence, convergence rate. 
\end{IEEEkeywords}

\section{Introduction} \label{section0}

Hidden Markov models are a broad class of stochastic processes capable of modeling 
complex correlated data and large-scale dynamical systems. 
These processes consist of two components: 
states and observations. 
The states are unobservable and form a Markov chain. 
The observations are independent conditionally on the states
and 
provide only available information about the state dynamics. 
Hidden Markov models have been formulated in the seminal paper \cite{baum&petrie}, 
and over last few decades, they have found a wide range of applications 
in diverse areas such as acoustics and signal processing, 
image analysis and computer vision, automatic control and robotics, 
economics and finance, computation biology and bioinformatics. 
Due to their practical relevance, these models have extensively been studied in 
a large number of papers and books 
(see e.g., \cite{cappe&moulines&ryden}, \cite{ephraim&merhav} 
and references cited therein). 

Besides the estimation of states given available observations
(also known as filtering), 
the identification of model parameters are probably the most important problem
associated with hidden Markov models. 
This problem can be described as the estimation (or approximation) of 
the state transition probabilities and the observations conditional distributions 
given available observations.  
The identification of hidden Markov models have been considered in numerous papers 
and several methods and algorithms have been developed 
(see \cite[Part II]{cappe&moulines&ryden}, \cite{ephraim&merhav} and references 
cited therein). 
Among them, the methods based on the maximum likelihood principle 
are probably one of the most important and popular. 
Their various asymptotic properties 
(asymptotic consistency, asymptotic normality, convergence rate)
have been analyzed in a number of papers 
(see 
\cite{baum&petrie}, 
\cite{bickel&ritov}, 
\cite{bickel&ritov&ryden}, 
\cite{douc&matias}, \cite{douc&moulines&ryden}, 
\cite{legland&mevel1} -- \cite{legland&mevel4}, 
\cite{leroux}, 
\cite{mevel&finesso}, 
\cite{ryden1}, \cite{ryden2}; 
see also \cite[Chapter 12]{cappe&moulines&ryden}, \cite{ephraim&merhav} and references 
cited therein). 
Although the existing results provide an excellent insight into 
the asymptotic behavior of maximum likelihood estimators for 
hidden Markov models, 
they all crucially rely on the assumption that the log-likelihood function 
has a strong maximum,
i.e., an isolated maximum at which the Hessian is strictly negative definite.
As the log-likelihood function admits no close-form expression and is 
fairly complex even for small-size hidden Markov models (four or more states), 
it is hard (if not impossible at all) to show the existence of an isolated maximum, 
let alone checking the definiteness of the Hessian. 

The differentiability, analyticity and other analytic properties 
of functionals of hidden Markov models similar to the asymptotic likelihood 
(mainly entropy rate) 
have recently been studied in 
\cite{han&marcus1}, 
\cite{han&marcus2}, 
\cite{holliday&goldsmith&glynn}, 
\cite{ordentlitch&weissman}, 
\cite{peres}, 
\cite{schoenhuht}. 
Although very insightful and useful,
the results presented in these papers cover only models 
with discrete state and observation spaces 
and do not consider
the asymptotic behavior of   
the maximum likelihood estimation method. 

In this paper, we study the asymptotic behavior of 
the recursive maximum likelihood estimation in hidden Markov models
with a discrete state-space and continuous observations.  
We establish a link between the analyticity of the asymptotic log-likelihood 
on one side, 
and the point-convergence and convergence rate of 
the recursive maximum likelihood algorithm, on the other side. 
More specifically, relying on the principle of analytical continuation,
we show under mild conditions that the asymptotic log-likelihood function is analytical in 
the model parameters if the state transition probabilities and 
the observation conditional distributions are analytically parameterized. 
Using this fact and some results from differential geometry 
(Lojasiewicz inequality), 
we demonstrate that the recursive maximum likelihood algorithm for hidden Markov models
is almost surely point-convergent (i.e., it has a single accumulation point 
w.p.1). 
We also derive tight bounds on the almost sure convergence rate. 
As opposed to all existing results on the asymptotic behavior of 
maximum likelihood estimation in hidden Markov models, 
the results of this paper are obtained without assuming that the 
log-likelihood function has an isolated strong maximum. 

The paper is organized as follows. 
In Section \ref{section1}, the hidden Markov models and 
the corresponding recursive maximum likelihood algorithms are defined. 
The main results are presented in Section \ref{section1}, too. 
Section \ref{section2} provides several practically relevant examples of the main results. 
Section \ref{section1*} contains the proofs of the main results, 
while the results of Section \ref{section2} are shown in Section \ref{section2*}. 

\section{Main Results} \label{section1}

In order to state the problems of recursive identification 
and 
maximum likelihood estimation in  
hidden Markov models with finite state-spaces and continuous observations, 
we use the following notation. 
$N_{x} > 1$ is an integer, while 
${\cal X} = \{1,\dots,N_{x} \}$.  
$d_{y} \geq 1$ is also an integer, while  
${\cal Y}$ is a Borel-measurable set from $\mathbb{R}^{d_{y} }$. 
$\{p(x'|x) \}_{x,x' \in {\cal X} }$ 
are non-negative real numbers such that  
$\sum_{x'\in {\cal X} } p(x'|x) = 1$
for each $x\in {\cal X}$.   
$\{Q(\cdot|x) \}_{x \in {\cal X} }$ 
are probability measures on ${\cal Y}$. 
$\{ (X_{n}, Y_{n} ) \}_{n\geq 0}$ 
is an ${\cal X} \times {\cal Y}$-valued stochastic process 
which is defined on a (canonical) probability space 
$(\Omega, {\cal F}, P)$
and 
satisfies 
\begin{align*} 
	P(Y_{n+1} \in B, X_{n+1} = x| X_{0}, Y_{0}, \dots, X_{n}, Y_{n} ) 
	=
	Q(B|x) p(x|X_{n} )
\end{align*}
w.p.1 
for all $x\in {\cal X}$, $n\geq 0$, and any 
Borel measurable set $B$ from ${\cal Y}$. 
On the other side, 
$d_{\theta }$ is a positive integer, while 
$\Theta$ is an open set from $\mathbb{R}^{d_{\theta } }$.  
$\{p_{\theta }(x'|x) \}_{x,x' \in {\cal X} }$ 
are Borel-measurable functions of $\theta \in \Theta$ 
such that 
$p_{\theta }(x'|x) \geq 0$
and
$\sum_{x''\in {\cal X} } p_{\theta }(x''|x) = 1$
for all $\theta \in \Theta$, $x,x' \in {\cal X}$. 
$\{q_{\theta }(y|x) \}_{x \in {\cal X} }$ 
are Borel-measurable functions of $(\theta, y ) \in \Theta \times {\cal Y}$ 
such that 
$q_{\theta }(y|x) \geq 0$
and 
$\int_{{\cal Y} } q_{\theta }(y'|x) dy' = 1$
for all $\theta \in \Theta$, $x\in {\cal X}$, $y \in {\cal Y}$. 
For $\theta \in \Theta$, 
$\{ (X_{n}^{\theta }, Y_{n}^{\theta } ) \}_{n\geq 0}$ 
is an ${\cal X} \times {\cal Y}$-valued stochastic process 
which is 
defined on a (canonical) probability space 
$(\Omega, {\cal F}, P_{\theta } )$
and admits 
\begin{align*} 
	P_{\theta }(Y_{n+1}^{\theta } \in B, X_{n+1}^{\theta } = x | 
	X_{0}^{\theta }, Y_{0}^{\theta }, \dots, X_{n}^{\theta }, Y_{n}^{\theta } ) 
	=
	\int_{B}
	q_{\theta }(y|x) p_{\theta }(x|X_{n}^{\theta } ) dy
\end{align*}
w.p.1 for each $x\in {\cal X}$, $n\geq 0$, 
and any Borel measurable set $B$ from {\cal Y}. 
Finally, 
$f(\cdot )$ stands for 
the asymptotic value of the log-likelihood function 
associated with data $\{Y_{n} \}_{n\geq 0}$.  
It is defined by 
\begin{align*}
	f(\theta )
	= 
	\lim_{n\rightarrow \infty } 
	E\left(
	\frac{1}{n} 
	\log p_{\theta }^{n}(Y_{1},\dots,Y_{n} ) 
	\right)
\end{align*}
for $\theta\in \Theta$, 
where 
\begin{align*}
	p_{\theta }^{n}(y_{1},\dots,y_{n} )
	=
	\sum_{x_{0},\dots,x_{n} \in {\cal X} }
	P_{\theta }(X_{0}^{\theta } = x_{0} )
	\prod_{i=1}^{n}
	\big(
	q_{\theta }(y_{k}|x_{k} ) p_{\theta }(x_{k}|x_{k-1} ) 
	\big) 
\end{align*}
for $\theta \in \Theta$, 
$y_{1}, \dots, y_{n} \in {\cal Y}$, $n\geq 0$.

In the statistics and engineering literature, 
$\{ (X_{n}, Y_{n} ) \}_{n\geq 0}$  
(as well as $\{ (X_{n}^{\theta }, Y_{n}^{\theta } ) \}_{n\geq 0}$)
is commonly referred to as a hidden Markov model with 
a finite state-space and continuous observations, 
while 
$X_{n}$ and $Y_{n}$ are considered as 
the (unobservable) state and (observable) output at discrete-time $n$. 
On the other hand, 
the identification of $\{ (X_{n}, Y_{n} ) \}_{n\geq 0}$
is regarded to as the estimation (or approximation) of  
$\{p(x'|x) \}_{x,x' \in {\cal X} }$ and
$\{Q(\cdot|x) \}_{x \in {\cal X}, y \in {\cal Y} }$ 
given the output sequence $\{Y_{n} \}_{n\geq 0}$. 
If the identification is based on the maximum likelihood principle
and the parameterized model 
$\{p_{\theta }(x'|x) \}_{x,x' \in {\cal X} }$, 
$\{q_{\theta }(y|x) \}_{x \in {\cal X}, y \in {\cal Y} }$, 
the estimation reduces to the maximization of 
the likelihood function $f(\cdot )$
over $\Theta$. 
In that context, 
$\{ (X_{n}^{\theta }, Y_{n}^{\theta } ) \}_{n\geq 0}$
is considered as a candidate model of 
$\{ (X_{n}, Y_{n} ) \}_{n\geq 0}$. 
For more details on hidden Markov models and their identification see 
\cite[Part II]{cappe&moulines&ryden} and references cited therein. 

Since the asymptotic mean of  
$\log p_{\theta }^{n}(Y_{1},\dots,Y_{n} ) / n$
is rarely available analytically, 
$f(\cdot )$ is usually maximized by a stochastic gradient algorithm, 
which itself is a special case of stochastic approximation
(for details see \cite{benveniste&metivier&priouret}, \cite{kushner&yin}, 
\cite{polyak} and references cited therein). 
To define such an algorithm, we introduce some further notation. 
For $\theta\in \mathbb{R}^{d_{\theta } }$, $x,x' \in {\cal X}$, $y\in {\cal Y}$, 
let 
\begin{align*}
	r_{\theta }(y|x',x)
	=
	q_{\theta }(y|x') p_{\theta }(x'|x),  
\end{align*}
while 
$R_{\theta }(y)$ is an $\mathbb{R}^{N_{x} \times N_{x} }$ matrix whose
$(i,j)$ entry is 
$r_{\theta }(y|i,j)$
(i.e., $R_{\theta }(y) = [r_{\theta }(y|i,j) ]_{i,j \in {\cal X} }$). 
%Moreover, 
%for $\theta \in \mathbb{R}^{d_{\theta } }$, 
%$1\leq k \leq d_{\theta }$, $y\in {\cal Y}$, 
%$S_{\theta }^{k}(y)$ is the partial derivative of 
%$\mathbb{R}_{\theta }(y)$ with respect to the $l$-th component of $\theta$. 
On the other side, 
for $\theta \in \mathbb{R}^{d_{\theta } }$, 
$u \in [0,\infty )^{N_{x} } \setminus \{0\}$, 
$V \in \mathbb{R}^{d_{\theta } \times N_{x} }$, $y \in {\cal Y}$, 
$1\leq k \leq d_{\theta }$, 
let 
\begin{align*}
	&
	\phi_{\theta }(u,y) 
	=
	\log(e^{T} R_{\theta }(y) u ), 
	\\
	&
	F_{\theta }(u,V,y)
	=
	\nabla_{\theta } \phi_{\theta }(u,y) 
	+
	V \nabla_{u} \phi_{\theta }(u,y), 
	\\
	&
	G_{\theta }(u,y)
	=
	\frac{R_{\theta }(y) u}{e^{T} R_{\theta }(y) u}, 
	\\
	&
	H_{\theta }(u,V,y)
	=
	\nabla_{\theta } G_{\theta }(u,y) 
	+
	V \nabla_{u} G_{\theta }(u,y)
\end{align*}
where 
$e=[1\dots 1]^{T} \in \mathbb{R}^{N_{x} }$. 
With this notation, 
a stochastic gradient algorithm  
for maximizing $f(\cdot )$
can be defined as 
\begin{align}
	& \label{1.1}
	\theta_{n+1}
	=
	\theta_{n} 
	+
	\alpha_{n} F_{\theta_{n} }(U_{n}, V_{n}, Y_{n+1} ), 
	\\ \label{1.3}
	&
	U_{n+1}
	=
	G_{\theta_{n+1} }(U_{n}, Y_{n+1} ), 
	\\ \label{1.5}
	&
	V_{n+1}
	=
	H_{\theta_{n+1} }(U_{n}, V_{n}, Y_{n+1} ), 
	\;\;\; n\geq 0. 
\end{align}
In this recursion, 
$\{\alpha_{n} \}_{n\geq 0}$ denotes a sequence of positive reals. 
$\theta_{0} \in \mathbb{R}^{d_{\theta } }$, 
$U_{0} \in \mathbb{R}^{N_{x} }$ and 
$V_{0} \in \mathbb{R}^{d_{\theta } \times N_{x} }$ 
are random variables which are defined on the probability space
$(\Omega, {\cal F}, P )$
and independent of $\{Y_{n} \}_{n\geq 0}$. 

In the literature on hidden Markov models and system identification, 
recursion (\ref{1.1}) -- (\ref{1.5}) is known as 
the recursive maximum likelihood algorithm, 
while subrecursions (\ref{1.3}) and (\ref{1.5}) are 
referred to as the optimal filter
and 
the optimal filter derivatives, respectively
(see \cite{cappe&moulines&ryden} for further details). 
Recursion (\ref{1.1}) -- (\ref{1.5})
usually includes a projection (or truncation) device 
which prevents estimates $\{\theta_{n} \}_{n\geq 0}$
from leaving $\Theta$ (see \cite{ljung} for further details). 
However, in order to avoid unnecessary technical details 
and to keep the exposition as simple as possible, 
this aspect of algorithm (\ref{1.1}) -- (\ref{1.5})
is not considered here. 
Instead, similarly as in 
\cite{ljung}, 
our results on the asymptotic behavior of 
algorithm (\ref{1.1}) -- (\ref{1.5}) 
(Theorems \ref{theorem2} and \ref{theorem3}) are expressed in a local form. 

Throughout the paper, unless stated otherwise, the following notation is used. 
For an integer $d\geq 1$, 
${\cal P}^{d}$ denotes the set of 
$d$-dimensional probability vectors 
(i.e., 
${\cal P}^{d} = 
\{u \in [0,\infty )^{d}: e^{T} u = 1 \}$), 
while $\mathbb{C}^{d}$ and $\mathbb{C}^{d\times d}$
are the sets of 
$d$-dimensional complex vectors 
and $d\times d$ complex matrices (respectively). 
$\|\cdot \|$ is the Euclidean norm in $\mathbb{R}^{d}$ or $\mathbb{C}^{d}$, 
while $d(\cdot, \cdot )$ is the distance induced by this norm. 
For a real number $\delta \in (0,\infty )$ and a set $A \subseteq \mathbb{C}^{d}$, 
$V_{\delta }(A)$ is the (complex) $\delta$-vicinity of $A$
induced by distance $d(\cdot, \cdot )$, i.e., 
\begin{align*}
	V_{\delta }(A) 
	= 
	\{w \in \mathbb{C}^{d}: d(w,A) \leq \delta \}. 
\end{align*}
$S$ is the set of stationary points of $f(\cdot )$, i.e., 
\begin{align*}
	S = \{\theta \in \Theta: \nabla f(\theta ) = 0 \}. 
\end{align*}

Algorithm (\ref{1.1}) -- (\ref{1.5}) is analyzed under the following assumptions. 

\begin{assumption} \label{a1}
$\lim_{n\rightarrow \infty } \alpha_{n} = 0$, 
$\limsup_{n\rightarrow \infty } |\alpha_{n+1}^{-1} - \alpha_{n}^{-1} | < \infty$ 
and 
$\sum_{n=0}^{\infty } \alpha_{n} = \infty$. 
Moreover, 
there exists a real number $r\in (1,\infty )$
such that 
$\sum_{n=0}^{\infty } \alpha_{n}^{2} \gamma_{n}^{2r} < \infty$. 
\end{assumption} 

\begin{assumption} \label{a2}
$\{X_{n} \}_{n\geq 0}$ is geometrically ergodic. 
\end{assumption}

\begin{assumption} \label{a3}
There exists a function 
$s_{\theta }(y|x)$ mapping $(\theta, x, y ) \in \Theta \times {\cal X} \times {\cal Y}$ 
into $[0, \infty )$, 
and for any compact set $Q \subset \Theta$,  
there exists a real number 
$\varepsilon_{Q} \in (0,1)$ 
such that 
\begin{align*}
	\varepsilon_{Q} s_{\theta }(y|x') 
	\leq 
	r_{\theta }(y|x',x) 
	\leq 
	\varepsilon_{Q}^{-1} s_{\theta }(y|x')  
\end{align*}
for all $\theta \in Q$, $x,x' \in {\cal X}$, $y\in {\cal Y}$. 
\end{assumption}

\begin{assumption} \label{a4}
For each $y \in {\cal Y}$, 
$\phi_{\theta }(u,y)$ and 
$G_{\theta }(u,y)$ are real-analytic functions of $(\theta,u)$
on entire $\Theta\times {\cal P}^{N_{x} }$. 
Moreover, $\phi_{\theta }(u,y)$ and $G_{\theta }(u,y)$
have (complex-valued) analytical continuations 
$\hat{\phi}_{\eta}(w,y)$ and $\hat{G}_{\eta }(w,y)$ (respectively)
with the following properties: 
\begin{enumerate}
\item
$\hat{\phi}_{\eta}(w,y)$ and
$\hat{G}_{\eta }(w,y)$ map
$(\eta, w, y ) \in \mathbb{C}^{d_{\theta } } \times \mathbb{C}^{N_{x} } \times {\cal Y}$
into $\mathbb{C}$ and $\mathbb{C}^{N_{x} }$ (respectively). 
\item
$\hat{\phi}_{\theta }(u,y) = \phi_{\theta }(u,y)$ and
$\hat{G}_{\theta }(u,y) = G_{\theta }(u,y)$
for all $\theta \in \Theta$, $u \in {\cal P}^{N_{x} }$, $y \in {\cal Y}$. 
\item
For any compact set $Q\subset \Theta$, 
there exist real numbers 
$\delta_{Q} \in (0,1)$, $K_{Q} \in [1, \infty )$ 
and a Borel-measurable function 
$\psi_{Q}: {\cal Y} \rightarrow [1, \infty )$ 
such that 
$\hat{\phi}_{\eta }(w,y)$ and 
$\hat{G}_{\eta }(w,y)$ are analytical in $(\eta,w)$ on 
$V_{\delta_{Q} }(Q) \times V_{\delta_{Q} }({\cal P}^{N_{x} } )$ 
for each $y \in {\cal Y}$, 
and such that 
\begin{align*}
	&
	|\hat{\phi}_{\eta}(w,y) | \leq \psi_{Q}(y), 
	\\
	&
	\|\hat{G}_{\eta}(w,y) \|
	\leq 
	K_{Q}, 
	\\
	&
	\int \psi_{Q}^{2}(y') Q(dy'|x) < \infty
\end{align*}
for all
$\eta \in V_{\delta_{Q} }(Q)$, 
$w \in V_{\delta_{Q} }({\cal P}^{N_{x} } )$, 
$x \in {\cal X}$, $y \in {\cal Y}$. 
\end{enumerate}
\end{assumption}

Assumption \ref{a1} corresponds to the properties 
of step-size sequence 
$\{\alpha_{n} \}_{n\geq 0}$
and is commonly used in the asymptotic analysis of 
stochastic approximation algorithms. 
It holds if 
$\alpha_{n} = 1/n^{a}$ for $n\geq 1$, where $a\in (3/4,1]$. 

Assumptions \ref{a2} and \ref{a3} 
are related to the stability of the model 
$\{ (X_{n}, Y_{n} ) \}_{n\geq 0}$ and its optimal filter. 
In this or similar form, 
they are involved in the analysis of various 
aspects of optimal filtering 
and parameter estimation in hidden Markov models
(see e.g., 
\cite{bickel&ritov}, 
\cite{bickel&ritov&ryden}, 
\cite{douc&matias}, \cite{douc&moulines&ryden}, 
\cite{legland&mevel1} -- \cite{legland&mevel4}, 
\cite{leroux}, 
\cite{mevel&finesso}, 
\cite{ryden1}, \cite{ryden2}, 
\cite{tadic&doucet}; 
see also \cite[Part II]{cappe&moulines&ryden} and references 
cited therein). 

Assumption \ref{a4} corresponds to the parametrization 
of candidate models 
$\{ (X_{n}^{\theta }, Y_{n}^{\theta } ) \}_{n\geq 0}$. 
Basically, Assumption \ref{a4} requires 
transition probabilities 
$p_{\theta }(x'|x)$
and observation conditional densities $q_{\theta }(y|x)$
to be analytic in $\theta$. 
It also requires 
$p_{\theta }(x'|x)$
and $q_{\theta }(y|x)$
can be analytically continuable to a complex domain
such that the corresponding continuation of 
the optimal filter transfer function 
$G_{\theta}(u,y)$
is analytic and uniformly bounded in 
$(\theta, u)$. 
Although these requirements are restrictive, 
they still hold in many practically relevant cases and situations. 
Several examples are provided in the next section. 

In order to state our main results 
we rely on the following notation. 
$\{\gamma_{n} \}_{n\geq 0}$ is a sequence of real numbers defined by 
$\gamma_{0}=1$
and 
\begin{align*}
	\gamma_{n} = 1+\sum_{i=0}^{n-1} \alpha_{i}
\end{align*}
for $n \geq 1$. 
Event $\Lambda$
is defined as 
\begin{align*}
	\Lambda 
	=
	\left\{
	\sup_{n\geq 0} \|\theta_{n} \| < \infty, 
	\inf_{n\geq 0} d(\theta_{n}, \partial \Theta ) > 0
	\right\}.  
\end{align*}
With this notation, our main results on the properties of 
objective function $f(\cdot )$
and algorithm 
(\ref{1.1}) -- (\ref{1.5}) can be stated as follows. 

\begin{theorem}[Analyticity] \label{theorem1}
Let Assumptions \ref{a2} -- \ref{a4} hold.
Then, the following is true: 
\begin{enumerate}
\item
$f(\cdot )$ is analytic on entire $\Theta$. 
\item
For each $\theta \in \Theta$, there exist 
real numbers 
$\delta_{\theta } \in (0,1)$, $\mu_{\theta } \in (1,2]$, 
$M_{\theta } \in [1, \infty )$
such that 
\begin{align*}
	|f(\theta' ) - f(\theta ) |
	\leq 
	M_{\theta } \|\nabla f(\theta' ) \|^{\mu_{\theta } }
\end{align*}
for all $\theta' \in \Theta$ satisfying 
$\|\theta - \theta' \|\leq \delta_{\theta }$. 
\end{enumerate}
\end{theorem}

\begin{theorem}[Convergence] \label{theorem2}
Let Assumption \ref{a1} -- \ref{a4} hold. 
Then, 
$\hat{\theta } = \lim_{n\rightarrow \infty } \theta_{n}$
exists and satisfies 
$\nabla f(\hat{\theta} ) = 0$
w.p.1 on 
event $\Lambda$. 
\end{theorem}

\begin{theorem}[Convergence Rate] \label{theorem3}
Let Assumptions \ref{a1} -- \ref{a4} hold. 
Then, 
\begin{align}\label{t3.1*}
	\|\nabla f(\theta_{n} ) \|^{2}
	=
	O\big(\gamma_{n}^{-\hat{p} } \big),
	\;\;\;\;\; 
	|f(\theta_{n} ) - f(\hat{\theta } ) |
	=
	O\big(\gamma_{n}^{-\hat{p} } \big), 
	\;\;\;\;\; 
	\|\theta_{n} - \hat{\theta } \| 
	=
	O\big(\gamma_{n}^{-\hat{q} } \big)
\end{align}
w.p.1 on $\Lambda$, 
where $\hat{\mu} = \mu_{\hat{\theta } }$ and 
\begin{align}
	&\label{t3.3*}
	\hat{r} 
	= 
	\begin{cases}
	1/(2 - \hat{\mu } ), 
	&\text{ if } \hat{\mu } < 2
	\\
	\infty, 
	&\text{ otherwise}
	\end{cases}, 
	\;\;\;\;\;
	\hat{p} 
	=
	\hat{\mu } \min\{r, \hat{r} \}, 
	\;\;\;\;\; 
	\hat{q}
	=
	\min\{(\hat{p} - 1 )/2,r-1 \}.   
\end{align}  
\end{theorem}

Proofs of the Theorems \ref{theorem1} -- \ref{theorem3} 
are provided in Section \ref{section1*}. 

In the literature on deterministic and stochastic optimization
(notice that recursion (\ref{1.1}) -- (\ref{1.5}) belongs 
to the class of stochastic gradient algorithms), 
the convergence of gradient search is usually characterized 
by gradient, objective and estimate convergence, 
i.e., by the convergence of sequences
$\{\nabla f(\theta_{n} ) \}_{n\geq 0}$, 
$\{f(\theta_{n} ) \}_{n\geq 0}$
and 
$\{\theta_{n} \}_{n\geq 0}$
(see e.g., \cite{bertsekas}, \cite{bertsekas&tsitsiklis}, 
\cite{polyak&tsypkin}, \cite{polyak} and references cited therein). 
Similarly, the convergence rate can be described by the rates 
at which sequences 
$\{\nabla f(\theta_{n} ) \}_{n\geq 0}$, 
$\{f(\theta_{n} ) \}_{n\geq 0}$
and 
$\{\theta_{n} \}_{n\geq 0}$
converge to the sets of their accumulation points. 
In the case of algorithm (\ref{1.1}) -- (\ref{1.5}), 
this kind of information is provided by 
Theorems \ref{theorem2} and \ref{theorem3}. 
Basically, Theorem \ref{theorem2} claims that 
recursion (\ref{1.1}) -- (\ref{1.5}) is point-convergent w.p.1
(i.e.,
the set of accumulation points of 
$\{\theta_{n} \}_{n\geq 0}$
is almost surely a singleton), 
while Theorem \ref{theorem3} provides relatively tight bounds on 
convergence rate in the terms of Lojasiewicz exponent $\mu_{\hat{\theta} }$ 
and the convergence rate of step-sizes 
$\{\alpha_{n} \}_{n\geq 0}$
(expressed through $r$ and $\{\gamma_{n} \}_{n\geq 0}$). 
Theorem \ref{theorem1}, on the other side, deals with the properties of 
the asymptotic log-likelihood $f(\cdot )$
and is a crucial prerequisite for 
Theorems \ref{theorem2} and \ref{theorem3}. 
Apparently, the results of Theorems \ref{theorem2} and \ref{theorem3} 
are of local nature: 
They hold on the event where 
algorithm (\ref{1.1}) -- (\ref{1.5}) is stable
(i.e., where 
$\{\theta_{n} \}_{n\geq 0}$ is contained in a compact subset of 
$\Theta$). 
Stating asymptotic results in such a form is quite common   
for stochastic recursive algorithms 
(see e.g., \cite{kushner&yin}, \cite{ljung} and references cited therein).   
Moreover, a global version of Theorems 
\ref{theorem2} and \ref{theorem3} 
can be obtained easily by 
combining them with methods used to verify or ensure stability 
(e.g., with \cite{borkar&meyn}, \cite{chen} or \cite{ljung}). 

Various asymptotic properties of maximum likelihood estimation in 
hidden Markov models have been analyzed thoroughly in 
a number of papers 
%\cite{araposthatis&marcus}, 
\cite{baum&petrie}, 
\cite{bickel&ritov}, 
\cite{bickel&ritov&ryden}, 
\cite{douc&matias}, \cite{douc&moulines&ryden}, 
\cite{legland&mevel1} -- \cite{legland&mevel4}, 
\cite{leroux}, 
\cite{mevel&finesso}, 
\cite{ryden1}, \cite{ryden2}; 
(see also \cite[Chapter 12]{cappe&moulines&ryden}, \cite{ephraim&merhav} and references 
cited therein). 
Although these results offer a deep insight into 
the asymptotic behavior of this estimation method, 
they can hardly be applied to complex hidden Markov models. 
The reason comes out of the fact that 
all existing results on the point-convergence and convergence 
rate of stochastic gradient search 
(which includes recursive maximum likelihood estimation as a special case)
require objective function to have an isolated maximum (or minimum)
at which the Hessian is strictly negative definite. 
Since $f(\cdot )$, the objective function of recursion (\ref{1.1}) -- (\ref{1.5}), 
is rather complex even when the observation space is finite 
(i.e., ${\cal Y} = \{1,\dots,N_{y} \}$) and 
$N_{x}$, $N_{y}$, the numbers of states and observations, are relatively small 
(three and above), 
it is hard (if possible at all) to show the existence of isolated maxima, 
let alone checking the definiteness of $\nabla^{2} f(\cdot )$. 
Exploiting the analyticity of $f(\cdot )$ and Lojasiewicz inequality, 
Theorems \ref{theorem2} and \ref{theorem3} overcome these difficulties: 
They both neither require the existence of an isolated maximum, 
nor impose any restriction on the definiteness of the Hessian
(notice that the Hessian cannot be strictly definite at 
a non-isolated maximum or minimum). 
In addition to this, the theorems cover a relatively broad class 
of hidden Markov models 
(see the next section). 
To the best of our knowledge, asymptotic results with similar
features do not exist in the literature on hidden Markov models
or stochastic optimization. 

The differentiability, analyticity and other analytic properties of 
the entropy rate of hidden Markov models, a functional similar to 
the asymptotic likelihood, have been studied thoroughly in several papers 
\cite{han&marcus1}, 
\cite{han&marcus2}, 
\cite{holliday&goldsmith&glynn}, 
\cite{ordentlitch&weissman}, 
\cite{peres}, 
\cite{schoenhuht}. 
The results presented therein 
cover only models with discrete state and observation spaces
and do not pay any attention to maximum likelihood estimation. 
Motivated by the problem of the point-convergence and convergence rate 
of recursive maximum likelihood estimators for hidden Markov models, 
we extend these results in Theorem \ref{theorem1} 
to models with continuous observations
and their likelihood functions. 
The approach we use to demonstrate the analyticity of the asymptotic likelihood 
is based on the principle of analytical continuation 
and is similar to the methodology formulated 
in \cite{han&marcus1}. 

\section{Examples} \label{section2}

In this section, we consider several practically relevant examples of 
the results presented in Section \ref{section1}. 
Analyzing these examples, we also provide a direction how the assumptions 
adopted in Section \ref{section1} can be verified in practice. 

\subsection{Finite Observation Space}

Hidden Markov models with finite state and observation spaces are studied in 
this subsection. 
For these models, we show that the conclusion of Theorems \ref{theorem1} -- \ref{theorem3} 
hold whenever the parameterization of candidate models is analytic. 

Let $N_{y}>2$ be an integer, while 
${\cal Y} = \{1,\dots,N_{y} \}$. 
Then, the following results hold. 

\begin{proposition} \label{proposition1}
Assumptions \ref{a3} and \ref{a4} are true if the following conditions are 
satisfied: 
\begin{enumerate}
\item
For each $x,x'\in {\cal X}$, $y\in {\cal Y}$, 
$r_{\theta }(y|x',x)$ is analytical in $\theta$ on entire $\Theta$.
\item
$r_{\theta }(y|x',x) > 0$ for all $\theta \in \Theta$, 
$x,x' \in {\cal X}$, $y \in {\cal Y}$. 
\end{enumerate}
\end{proposition}

\begin{corollary}
Let Assumptions \ref{a1}, \ref{a2} and the conditions of Proposition \ref{proposition1} 
hold. Then, the conclusions of Theorems \ref{theorem1} -- \ref{theorem3} are true. 
\end{corollary}

The proof is provided in Section \ref{section2*}. 

\begin{remark}
The conditions of Proposition \ref{proposition1} correspond to the way the candidate models
are parameterized. 
They hold for the natural\footnote{
The natural parameterization can be defined as follows:  
$\theta = 
[\alpha_{1,1} \cdots \alpha_{N_{x},N_{x} } \: 
\beta_{1,1} \cdots \beta_{N_{x},N_{y} } ]^{T}$
and $p_{\theta}(x'|x)=\alpha_{x,x'}$, 
$q_{\theta}(y|x)=\beta_{x,y}$ 
for $x,x' \in {\cal X}$, $y \in {\cal Y}$, 
while $\Theta$ is the set of vectors 
$[\alpha_{1,1} \cdots \alpha_{N_{x},N_{x} } \: 
\beta_{1,1} \cdots \beta_{N_{x},N_{y} } ]^{T} \in (0,1)^{N_{x}(N_{x} + N_{y} ) }$
satisfying 
$\sum_{l=1}^{N_{x} } \alpha_{x,l} 
=
\sum_{l=1}^{N_{y} } \beta_{x,l} = 1$
for each  $x \in {\cal X}$.}, 
exponential\footnote{
In the case of the exponential parameterization, 
we have 
$\theta = 
[\alpha_{1,1} \cdots \alpha_{N_{x},N_{x} } \: 
\beta_{1,1} \cdots \beta_{N_{x},N_{y} } ]^{T}$, 
and 
\begin{align*}
	p_{\theta }(x'|x)
	=
	\frac{\exp(\alpha_{x,x'} ) }{\sum_{l=1}^{N_{x} } \exp(\alpha_{x,l} ) }, 
	\;\;\; 
	q_{\theta }(y|x)
	=
	\frac{\exp(\beta_{x,y} ) }{\sum_{l=1}^{N_{y} } \exp(\beta_{x,l} ) }
\end{align*}
for $x,x' \in {\cal X}$, $y \in {\cal Y}$, 
while $\Theta = \mathbb{R}^{N_{x}(N_{x}+N_{y} ) }$.} 
and 
trigonometric\footnote{
The trigonometric parameterization is defined as 
$\theta = [\alpha_{1,1} \cdots \alpha_{N_{x},N_{x} } \: 
\beta_{1,1} \cdots \beta_{N_{x},N_{y} } ]^{T}$ and 
\begin{align*}
	&
	p_{\theta}(1|x)
	=
	\cos^{2} \alpha_{x,1}, 
	\;\;\;\;\; 
%	&
	p_{\theta}(x'|x)
	=
	\cos^{2} \alpha_{x,x'} \prod_{l=1}^{x'-1} \sin^{2} \alpha_{x,l}, 
	\;\;\;\;\; 
%	&
	p_{\theta}(N_{x}|x)
	=
	\prod_{l=1}^{N_{x} } \sin^{2} \alpha_{x,l}, 
	\\
	&
	q_{\theta}(1|x)
	=
	\cos^{2} \beta_{x,1}, 
	\;\;\;\;\; 
%	&
	q_{\theta}(y|x)
	=
	\cos^{2} \beta_{x,y} \prod_{l=1}^{y-1} \sin^{2} \beta_{x,l}, 
	\;\;\;\;\; 
%	&
	q_{\theta}(N_{y}|x)
	=
	\prod_{l=1}^{N_{y} } \sin^{2} \beta_{x,l}
\end{align*}
for $x \in {\cal X}$, $x' \in {\cal X}\setminus \{1,N_{x} \}$, 
$y \in {\cal Y}\setminus \{1,N_{y} \}$, 
while $\Theta = (0,\pi/2)^{N_{x} (N_{x} + N_{y} ) }$. } 
parameterizations. 
\end{remark}

\subsection{Compactly Supported Observations}

In this subsection, we consider hidden Markov models with 
a finite number of states and compactly supported observations. 
More specifically, we assume that ${\cal Y}$ is a compact set from $\mathbb{R}^{d_{y} }$. 
For such models, the following results can be shown. 

\begin{proposition} \label{proposition2}
Assumptions \ref{a3} and \ref{a4} are true if the following conditions are satisfied:
\begin{enumerate}
\item
For each $x,x'\in {\cal X}$, 
$r_{\theta }(y|x',x)$ is analytical in $(\theta,y)$ on entire $\Theta\times {\cal Y}$.
\item
$r_{\theta }(y|x',x) > 0$ for all $\theta \in \Theta$, 
$x,x' \in {\cal X}$, $y \in {\cal Y}$. 
\end{enumerate} 
\end{proposition}

\begin{corollary}
Let Assumptions \ref{a1}, \ref{a2} and the conditions of Proposition \ref{proposition2} 
hold. Then, the conclusions of Theorems \ref{theorem1} -- \ref{theorem3} are true. 
\end{corollary}

The proof is provided in Section \ref{section2*}. 

\begin{remark}
The conditions of Proposition \ref{proposition2} are fulfilled 
if the natural, exponential or trigonometric parameterization 
(see the previous subsection) is applied to the state transition probabilities 
$\{p_{\theta}(x'|x) \}_{x,x'\in {\cal X} }$, 
and 
if the observation likelihoods
$\{q_{\theta }(\cdot|x) \}_{x\in {\cal X} }$ are analytic jointly in 
$\theta$ and $y$. 
The later holds when  
$\{q_{\theta }(\cdot |x) \}_{x\in {\cal X} }$ 
are compactly truncated mixtures of 
beta, exponential, gamma, logistic, normal, log-normal, 
Pareto, uniform, Weibull distributions, and 
when 
each of these mixtures is indexed by
its weights 
and by 
the `natural' parameters of its ingredient distributions. 
\end{remark}

\subsection{Mixture of Observation Likelihoods} 

In this subsection, we consider the case when 
the observation likelihoods 
$\{q_{\theta}(\cdot|x) \}_{x\in {\cal X} }$
are mixtures of known probability density functions. 
More specifically, let $d_{\alpha } \geq 1$, $N_{\beta } > 1$ be integers, 
while ${\cal A}\subseteq \mathbb{R}^{d_{\alpha } }$ is an open set and  
\begin{align*}
	{\cal B} 
	=
	\left\{
	[\beta_{1,1} \cdots \beta_{N_{x},N_{\beta } } ]^{T} \in (0, 1 )^{N_{x} N_{\beta } }: 
	\sum_{i=1}^{N_{\beta } } \beta_{x,k} = 1
	\text{ for each } x\in {\cal X}
	\right\}. 
\end{align*}
We assume that 
the state transition probabilities and the observation likelihoods 
are parameterized by vectors 
$\alpha \in {\cal A}$ and $\beta \in {\cal B}$ (respectively), i.e., 
$p_{\theta }(x'|x) = p_{\alpha }(x'|x)$, 
$q_{\theta }(y|x) = q_{\beta }(y|x)$ 
for $\alpha \in {\cal A}$, $\beta \in {\cal B}$,
$\theta = [\alpha^{T} \: \beta^{T} ]^{T}$, 
$x,x' \in {\cal X}$, $y \in {\cal Y}$. 
We also assume
\begin{align*}
	q_{\beta }(y|x)
	=
	\sum_{k=1}^{N_{\beta } } \beta_{x,k} f_{k}(y|x),
\end{align*}
where 
$\beta = [\beta_{1,1} \cdots \beta_{N_{x},N_{\beta } } ]^{T} \in {\cal B}$, 
$x \in {\cal X}$, $y \in {\cal Y}$, 
while 
$\{f_{k}(\cdot |x)\}_{x \in {\cal X}, 1 \leq k \leq N_{\beta } }$
are known probability density functions. 

For the models specified in this subsection, the following results hold. 

\begin{proposition} \label{proposition3}
Assumptions \ref{a3} and \ref{a4} are true if the following conditions are satisfied: 
\begin{enumerate}
\item
For each $x,x'\in {\cal X}$, 
$p_{\alpha }(x'|x)$ is analytical in $\alpha$ on entire ${\cal A}$.
\item
$p_{\alpha }(x'|x) > 0$ for all $\alpha \in {\cal A}$, 
$x,x' \in {\cal X}$. 
\item
$\psi(y) > 0$ and 
$
	\int \log^{2} \psi(y') Q(dy'|x ) < \infty
$
for all $x\in {\cal X}$, $y \in {\cal Y}$, 
where 
$\psi(y) = 
\sum_{x \in {\cal X} }
\sum_{k=1}^{N_{\beta } }
f_{k}(y|x)$. 
\end{enumerate}
\end{proposition}

\begin{corollary}
Let Assumptions \ref{a1}, \ref{a2} and the conditions of Proposition \ref{proposition3} 
hold. Then, the conclusions of Theorems \ref{theorem1} -- \ref{theorem3} are true. 
\end{corollary}

The proof is provided in Section \ref{section2*}. 

\subsection{Gaussian Observations}

This subsection is devoted to hidden Markov models with a finite number of states 
and with Gaussian observations. 
More specifically, let $d_{\alpha }$ and ${\cal A}$ have the same meaning as in 
the previous section, while ${\cal Y} = \mathbb{R}$ and 
\begin{align} \label{2.1}
	{\cal B}
	=
	\left\{
	[\lambda_{1} \cdots \lambda_{N_{x} } \: \mu_{1} \cdots \mu_{N_{x} } ]^{T} 
	\in (0,\infty )^{N_{x} } \times \mathbb{R}^{N_{x} }
	: \lambda_{x} \neq \lambda_{x} \text{ for } x\neq x', x,x' \in {\cal X} 
	\right\}.
\end{align}
Similarly as in the previous subsection, 
we assume that the state transition probabilities and 
the observation likelihoods are indexed by 
vectors $\alpha \in {\cal A}$ and $\beta \in {\cal B}$ (respectively). 
We also assume 
\begin{align*}
	q_{\beta }(y|x) 
	=
	\sqrt{\lambda_{x}/\pi }
	\exp(-\lambda_{x} (y-\mu_{x} )^{2} ), 
\end{align*}
where 
$\beta = 
[\lambda_{1} \cdots \lambda_{N_{x} } \: \mu_{1} \cdots \mu_{N_{x} }]^{T} \in {\cal B}$, 
$x\in {\cal X}$, $y \in {\cal Y}$. 

For the models described in this subsection, the following results can 
be shown. 

\begin{proposition} \label{proposition4}
Assumptions \ref{a3} and \ref{a4} are true if the following conditions are satisfied: 
\begin{enumerate}
\item
For each $x,x'\in {\cal X}$, 
$p_{\alpha }(x'|x)$ is analytical in $\alpha$ on entire ${\cal A}$.
\item
$p_{\alpha }(x'|x) > 0$ for all $\alpha \in {\cal A}$, 
$x,x' \in {\cal X}$. 
\item
$
	\int y^{4} Q(dy|x ) < \infty
$
for all $x\in {\cal X}$. 
\end{enumerate}
\end{proposition}

\begin{corollary} \label{corollary4}
Let Assumptions \ref{a1}, \ref{a2} and the conditions of Proposition \ref{proposition4} 
hold. Then, the conclusions of Theorems \ref{theorem1} -- \ref{theorem3} are true. 
\end{corollary}

The proof is provided in Section \ref{section2*}. 

\begin{remark}
Unfortunately, Proposition \ref{proposition4} and Corollary \ref{corollary4} cannot 
be extended to the case 
${\cal B} = (0,\infty )^{N_{x} } \times \mathbb{R}^{N_{x} }$, 
since the models specified in the subsection do not satisfy 
Assumption \ref{a4} 
without the condition 
$\lambda_{x} \neq \lambda_{x'}$ for $x\neq x'$
(which appears in (\ref{2.1})).\footnote
{Let 
$h_{\alpha, y, u }(\beta ) = e^{T} R_{\theta }(y) u$ for 
$\alpha \in {\cal A}$, $\beta \in {\cal B}$, 
$\theta = [\alpha^{T} \beta^{T} ]^{T}$, 
$y\in {\cal Y}$, $u \in {\cal P}^{N_{x} }$. 
Obviously, for any 
$\alpha \in {\cal A}$, 
$y\in {\cal Y}$, $u \in {\cal P}^{N_{x} }$, 
$h_{\alpha, y, u }(\cdot ) $ has a unique (complex-valued) 
analytical continuation, which can be defined as
\begin{align*}
	\hat{h}_{\alpha, y, u }(b) 
	=
	\sum_{x,x' \in {\cal X} }
	\sqrt{l_{x'}/\pi }
	\exp(-l_{x'} (y - m_{x'} )^{2} )
	p_{\alpha }(x'|x) u_{x}
\end{align*}
where 
$b = [l_{1} \cdots l_{N_{x} } \: m_{1} \cdots m_{N_{x} } ]^{T} \in \mathbb{C}^{2N_{x} }$. 
Let 
$\beta = 
[\lambda_{1} \cdots \lambda_{N_{x} } \: \mu_{1} \cdots \mu_{N_{x} } ]^{T} \in 
(0,\infty )^{N_{x} } \times \mathbb{R}^{N_{x} }$
be any vector 
satisfying 
$\lambda_{x} = \lambda_{x'}$
for some $x\neq x'$, $x,x' \in {\cal X}$. 
Then, it is not hard to deduce that 
there exist
$\alpha\in {\cal A}$, $y \in {\cal Y}$, $u \in {\cal P}^{N_{x} }$ 
(depending on $\beta$)
such that 
$\hat{h}_{\alpha,y,u}(\cdot )$ has a zero in any (complex) vicinity of 
$\beta$. 
Since the zeros of the analytical continuation of 
$e^{T} R_{\theta }(y) u$ would be the poles of the analytical continuation of 
$G_{\theta }(u,y)$, 
it is not possible to continue $G_{\theta }(u,y)$ analytically in 
any vicinity of point 
$(\theta, u )$, where  
$\theta = [\alpha^{T} \beta^{T} ]^{T}$. 
Hence, Proposition \ref{proposition4} and Corollary \ref{corollary4} 
cannot be extended to the case 
${\cal B} = (0, \infty )^{N_{x} } \times \mathbb{R}^{N_{x}}$. } 
However, this condition is not so restrictive in practice  
as ${\cal B}$ is dense
in $(0,\infty )^{N_{x} } \times \mathbb{R}^{N_{x} }$ 
and provides an arbitrarily close approximation to 
$(0,\infty )^{N_{x} } \times \mathbb{R}^{N_{x} }$. 
\end{remark}

\section{Proof of Main Results} \label{section1*}

\subsection{Optimal Filter and Its Properties}\label{subsection1.1*}

The stability properties (forgetting and ergodicity) of the optimal filter
(\ref{1.3}), its derivatives (\ref{1.5}) and its analytical continuation 
(to be defined in the next paragraph) are studied in 
this subsection. 
The analysis mainly follows the ideas and results of 
\cite{legland&mevel3}, \cite{legland&mevel4} and \cite{legland&oudjane}. 
The results presented in the subsection are 
an essential prerequisite for the analysis carried out 
in Subsections \ref{subsection1.2*} and \ref{subsection1.3*}. 

Throughout this subsection, we rely on the following notation. 
${\cal Q}^{N_{x} }$ denotes the set 
\begin{align*}
	{\cal Q}^{N_{x} }
	=
	\{u\in [0,\infty )^{N_{x} }: 
	e^{T} u \geq 1/2 \},  
\end{align*}
where $e = [1 \cdots 1]^{T} \in R^{N_{x} }$
(${\cal Q}^{N_{x} }$ can be any compact set from $[0,\infty )^{N_{x} }$
satisfying $0\not\in {\cal Q}^{N_{x} }$, 
$\text{int} {\cal P}^{N_{x} } \subset {\cal Q}^{N_{x} }$, 
but the above one is selected for analytical convenience). 
For $n\geq m\geq 0$ 
and a sequence
$\boldsymbol y = \{y_{n} \}_{n\geq 0}$
from ${\cal Y}$, 
$y_{m:n}$ denotes finite subsequence 
$(y_{m},\dots,y_{n})$. 
For $u \in [0,\infty )^{N_{x} } \setminus \{0\}$, 
$w \in \mathbb{C}^{N_{x} }$, 
$V \in \mathbb{R}^{d_{\theta } \times N_{x} }$, $n\geq m\geq 0$
and sequences 
$\boldsymbol \vartheta = \{\vartheta_{n} \}_{n\geq 0}$, 
$\boldsymbol \eta = \{\eta_{n} \}_{n\geq 0}$, 
$\boldsymbol y = \{y_{n} \}_{n\geq 1}$
from $\Theta$, $\mathbb{C}^{d_{\theta } }$, ${\cal Y}$ (respectively), 
let 
$G_{\boldsymbol \vartheta, \boldsymbol y}^{m:m}(u) = u$, 
$\hat{G}_{\boldsymbol \eta, \boldsymbol y}^{m:m}(w) = w$, 
$H_{\boldsymbol \vartheta, \boldsymbol y}^{m:m}(u,V) = V$
and 
\begin{align*}
	&
	G_{\boldsymbol \vartheta, \boldsymbol y}^{m:n+1}(u)
	=
	G_{\vartheta_{n+1} }(G_{\boldsymbol \vartheta, \boldsymbol y}^{m:n}(u),y_{n+1} ), 
	\\
	&
	\hat{G}_{\boldsymbol \eta, \boldsymbol y}^{m:n+1}(w)
	=
	\hat{G}_{\eta_{n+1} }(\hat{G}_{\boldsymbol \eta, \boldsymbol y}^{m:n}(w),y_{n+1} ), 
	\\
	&
	H_{\boldsymbol \vartheta, \boldsymbol y}^{m:n+1}(u,V)
	=
	H_{\vartheta_{n+1} }(G_{\boldsymbol \vartheta, \boldsymbol y}^{m:n}(u),
	H_{\boldsymbol \vartheta, \boldsymbol y}^{m:n}(u,V), 
	y_{n+1} )   
\end{align*}
($G_{\theta }(u,y)$, $\hat{G}_{\eta}(w,y)$, $H_{\theta}(u,V,y)$
are defined in Section \ref{section1}).
If $\boldsymbol \vartheta = \{\theta \}_{n\geq 0}$ 
(i.e., $\vartheta_{n} = \theta$), 
we also use notation 
$G_{\theta, \boldsymbol y}^{m:n}(u) = 
G_{\boldsymbol \vartheta, \boldsymbol y}^{m:n}(u)$, 
$H_{\theta, \boldsymbol y}^{m:n}(u,V) = 
H_{\boldsymbol \vartheta, \boldsymbol y}^{m:n}(u,V)$, 
as well as 
$G_{\theta}^{0:n}(u,y_{1:n} ) = G_{\boldsymbol \vartheta, \boldsymbol y}^{0:n}(u)$, 
$H_{\theta}^{0:n}(u,V,y_{1:n} ) = H_{\boldsymbol \vartheta, \boldsymbol y}^{0:n}(u,V)$. 
Similarly, 
if $\boldsymbol \eta = \{\eta \}_{n\geq 0}$ 
(i.e., $\eta_{n} = \eta$), 
we rely on notation 
$\hat{G}_{\eta, \boldsymbol y}^{m:n}(w) = 
\hat{G}_{\boldsymbol \eta, \boldsymbol y}^{m:n}(w)$
and
$\hat{G}_{\eta}^{0:n}(w,y_{1:n} ) = \hat{G}_{\boldsymbol \eta, \boldsymbol y}^{0:n}(w)$. 
Then, it straightforward to verify 
\begin{align*}
	&%\label{1*.1.1}
	G_{\boldsymbol \vartheta, \boldsymbol y }^{m:n}(u)
	=
	G_{\boldsymbol \vartheta, \boldsymbol y }^{k:n}
	(G_{\boldsymbol \vartheta, \boldsymbol y }^{m:k}(u)), 
	\\
	&%\label{1*.1.3}
	\hat{G}_{\boldsymbol \eta, \boldsymbol y }^{m:n}(w)
	=
	\hat{G}_{\boldsymbol \eta, \boldsymbol y }^{k:n}
	(\hat{G}_{\boldsymbol \eta, \boldsymbol y }^{m:k}(w)), 
	\\
	&%\label{1*.1.5}
	H_{\boldsymbol \vartheta, \boldsymbol y }^{m:n}(u,V)
	=
	H_{\boldsymbol \vartheta, \boldsymbol y }^{k:n}
	(G_{\boldsymbol \vartheta, \boldsymbol y }^{m:k}(u), 
	H_{\boldsymbol \vartheta, \boldsymbol y }^{m:k}(u,V) )  
\end{align*}
for each  
$u \in [0,\infty )^{N_{x} } \setminus \{0 \}$, 
$w \in \mathbb{C}^{N_{x} }$, 
$V \in \mathbb{R}^{d_{\theta } \times N_{x} }$, 
$0 \leq m \leq k \leq n$
and any sequences 
$\boldsymbol \vartheta = \{\vartheta_{n} \}_{n\geq 0}$, 
$\boldsymbol \eta = \{\eta_{n} \}_{n\geq 0}$, 
$\boldsymbol y = \{y_{n} \}_{n\geq 1}$
from $\Theta$, $\mathbb{C}^{d_{\theta } }$, ${\cal Y}$
(respectively). 
Moreover, it can be demonstrated easily 
\begin{align}\label{1*.1.7}
	\hat{G}_{\eta}^{0:n}(w,y_{1:n} ) 
	=
	\hat{G}_{\eta}^{0:k}(\hat{G}_{\eta}^{0:n-k}(w,y_{1:n-k}), y_{n-k+1:n} )
\end{align}
for all 
$\eta \in \mathbb{C}^{d_{\theta } }$, $w \in \mathbb{C}^{N_{x} }$, 
$0\leq k \leq n$ 
and any sequence 
$\boldsymbol y = \{y_{n} \}_{n\geq 1}$
from ${\cal Y}$. 
It is also easy to show
\begin{align}
	&%\label{1*.1.9}
	H_{\theta }^{0:n}(u,V,y_{1:n} )
	=
	V\: 
	(\nabla_{u} G_{\theta }^{0:n} )(u,y_{1:n} )
	+
	(\nabla_{\theta } G_{\theta }^{0:n} )(u,y_{1:n} ), 
	\nonumber\\
	&\label{1*.1.21}
	\begin{aligned}[b]
		F_{\theta }\left(
		G_{\theta }^{0:n}(u,y_{1:n} ), 
		H_{\theta }^{0:n}(u,V,y_{1:n} ), 
		y_{n+1} 
		\right)
		= &
		V\: 
		(\nabla_{u} G_{\theta }^{0:n} )(u,y_{1:n} ) \: 
		(\nabla_{u} \phi_{\theta} )(G_{\theta }^{0:n}(u,y_{1:n} ), y_{n+1} ) 
		\\
		&
		+
		\nabla_{\theta }
		\big(\phi_{\theta }(G_{\theta }^{0:n}(u,y_{1:n} ), y_{n+1} ) \big)
	\end{aligned}	
\end{align}
for each $\theta \in \Theta$, $u \in [0,\infty )^{N_{x} } \setminus \{0\}$, 
$V \in \mathbb{R}^{d_{\theta } \times N_{x} }$, $n \geq 0$
and any sequence 
$\boldsymbol y = \{y_{n} \}_{n\geq 1}$ from ${\cal Y}$
($\phi_{\theta}(u,y)$, $F_{\theta}(u,V,y)$ 
are Section \ref{section1}; 
$(\nabla_{u} G_{\theta}^{0:n} )(u,y)$, 
$(\nabla_{\theta} G_{\theta}^{0:n} )(u,y)$ denote
the Jacobians of 
$G_{\theta}^{0:n}(u,y)$ with respect to $u$, $\theta$, 
while $(\nabla_{u} \phi_{\theta } )(u,y)$ stands for the gradient of 
$\phi_{\theta }(u,y)$ with respect to $u$). 

Besides the previously introduced notation, 
the following notation is also used in this section. 
For $u \in [0,\infty )^{N_{x} } \setminus \{0\}$, 
$n > m \geq 0$
and sequences 
$\boldsymbol \vartheta = \{\vartheta_{n} \}_{n\geq 0}$, 
$\boldsymbol y = \{y_{n} \}_{n\geq 1}$
from $\Theta$, ${\cal Y}$ (respectively), 
let 
$A_{\boldsymbol \vartheta, \boldsymbol y }^{n:n}(u) = I \in \mathbb{R}^{N_{x} \times N_{x} }$
($I$ denotes a unit matrix) and 
\begin{align*}
	&
	A_{\boldsymbol \vartheta, \boldsymbol y }^{m:n}(u) 
	=
	(\nabla_{u} G_{\vartheta_{m+1} } )
	(G_{\boldsymbol \vartheta, \boldsymbol y }^{m:m}(u), y_{m+1} )
	\cdots 
	(\nabla_{u} G_{\vartheta_{n} } )
	(G_{\boldsymbol \vartheta, \boldsymbol y }^{m:n-1}(u), y_{n} ). 
\end{align*}
Then, it is easy to demonstrate
\begin{align}\label{1*.1.23}
	H_{\boldsymbol \vartheta, \boldsymbol y }^{m:n}(u,V)
	=
	V\: A_{\boldsymbol \vartheta, \boldsymbol y }^{m:n}(u) 
	+
	\sum_{i=m}^{n-1} 
	(\nabla_{\theta } G_{\vartheta_{i+1} } )
	(G_{\boldsymbol \vartheta, \boldsymbol y }^{m:i}(u), y_{i+1} )
	\:
	A_{\boldsymbol \vartheta, \boldsymbol y }^{i+1:n}
	(G_{\boldsymbol \vartheta, \boldsymbol y }^{m:i+1}(u)) 
\end{align}
for each $u \in [0,\infty )^{N_{x} } \setminus \{0\}$, 
$V \in \mathbb{R}^{d_{\theta } \times N_{x} }$, $n\geq m \geq 0$
and any sequences 
$\boldsymbol \vartheta = \{\vartheta_{n} \}_{n\geq 0}$, 
$\boldsymbol y = \{y_{n} \}_{n\geq 1}$
from $\Theta$, ${\cal Y}$ (respectively). 

In this subsection, we also rely on the following notation. 
${\cal S}_{z}$ and ${\cal S}_{\zeta }$ 
denote sets 
${\cal S}_{z} = 
{\cal X} \times {\cal Y} \times {\cal P}^{N_{x} } \times \mathbb{R}^{d_{\theta } \times N_{x} }$
and 
${\cal S}_{\zeta} = 
{\cal X} \times {\cal Y} \times {\cal P}^{N_{x} }$. 
%${\cal M}^{N_{x} } = \{u \in [0,\infty )^{N_{x} }: e^{T} u > 0\}$
For $\theta \in \Theta$, 
$P_{\theta }(\cdot, \cdot )$ and 
$\tilde{P}_{\theta }(\cdot, \cdot )$ are the transition kernels of 
Markov chains 
\begin{align*}
\{X_{n+1}, Y_{n+1}, G_{\theta }^{0:n}(u, Y_{1:n} ), 
H_{\theta }^{0:n}(u,V,Y_{1:n} ) \}_{n\geq 0}
\;\;\; \text{ and } \;\;\; 
\{X_{n}, Y_{n}, G_{\theta }^{0:n}(u, Y_{1:n} ), 
H_{\theta }^{0:n}(u,V,Y_{1:n} ) \}_{n\geq 0} 
\end{align*}
(respectively), while 
$\Pi_{\theta }(\cdot, \cdot )$ and 
$\tilde{\Pi}_{\theta }(\cdot, \cdot )$ are the transition kernels of 
Markov chains 
\begin{align*}
\{X_{n+1}, Y_{n+1}, G_{\theta }^{0:n}(u, Y_{1:n} ) \}_{n\geq 0}
\;\;\; \text{ and } \;\;\; 
\{X_{n}, Y_{n}, G_{\theta }^{0:n}(u, Y_{1:n} ) \}_{n\geq 0}
\end{align*}
(notice that 
$P_{\theta }(\cdot, \cdot )$, $\tilde{P}_{\theta }(\cdot, \cdot )$, 
$\Pi_{\theta }(\cdot, \cdot )$, $\tilde{\Pi}_{\theta }(\cdot, \cdot )$
do not depend on $u$, $V$). 
For $\theta \in \Theta$, $z=(x,y,u,V) \in {\cal S}_{z}$, 
$\zeta = (x,y,u) \in {\cal S}_{\zeta }$, let 
\begin{align*}
	&
	\tilde{F}_{\theta }(u,V,x) 
	=
	E(F_{\theta }(u,V,Y_{2} )|X_{1}=x ), 
	\\
	&
	\tilde{\phi}_{\theta }(u,x) 
	=
	E(\phi_{\theta }(u,Y_{2} )|X_{1}=x )
\end{align*}
while 
\begin{align*}
F(\theta, z ) = F_{\theta }(u,V,y), 
\;\;\;
\tilde{F}(\theta, z ) = \tilde{F}_{\theta }(u,V,y), 
\;\;\;
\phi(\theta, \zeta ) = \phi_{\theta }(u,y), 
\;\;\; 
\tilde{\phi}(\theta, \zeta ) = \tilde{\phi}_{\theta }(u,y). 
\end{align*}
Then, it is straightforward to verify 
\begin{align}
	&\label{1*.1.25}
	\begin{aligned}[b]
		(P^{n} F)(\theta, z )
		= &
		E\left(
		F_{\theta }(G_{\theta }^{0:n}(u,Y_{1:n} ), H_{\theta }^{0:n}(u,V,Y_{1:n} ), Y_{n+1} )
		|X_{1}=x, Y_{1}=y
		\right)
		\\
		= &
		E\left(
		\tilde{F}_{\theta }(G_{\theta }^{0:n}(u,Y_{1:n} ), H_{\theta }^{0:n}(u,V,Y_{1:n} ), X_{n} )
		|X_{1}=x, Y_{1}=y
		\right)
		\\
		= &
		(\tilde{P}^{n-1} \tilde{F} )
		\Big(
		\theta, \big(\, x,y,G_{\theta }(u,y), H_{\theta }(u,V,y) \,\big)
		\Big), 
	\end{aligned} 
	\\
	&\label{1*.1.27}
	\begin{aligned}[b]
		(\Pi^{n} \phi)(\theta, \zeta )
		= &
		E\left(
		\phi_{\theta }(G_{\theta }^{0:n}(u,Y_{1:n} ), Y_{n+1} )
		|X_{1}=x, Y_{1}=y
		\right)
		\\
		= &
		E\left(
		\tilde{\phi}_{\theta }(G_{\theta }^{0:n}(u,Y_{1:n} ), X_{n} )
		|X_{1}=x, Y_{1}=y
		\right)
		\\
		= &
		(\tilde{\Pi}^{n-1} \tilde{\phi} )
		\Big(
		\theta, \big(\, x,y,G_{\theta }(u,y) \,\big)
		\Big)
	\end{aligned}
\end{align}
for all 
$\theta \in \Theta$, 
$z = (x,y,u,V ) \in {\cal S}_{z}$, 
$\zeta = (x,y,u) \in {\cal S}_{\zeta }$, 
$n>1$. 
It can also be concluded 
\begin{align}\label{1*.1.29}
	&
	E\left(\left.
	\frac{\log p_{\theta }^{n+1}(Y_{1},\dots, Y_{n+1} ) }{n+1}
	\right|X_{1}=x, Y_{1}=y 
	\right)
	\nonumber\\
	&
	=
	E\left(\left.
	\frac{1}{n+1}
	\sum_{i=0}^{n} \phi_{\theta}( G_{\theta}^{0:i}(u_{\theta}, Y_{1:i} ), Y_{i+1} )
	\right|X_{1}=x, Y_{1}=y 
	\right)
	\nonumber\\
	&
	=
	\frac{1}{n+1} 
	\sum_{i=1}^{n} 
	(\tilde{\Pi}^{i-1} \tilde{\phi} )
	\Big( \theta, \big(\, x,y,G_{\theta }(u_{\theta} , y ) \,\big) \Big)
	+
	\frac{\phi_{\theta }(u_{\theta},Y_{1} ) }{n+1}
\end{align}
for each 
$\theta \in \Theta$, 
$\zeta = (x,y,u) \in {\cal S}_{\zeta }$, 
$n>1$, 
where 
$u_{\theta } = 
[P(X_{1}^{\theta } = 1 ) \cdots P(X_{1}^{\theta } = N_{x} ) ]^{T}$. 

\begin{lemma} \label{lemma1.1} 
Suppose that Assumption \ref{a4} hold. 
Let $Q\subset \Theta$ be an arbitrary compact set. 
Then, there exist 
real numbers $\delta_{1,Q } \in (0,1)$, $C_{1,Q } \in [1,\infty )$ such 
that 
\begin{align}
	& \label{l1.1.1*}
	|\tilde{\phi}_{\theta}(u,x) |
	\leq 
	C_{1,Q }, 
	\\
	& \label{l1.1.7*}
	\|F_{\theta }(u,V,y)\| 
	\leq 
	C_{1,Q } \psi_{Q}(y) (1 + \|V \| ), 
	\\
	& \label{l1.1.9*}
	\|\tilde{F}_{\theta }(u,V,x) \|  
	\leq 
	C_{1,Q } (1 + \|V \| ), 
	\\
	& \label{l1.1.3*}
	|\tilde{\phi}_{\theta' }(u',x) - \tilde{\phi}_{\theta'' }(u'',x) |
	\leq 
	C_{1,Q } (\|\theta' - \theta'' \| + \|u' - u'' \|), 
	\\
	& \label{l1.1.5*}
	|\hat{\phi}_{\eta' }(w',y) - \hat{\phi}_{\eta'' }(w'',y) |
	\leq 
	C_{1,Q } \psi_{Q}(y) (\|\eta' - \eta'' \| + \|w' - w'' \|), 
	\\
	& \label{l1.1.21*}	
	\|F_{\theta'}(u',V',y) - F_{\theta''}(u'',V'',y) \|
	\nonumber\\
	&\;\;\; 
	\leq 
	C_{1,Q } 
	\psi_{Q}(y) 
	(1 + \|V' \| + \|V'' \| ) 
	(\|\theta' - \theta'' \| + \|u' - u'' \| + \|V' - V'' \| ), 
	 \\
	& \label{l1.1.23*}	
	\|\tilde{F}_{\theta'}(u',V',x) - \tilde{F}_{\theta''}(u'',V'',x) \|
	\nonumber\\
	&\;\;\; 
	\leq 
	C_{1,Q } 
	(1 + \|V' \| + \|V'' \| ) 
	(\|\theta' - \theta'' \| + \|u' - u'' \| + \|V' - V'' \| ) 
\end{align}
for all 
$\theta, \theta', \theta'' \in Q$, 
$\eta', \eta'' \in V_{\delta_{1,Q} }(Q)$, 
$u, u', u'' \in {\cal P}^{N_{x} }$, 
$w', w'' \in V_{\delta_{1,Q} }({\cal P}^{N_{x} } )$, 
$V, V', V'' \in \mathbb{R}^{d_{\theta } \times N_{x} }$, 
$x \in {\cal X}$, 
$y \in {\cal Y}$
($\psi_{Q}(\cdot )$ is specified in Assumption \ref{a4}).  
\end{lemma}

\begin{IEEEproof}
Let $\delta_{1,Q } = \delta_{Q }/2$
($\delta_{Q}$ is defined in Assumption \ref{a4}). 
Then, Cauchy inequality for analytic functions (see e.g., \cite[Proposition 2.1.3]{taylor}) 
and Assumption \ref{a4}
imply that there exists a real number 
$\tilde{C}_{1,Q } \in [1, \infty )$ such that 
\begin{align*}
	\max\{
	\|\nabla_{(\eta, w )}\hat{\phi}_{\eta}(w,y) \|, 
	\|\nabla_{(\eta, w )}^{2}\hat{\phi}_{\eta}(w,y) \|
	\}
	\leq
	\tilde{C}_{1,Q } \psi_{Q}(y) 
\end{align*}
for all $\eta \in V_{\delta_{1,Q } }(Q)$, 
$w \in V_{\delta_{1,Q } }({\cal P}^{N_{x} } )$, 
$y \in {\cal Y}$
($\nabla_{(\eta, w )}$, $\nabla_{(\eta, w )}^{2}$ 
denote the gradient and Hessian with respect to $(\eta, w )$). 
Consequently, there exists another real number 
$\tilde{C}_{2,Q } \in [1,\infty )$ such that 
\begin{align*}
	&
	\max\{
	\|\hat{\phi}_{\eta'}(w',y) 
	-
	\hat{\phi}_{\eta''}(w'',y) 
	\|, 
	\|\nabla_{w} \hat{\phi}_{\eta'}(w',y) 
	-
	\nabla_{w} \hat{\phi}_{\eta''}(w'',y) 
	\|
	\}
	\\
	&
	\leq 
	\tilde{C}_{2,Q} \psi_{Q}(y)
	(\|\eta' - \eta'' \| + \|w' - w'' \| )
\end{align*}
for any $\eta',\eta'' \in V_{\delta_{1,Q} }(Q)$, 
$w',w'' \in V_{\delta_{1,Q} }({\cal P}^{N_{x} } )$, 
$y \in {\cal Y}$. 
Therefore,
\begin{align*}
	&
	\begin{aligned}[t]
	\|F_{\theta}(u,V,y) \|
	\leq &
	\|\nabla_{\theta} \phi_{\theta}(u,y) \|
	+
	\|\nabla_{u} \phi_{\theta}(u,y) \| 
	\|V\| 
	\\
	\leq &
	\tilde{C}_{1,Q } \psi_{Q}(y) (1 + \|V\| ), 
	\end{aligned}
	\\
	&
	\begin{aligned}[b]
		\|F_{\theta'}(u',V',y) 
		-
		F_{\theta''}(u'',V'',y)\|
		\leq &
		\|\nabla_{\theta} \phi_{\theta'}(u',y) - \nabla_{\theta} \phi_{\theta''}(u'',y) \|
		+
		\|\nabla_{u} \phi_{\theta'}(u',y) - \nabla_{u} \phi_{\theta''}(u'',y) \|
		\|V'\| 
		\\
		&
		+
		\|\nabla_{u} \phi_{\theta'' } (u'',y) \|
		\|V' - V'' \| 
		\\
		\leq &
		\tilde{C}_{2,Q } \psi_{Q}(y) 
		(1 + \|V'\| + \|V''\| ) 
		(\|\theta' - \theta'' \| + \|u' - u'' \| )
		\\
		&
		+
		\tilde{C}_{1,Q } \psi_{Q}(y) \|V' - V'' \|
	\end{aligned}
\end{align*}
for each $\theta, \theta', \theta'' \in Q$, 
$u,u',u'' \in {\cal P}^{N_{x} }$, 
$V,V',V'' \in \mathbb{R}^{d_{\theta} \times N_{x} }$. 
We also have 
\begin{align*}
	&
	\|\tilde{F}_{\theta }(u,V,x) \| 
	\leq
	\tilde{C}_{1,Q } 
	(1 + \|V\| ) 
	\int \psi_{Q}(y) Q(dy|x) 
	\\
	&
	\begin{aligned}[b]
		&
		\|\tilde{F}_{\theta'}(u',V',x) 
		-
		\tilde{F}_{\theta''}(u'',V'',x)\|
		\\
		& \;\;\; 
		\leq 
		(\tilde{C}_{1,Q } + \tilde{C}_{2,Q } ) 
		(1 + \|V'\| + \|V''\| ) 
		(\|\theta' - \theta'' \| + \|u' - u'' \| + \|V' - V'' \| )
		\int \psi_{Q}(y) Q(dy|x) 
	\end{aligned}
\end{align*}
for all 
$\theta, \theta', \theta'' \in Q$, 
$u,u',u'' \in {\cal P}^{N_{x} }$, 
$V,V',V'' \in \mathbb{R}^{d_{\theta } \times N_{x} }$, 
$x\in {\cal X}$. 
Then, it can be deduced that there exists 
a real number $C_{1,Q } \in [1, \infty )$ such that 
(\ref{l1.1.1*}) -- (\ref{l1.1.23*}) hold for each 
$\theta, \theta', \theta'' \in Q$, 
$\eta', \eta'' \in V_{\delta_{1,Q} }(Q)$, 
$u, u', u'' \in {\cal P}^{N_{x} }$, 
$w', w'' \in V_{\delta_{1,Q} }({\cal P}^{N_{x} } )$, 
$V, V', V'' \in \mathbb{R}^{d_{\theta } \times N_{x} }$, 
$x \in {\cal X}$, 
$y \in {\cal Y}$.  
\end{IEEEproof}

\begin{lemma} \label{lemma1.2}
Suppose that Assumption \ref{a4} hold. 
Let $Q\subset \Theta$ be an arbitrary compact set. 
Then, there exist  
real numbers $\delta_{2,Q } \in (0,1)$, $C_{2,Q } \in [1,\infty )$ such that 
\begin{align}
	& \label{l1.2.1*}
	\|\nabla_{\eta} \hat{G}_{\eta }(w,y) \|
	\leq 
	C_{2,Q}, 
	\\
	& \label{l1.2.3*}
	\|H_{\theta }(u,V,y) \|
	\leq 
	C_{2,Q} (1 + \|V\| ), 
	\\
	& \label{l1.2.5*} 
	\begin{aligned}[b]
		&
		\max\{
		\|\hat{G}_{\eta' }(w',y) - \hat{G}_{\eta'' }(w'',y) \|, 
		\|\nabla_{w} \hat{G}_{\eta' }(w',y) - \nabla_{w} \hat{G}_{\eta'' }(w'',y) \|
		\}
		\\
		& \;\;\; 
		\leq 
		C_{2,Q } (\|\eta' - \eta'' \| + \|w' - w'' \|),  
	\end{aligned}
	\\
	& \label{l1.2.7*}
	\|H_{\theta'}(u,V,y) - H_{\theta''}(u,V,y) \| 
	\leq 
	C_{2,Q} (1 + \|V\| ) \|\theta' - \theta'' \| 
\end{align}
for all $\theta, \theta', \theta'' \in Q$, 
$\eta, \eta', \eta'' \in V_{\delta_{2,Q } }(Q )$, 
$u\in {\cal P}^{N_{x} }$, 
$w,w',w'' \in V_{\delta_{2,Q } }({\cal P}^{N_{x} } )$, 
$V\in \mathbb{R}^{d_{\theta } \times N_{x} }$, 
$y \in {\cal Y}$. 
\end{lemma}

\begin{IEEEproof}
Let $\delta_{2,Q } = \min\{\delta_{Q }/2,\delta_{1,Q } \}$. 
Owing to Cauchy inequality for analytic functions and Assumption \ref{a4}, 
there exists a real number 
$\tilde{C}_{1,Q } \in [1, \infty )$ such that 
\begin{align*}
	\max\{
	\|\nabla_{(\eta,w)} \hat{G}_{\eta}^{k}(w,y) \|, 
	\|\nabla_{(\eta,w)}^{2} \hat{G}_{\eta}^{k}(w,y) \|
	\}
	\leq 
	\tilde{C}_{1,Q }
\end{align*}
for any $\eta \in V_{\delta_{2,Q } }(Q )$, 
$w \in V_{\delta_{2,Q } }({\cal P}^{N_{x} } )$, 
$y \in {\cal Y}$
($\hat{G}_{\eta}^{k}(w,y)$ stands for the $k$-th component of 
$\hat{G}_{\eta}(w,y)$). 
Consequently, there exists another real number 
$\tilde{C}_{2,Q } \in [1, \infty )$ such that 
\begin{align*}
	&
	\max\{
	\|\hat{G}_{\eta'}(w',y) - \hat{G}_{\eta''}(w'',y) \|, 
	\|\nabla_{\eta} \hat{G}_{\eta'}(w',y) - \nabla_{\eta} \hat{G}_{\eta''}(w'',y) \|, 
	\|\nabla_{w} \hat{G}_{\eta'}(w',y) - \nabla_{w} \hat{G}_{\eta''}(w'',y) \|
	\}
	\\
	&\;\;\; 
	\leq 
	\tilde{C}_{2,Q} 
	(\|\eta' - \eta'' \| + \|w' - w'' \| )
\end{align*}
for all $\eta', \eta'' \in V_{\delta_{2,Q } }(Q )$, 
$w',w'' \in V_{\delta_{2,Q } }({\cal P}^{N_{x} } )$, 
$y \in {\cal Y}$. 
Therefore, 
\begin{align*}
	&
	\begin{aligned}
	\|H_{\theta}(u,V,y) \|
	\leq &
	\|\nabla_{\theta} G_{\theta}(u,y) \| 
	+
	\|\nabla_{u} G_{\theta }(u,y) \| \|V\| 
	\\
	\leq &
	\tilde{C}_{1,Q} N_{x} (1 + \|V\| ), 
	\end{aligned} 
	\\
	&
	\begin{aligned}
	\|H_{\theta'}(u,V,y) - H_{\theta''}(u,V,y) \|
	\leq &
	\|\nabla_{\theta} G_{\theta'}(u,y) - \nabla_{\theta} G_{\theta''}(u,y) \| 
	+
	\|\nabla_{u} G_{\theta'}(u,y) - \nabla_{u} G_{\theta''}(u,y) \| 
	\|V\| 
	\\
	\leq &
	\tilde{C}_{2,Q} (1 + \|V\| ) \|\theta' - \theta'' \|
	\end{aligned}
\end{align*}
for each $\theta, \theta', \theta'' \in Q$, 
$u \in {\cal P}^{N_{x} }$, 
$V \in \mathbb{R}^{d_{\theta } \times N_{x} }$. 
Then, it is clear that there exists a real number 
$C_{2,Q } \in [1, \infty )$ such that 
(\ref{l1.2.1*}) -- (\ref{l1.2.7*}) hold for all 
$\theta, \theta', \theta'' \in Q$, 
$\eta, \eta', \eta'' \in V_{\delta_{2,Q } }(Q )$, 
$u\in {\cal P}^{N_{x} }$, 
$w,w',w'' \in V_{\delta_{2,Q } }({\cal P}^{N_{x} } )$, 
$V\in \mathbb{R}^{d_{\theta } \times N_{x} }$, 
$y \in {\cal Y}$. 
\end{IEEEproof}

\begin{lemma} \label{lemma1.3}
Suppose that Assumptions \ref{a3} and \ref{a4} hold. 
Let $Q\subset \Theta$ be an arbitrary compact set. 
Then, 
the following is true: 
\begin{enumerate}
\item
There exist 
real numbers 
$\varepsilon_{1,Q } \in (0,1)$, 
$C_{3,Q } \in [1,\infty )$
such that 
\begin{align}
	& \label{l1.3.1*}
	\|A_{\boldsymbol \vartheta, \boldsymbol y }^{m:n}(u) \|
	\leq
	C_{3,Q } \varepsilon_{1,Q }^{n-m}, 
	\\
	& \label{l1.3.3*}
	\|A_{\boldsymbol \vartheta, \boldsymbol y }^{m:n}(u') 
	-
	A_{\boldsymbol \vartheta, \boldsymbol y }^{m:n}(u'') \|
	\leq
	C_{3,Q } \varepsilon_{1,Q }^{n-m} \|u' - u'' \|, 
	\\
	& \label{l1.3.5*}
	\|G^{m:n}_{\boldsymbol\vartheta, \boldsymbol y }(w') 
	- 
	G^{m:n}_{\boldsymbol\vartheta, \boldsymbol y }(w'') \|
	\leq 
	C_{3,Q } \varepsilon_{1,Q }^{n-m} \|w' - w'' \| 
\end{align}
for all $u,u',u'' \in {\cal P}^{N_{x} }$, 
$w', w'' \in {\cal Q}^{N_{x} }$, 
$n\geq m\geq 0$ 
and 
any sequences 
$\boldsymbol\vartheta = \{\vartheta_{n} \}_{n\geq 0}$, 
$\boldsymbol y = \{y_{n} \}_{n\geq 1}$
from $Q$, ${\cal Y}$ (respectively). 
\item
There exist real numbers $\varepsilon_{2,Q } \in (0,1)$, 
$C_{4,Q } \in [1,\infty )$ such that 
\begin{align} 
	& \label{l1.3.7*}
	\|H^{m:n}_{\boldsymbol\vartheta, \boldsymbol y }(u,V) \|
	\leq 
	C_{4,Q } (1 + \|V \| )
	\\
	& \label{l1.3.9*}	
	\|H^{m:n}_{\boldsymbol\vartheta, \boldsymbol y }(u',V') 
	- 
	H^{m:n}_{\boldsymbol\vartheta, \boldsymbol y }(u'',V'') \|
	\leq 
	C_{4,Q } \varepsilon_{2,Q }^{n-m}
	\big(\|u' - u'' \|(1 + \|V' \| + \|V'' \| ) + \|V' - V'' \| \big) 
\end{align}
for all $u,u',u'' \in {\cal P}^{N_{x} }$, 
$V, V', V'' \in \mathbb{R}^{d_{\theta } \times N_{x} }$, 
$n\geq m\geq 0$ 
and 
any sequences 
$\boldsymbol\vartheta = \{\vartheta_{n} \}_{n\geq 0}$, 
$\boldsymbol y = \{y_{n} \}_{n\geq 1}$
from $Q$, ${\cal Y}$ (respectively). 
\end{enumerate}
\end{lemma}

\begin{IEEEproof}
Using \cite[Theorem 3.1, Lemmas 6.6, 6.7]{tadic&doucet}
(with a few straightforward modifications), 
it can be deduced from Assumption \ref{a3} that there exist 
real numbers $\varepsilon_{1,Q} \in (0,1)$, 
$C_{3,Q} \in [1,\infty )$ such that 
(\ref{l1.3.1*}), (\ref{l1.3.3*}) and 
\begin{align*}
	\|G^{m:n}_{\boldsymbol\vartheta, \boldsymbol y }(w') 
	- 
	G^{m:n}_{\boldsymbol\vartheta, \boldsymbol y }(w'') \|
	\leq 
	2^{-1} (N_{x} + 1 )^{-1} C_{3,Q } \varepsilon_{1,Q }^{n-m} 
	\left\|
	\frac{w'}{e^{T} w' } - \frac{w''}{e^{T} w'' }
	\right\| 
\end{align*}
hold for all $u,u',u'' \in {\cal P}^{N_{x} }$, 
$w', w'' \in [0,\infty )^{N_{x} } \setminus \{0\}$, 
$n\geq m\geq 0$ 
and 
any sequences 
$\boldsymbol\vartheta = \{\vartheta_{n} \}_{n\geq 0}$, 
$\boldsymbol y = \{y_{n} \}_{n\geq 1}$
from $Q$, ${\cal Y}$.\footnote{
To deduce this, note that 
$u$, $V$, $y_{0:n}$, $G_{\boldsymbol\vartheta, \boldsymbol y }^{0:n}(u)$,  
$A_{\boldsymbol\vartheta, \boldsymbol y }^{0:n}(u) V$
have the same meaning respectively as quantities 
$\mu$, $\tilde{\mu}$, $y^{n}$, 
$F_{\theta }^{n}(\mu, y^{n} )$, 
$\tilde{G}_{\theta}^{n}(\mu,\tilde{\mu},y^{n} )$ 
appearing in \cite{tadic&doucet}.}\footnote{
Inequality (\ref{l1.3.5*}) can also be obtained 
from \cite[Theorem 2.1]{legland&mevel2} or 
\cite[Theorem 4.1]{legland&oudjane}. 
Similarly, (\ref{l1.3.1*}), (\ref{l1.3.3*}) 
can be deduced from 
\cite[Lemmas 3.4, 4.3, Proposition 5.2]{legland&mevel1}
(notice that $G^{m:n}_{\boldsymbol\vartheta, \boldsymbol y }(u)$, 
$A^{m:n}_{\boldsymbol\vartheta, \boldsymbol y }(u)$
have the same meaning respectively as 
$M_{m,n}$, $V[M_{m,n},p_{m}]$ specified in 
\cite[Section 5]{legland&mevel1}).}  
Since 
\begin{align*}
	\left\|
	\frac{w'}{e^{T} w' } - \frac{w''}{e^{T} w'' }
	\right\| 
	\leq 
	\frac{\|w' - w'' \| (e^{T} w'' )  
	+
	\|w'' \| \: |e^{T} (w' - w'' ) | }
	{(e^{T} w' ) (e^{T} w'' ) }
	\leq 
	2 (N_{x} + 1 ) \|w' - w'' \|
\end{align*}
for any $w', w'' \in {\cal Q}^{N_{x} }$, 
we have that 
(\ref{l1.3.1*}) is satisfied for all 
$w', w'' \in {\cal Q}^{N_{x} }$, 
$n\geq m\geq 0$ 
and 
any sequences 
$\boldsymbol\vartheta = \{\vartheta_{n} \}_{n\geq 0}$, 
$\boldsymbol y = \{y_{n} \}_{n\geq 1}$
from $Q$, ${\cal Y}$. 
Hence, (i) is true. 

Now, we shaw that (ii) is true, too. 
Let  
$\boldsymbol\vartheta = \{\vartheta_{n} \}_{n\geq 0}$, 
$\boldsymbol y = \{y_{n} \}_{n\geq 1}$
be arbitrary sequences 
from $Q$, ${\cal Y}$ (respectively). 
As a consequence of Lemma \ref{lemma1.2}, (i) and (\ref{1*.1.23}), 
we get
\begin{align*}
	\|H_{\boldsymbol \vartheta, \boldsymbol y}^{m:n}(u,V) \| 
	\leq 
	C_{3,Q } \varepsilon_{1,Q }^{n-m} \|V\| 
	+
	C_{2,Q } C_{3,Q } 
	\sum_{i=m}^{n-1} \varepsilon_{1,Q }^{n-i-1} 
	\leq 
	C_{3,Q } \|V\| 
	+
	C_{2,Q } C_{3,Q } (1 - \varepsilon_{1,Q } )^{-1}
\end{align*}
for all $u \in {\cal P}^{N_{x} }$, 
$V \in \mathbb{R}^{d_{\theta } \times N_{x} }$, 
$n\geq m \geq 0$. 
Due to the same arguments, we have 
\begin{align*}
	&
	\|H_{\boldsymbol \vartheta, \boldsymbol y }^{m:n}(u',V') 
	-
	H_{\boldsymbol \vartheta, \boldsymbol y }^{m:n}(u'',V'')  \| 
	\\
	&\;\;\; 
	\begin{aligned}[b]
	\leq &
	\|A_{\boldsymbol \vartheta, \boldsymbol y }^{m:n}(u') 
	-
	A_{\boldsymbol \vartheta, \boldsymbol y }^{m:n}(u'') \| \|V' \|
	+
	\|A_{\boldsymbol \vartheta, \boldsymbol y }^{m:n}(u'') \| \|V' - V'' \|
	\\
	&
	+
	\sum_{i=m}^{n-1} 
	\|(\nabla_{\theta } G_{\vartheta_{i+1} } )
	(G_{\boldsymbol \vartheta, \boldsymbol y }^{m:i}(u'),y_{i+1} ) 
	- 
	(\nabla_{\theta } G_{\vartheta_{i+1} } )
	(G_{\boldsymbol \vartheta, \boldsymbol y }^{m:i}(u''),y_{i+1} ) \|
	\|A_{\boldsymbol \vartheta, \boldsymbol y }^{i+1:n}
	(G_{\boldsymbol \vartheta, \boldsymbol y }^{m:i+1}(u') ) \|
	\\
	&
	+
	\sum_{i=m}^{n-1} 
	\|(\nabla_{\theta } G_{\vartheta_{i+1} } )
	(G_{\boldsymbol \vartheta, \boldsymbol y }^{m:i}(u''),y_{i+1} ) \|
	\|A_{\boldsymbol \vartheta, \boldsymbol y }^{i+1:n}
	(G_{\boldsymbol \vartheta, \boldsymbol y }^{m:i+1}(u') ) 
	-
	A_{\boldsymbol \vartheta, \boldsymbol y }^{i+1:n}
	(G_{\boldsymbol \vartheta, \boldsymbol y }^{m:i+1}(u'') ) \|
	\end{aligned}
	\\
	& \;\;\; 
	\begin{aligned}[b]
	\leq &
	C_{3,Q } 
	\varepsilon_{1,Q }^{n-m} \|V'\| \|u' - u'' \| 
	+
	C_{3,Q } \varepsilon_{1,Q }^{n-m} \|V' - V'' \| 
%	\\
%	&
	+
	C_{2,Q } C_{3,Q } 
	\sum_{i=m}^{n-1} 
	\varepsilon_{1,Q }^{n-i-1} 
	\|G_{\boldsymbol \vartheta, \boldsymbol y }^{m:i}(u') 
	-
	G_{\boldsymbol \vartheta, \boldsymbol y }^{m:i}(u'') \|
	\\
	&
	+
	C_{2,Q } C_{3,Q } 
	\sum_{i=m}^{n-1} 
	\varepsilon_{1,Q }^{n-i-1} 
	\|G_{\boldsymbol \vartheta, \boldsymbol y }^{m:i+1}(u') 
	-
	G_{\boldsymbol \vartheta, \boldsymbol y }^{m:i+1}(u'') \|
	\end{aligned} 
	\\
	& \;\;\; 
	\leq 
	C_{3,Q } \varepsilon_{1,Q }^{n-m} 
	(\|u' - u'' \| \|V'\| + \|V' - V'' \| )
	+
	2 C_{2,Q } C_{3,Q }^{2} \varepsilon_{1,Q }^{n-m-1} (n-m)
\end{align*}
for each $u',u'' \in {\cal P}^{N_{x} }$, 
$V', V'' \in \mathbb{R}^{d_{\theta } \times N_{x} }$, 
$n\geq m \geq 0$. 
Then, it is clear that there exist real numbers 
$\varepsilon_{2,Q } \in (0,1)$, 
$C_{4,Q } \in [1,\infty )$
such that (\ref{l1.3.7*}), (\ref{l1.3.9*}) 
hold for all 
$u,u',u'' \in {\cal P}^{N_{x} }$, 
$V,V',V'' \in \mathbb{R}^{d_{\theta } \times N_{x} }$
and any sequence 
$\boldsymbol \vartheta = \{\vartheta \}_{n\geq 0}$, 
$\boldsymbol y = \{y_{n} \}_{n\geq 1}$
from $Q$, ${\cal Y}$
(respectively). 
\end{IEEEproof}

\begin{lemma} \label{lemma1.4}
Suppose that Assumptions \ref{a3} and \ref{a4} hold. 
Let $Q\subset \Theta$ be an arbitrary compact set. 
Then, there exists 
a real number $C_{5,Q } \in [1, \infty )$ such that 
\begin{align}
	& \label{l1.4.1*}
	\|G_{\theta', \boldsymbol y }^{0:n}(u) 
	-
	G_{\theta'', \boldsymbol y }^{0:n}(u) 
	\|
	\leq
	C_{5,Q } \|\theta' - \theta'' \|, 
	\\
	& \label{l1.4.3*}
	\|H_{\theta', \boldsymbol y }^{0:n}(u,V) 
	-
	H_{\theta'', \boldsymbol y }^{0:n}(u,V) 
	\|
	\leq
	C_{5,Q } \|\theta' - \theta'' \| (1 + \|V\| ) 
\end{align}
for all $\theta', \theta'' \in Q$, 
$u\in {\cal P}^{N_{x} }$, $V\in \mathbb{R}^{d_{\theta } \times N_{x} }$, 
$n\geq 1$
and any sequence 
$\boldsymbol y = \{y_{n} \}_{n\geq 1}$ from ${\cal Y}$. 
\end{lemma}

\begin{IEEEproof}
Let 
$\tilde{C}_{Q } = C_{2,Q } C_{3,Q } C_{4,Q }^{2}$, 
while 
$\boldsymbol y = \{y_{n} \}_{n\geq 0}$
is an arbitrary sequence from ${\cal Y}$. 
It is straightforward to verify 
\begin{align}
	& \label{l1.4.1}
	G_{\theta', \boldsymbol y }^{0:n}(u) 
	-
	G_{\theta'', \boldsymbol y }^{0:n}(u) 
	=
	\sum_{i=0}^{n-1}
	\left(
	G_{\theta', \boldsymbol y }^{i:n}(G_{\theta'', \boldsymbol y }^{0:i}(u) ) 
	-
	G_{\theta', \boldsymbol y }^{i+1:n}(G_{\theta'', \boldsymbol y }^{0:i+1}(u) ) 
	\right), 
	\\
	& \label{l1.4.3}
	H_{\theta', \boldsymbol y }^{0:n}(u,V) 
	-
	H_{\theta'', \boldsymbol y }^{0:n}(u,V) 
%	\nonumber \\
%	&\;\;\; 
	=
	\sum_{i=0}^{n-1}
	\left(
	H_{\theta', \boldsymbol y }^{i:n}(
	G_{\theta'', \boldsymbol y }^{0:i}(u), H_{\theta'', \boldsymbol y }^{0:i}(u,V) ) 
	-
	H_{\theta', \boldsymbol y }^{i+1:n}(
	G_{\theta'', \boldsymbol y }^{0:i+1}(u), H_{\theta'', \boldsymbol y }^{0:i+1}(u,V) ) 
	\right)
\end{align}
for 
all $\theta',\theta'' \in Q$, 
$u \in {\cal P}^{N_{x} }$, $V \in \mathbb{R}^{d_{\theta } \times N_{x} }$, 
$n\geq 0$. 
On the other side, Lemmas \ref{lemma1.2} and \ref{lemma1.3} yield
\begin{align}
	\label{l1.4.5}
	\|
	G_{\theta', \boldsymbol y }^{i:n}(G_{\theta'', \boldsymbol y }^{0:i}(u) ) 
	-
	G_{\theta', \boldsymbol y }^{i+1:n}(G_{\theta'', \boldsymbol y }^{0:i+1}(u) ) 
	\|
%	\nonumber\\
%	&\;\;\; 
%	\begin{aligned}[t]
		=&
		\left\|
		G_{\theta', \boldsymbol y }^{i+1:n}\left(
		G_{\theta',\boldsymbol y }^{i:i+1}(
		G_{\theta'', \boldsymbol y }^{0:i}(u) ) \right)
		-
		G_{\theta', \boldsymbol y }^{i+1:n}\left(
		G_{\theta'',\boldsymbol y }^{i:i+1}(
		G_{\theta'', \boldsymbol y }^{0:i}(u) ) \right)
		\right\|
		\nonumber\\
		\leq &
		C_{3,Q } \varepsilon_{1,Q }^{n-i-1}
		\left\|
		G_{\theta',\boldsymbol y }^{i:i+1}(
		G_{\theta'', \boldsymbol y }^{0:i}(u) ) 
		-
		G_{\theta'',\boldsymbol y }^{i:i+1}(
		G_{\theta'', \boldsymbol y }^{0:i}(u) ) 
		\right\|
		\nonumber\\
		\leq &
		\tilde{C}_{Q } \varepsilon_{1,Q }^{n-i-1} \|\theta' - \theta'' \| 
%	\end{aligned}
\end{align}
for any 
$\theta',\theta'' \in Q$, 
$u \in {\cal P}^{N_{x} }$, $V \in \mathbb{R}^{d_{\theta } \times N_{x} }$, 
$0\leq i < n$. Using the same lemmas, we also get 
\begin{align}	
	&\label{l1.4.7}
	\|
	H_{\theta', \boldsymbol y }^{i:n}(
	G_{\theta'', \boldsymbol y }^{0:i}(u), H_{\theta'', \boldsymbol y }^{0:i}(u,V) ) 
	-
	H_{\theta', \boldsymbol y }^{i+1:n}(
	G_{\theta'', \boldsymbol y }^{0:i+1}(u), H_{\theta'', \boldsymbol y }^{0:i+1}(u,V) ) 
	\|
	\nonumber\\
	&\hspace{0.41em} 
	\begin{aligned}[b]
		=
		&
		\left\|
		H_{\theta', \boldsymbol y }^{i+1:n}\left(
		G_{\theta',\boldsymbol y }^{i:i+1}\left(
		G_{\theta'', \boldsymbol y }^{0:i}(u) \right), 
		H_{\theta',\boldsymbol y }^{i:i+1}\left(
		G_{\theta'', \boldsymbol y }^{0:i}(u) ), H_{\theta'', \boldsymbol y }^{0:i}(u,V) \right) 	
		\right) 
		\right.
		\\
		&			
		-
		\left.
		H_{\theta', \boldsymbol y }^{i+1:n}\left(
		G_{\theta'',\boldsymbol y }^{i:i+1}\left(
		G_{\theta'', \boldsymbol y }^{0:i}(u) \right), 
		H_{\theta'',\boldsymbol y }^{i:i+1}\left(
		G_{\theta'', \boldsymbol y }^{0:i}(u) ), H_{\theta'', \boldsymbol y }^{0:i}(u,V) \right) 	
		\right) 
		\right\|
		\\
		\leq &
		C_{4,Q } \varepsilon_{2,Q }^{n-i-1}
		\left\|		
		G_{\theta',\boldsymbol y }^{i:i+1}(
		G_{\theta'', \boldsymbol y }^{0:i}(u) ) 
		-
		G_{\theta'',\boldsymbol y }^{i:i+1}(
		G_{\theta'', \boldsymbol y }^{0:i}(u) ) 
		\right\|
		\\
		&
		\cdot 
		\left(
		1
		+
		\left\|
		H_{\theta',\boldsymbol y }^{i:i+1}\left(
		G_{\theta'', \boldsymbol y }^{0:i}(u) ), H_{\theta'', \boldsymbol y }^{0:i}(u,V) \right) 
		\right\|
		+
		\left\|
		H_{\theta'',\boldsymbol y }^{i:i+1}\left(
		G_{\theta'', \boldsymbol y }^{0:i}(u) ), H_{\theta'', \boldsymbol y }^{0:i}(u,V) \right) 
		\right\|
		\right)
		\\
		&
		+
		C_{4,Q } \varepsilon_{2,Q }^{n-i-1}
		\left\|
		H_{\theta',\boldsymbol y }^{i:i+1}\left(
		G_{\theta'', \boldsymbol y }^{0:i}(u) ), H_{\theta'', \boldsymbol y }^{0:i}(u,V) \right) 
		-
		H_{\theta'',\boldsymbol y }^{i:i+1}\left(
		G_{\theta'', \boldsymbol y }^{0:i}(u) ), H_{\theta'', \boldsymbol y }^{0:i}(u,V) \right) 
		\right\|
		\\
		\leq &
		3 C_{2,Q } C_{4,Q }^{2} \varepsilon_{2,Q }^{n-i-1}
		\|\theta' - \theta'' \|
		(1 + \|V\| ) 
		+
		C_{2,Q } C_{4,Q } \varepsilon_{2,Q }^{n-i-1}
		\|\theta' - \theta'' \|
		(1 + \|H_{\theta'', \boldsymbol y }^{0:i}(u,V) \|)
	\end{aligned}
	\nonumber\\
	&\hspace{0.41em}
	\leq 
	5 \tilde{C}_{Q } \varepsilon_{2,Q}^{n-i-1} \|\theta' - \theta'' \| 
	(1 + \|V\| )
\end{align}
for each 
$\theta',\theta'' \in Q$, 
$u \in {\cal P}^{N_{x} }$, $V \in \mathbb{R}^{d_{\theta } \times N_{x} }$, 
$0\leq i < n$. Combining (\ref{l1.4.1}) -- (\ref{l1.4.7}), 
we conclude that there exists a real number 
$C_{5,Q } \in [1, \infty )$ such that 
(\ref{l1.4.1*}), (\ref{l1.4.3*}) hold 
for all $\theta', \theta'' \in Q$, 
$u\in {\cal P}^{N_{x} }$, $V\in \mathbb{R}^{d_{\theta } \times N_{x} }$, 
$n\geq 1$
and any sequence 
$\boldsymbol y = \{y_{n} \}_{n\geq 1}$ from ${\cal Y}$.  
\end{IEEEproof}

\begin{lemma} \label{lemma1.5}
Suppose that 
Assumptions \ref{a2} -- \ref{a4} hold. 
Let $Q\subset \Theta$ be an arbitrary compact set. 
Then, the following is true: 
\begin{enumerate}
\item
$f(\cdot )$ is well-defined and differentiable
on $Q$. 
\item
There exist real numbers $\varepsilon_{3,Q } \in (0,1)$, 
$C_{6,Q } \in [1, \infty )$ such that 
\begin{align*}
	&
	\|(P^{n} F)(\theta, z ) - \nabla f(\theta ) \|
	\leq 
	C_{6,Q } \varepsilon_{3,Q }^{n} 
	(1 + \|V\|^{2} ), 
	\\
	&
	|(\Pi^{n} \phi)(\theta, \zeta ) - f(\theta ) |
	\leq 
	C_{6,Q } \varepsilon_{3,Q }^{n}
\end{align*}
for all $\theta \in Q$, 
$z = (x,y,u,V) \in {\cal S}_{z}$, 
$\zeta = (x,y,u) \in {\cal S}_{\zeta }$, 
$n\geq 1$. 
\end{enumerate}
\end{lemma}

\begin{IEEEproof}
Using \cite[Theorems 4.1, 4.2]{tadic&doucet}
(with a few straightforward modifications), 
it can be deduced from Lemma \ref{lemma1.1} 
that there exist functions 
$g: \Theta \rightarrow \mathbb{R}^{d_{\theta } }$, 
$\psi: \Theta \rightarrow \mathbb{R}$ and real numbers 
$\varepsilon_{3,Q } \in (0,1)$, 
$C_{6,Q} \in [1,\infty )$ such that 
\begin{align}
	& \label{l1.5.1}
	\|(\tilde{P}^{n} \tilde{F} )(\theta, z ) - g(\theta ) \|
	\leq 
	C_{6,Q } \varepsilon_{3,Q }^{n} 
	(1 + \|V\|^{2} ), 
	\\
	& \label{l1.5.3} 
	|(\tilde{\Pi}^{n} \tilde{\phi} )(\theta, \zeta ) - \psi(\theta ) |
	\leq 
	C_{6,Q } \varepsilon_{3,Q }^{n}
\end{align}
for all 
$\theta \in Q$, 
$z = (x,y,u,V) \in {\cal S}_{z}$, 
$\zeta = (x,y,u) \in {\cal S}_{\zeta }$, 
$n\geq 1$.\footnote
{The same result can also be obtained from \cite[Theorem 5.4]{legland&mevel2}}  
Since 
$E|\phi_{\theta }(u_{\theta}, Y_{1} ) | < \infty $
for any $\theta \in Q$
(due to Assumption \ref{a4}), 
it follows from (\ref{1*.1.29}), (\ref{l1.5.3})
that $f(\cdot )$ is well-defined and identical to 
$\psi(\cdot )$ on $Q$. 
On the other side, 
Lemmas \ref{lemma1.1}, \ref{lemma1.3} yield 
\begin{align*}
	\|F_{\theta}(G_{\theta }^{0:n}(u,y_{1:n} ), H_{\theta}^{0:n}(u,V,y_{1:n} ), y_{n+1} ) \|
	\leq &
	C_{1,Q } \psi_{Q}(y_{n+1} ) 
	(1 + \|H_{\theta}^{0:n}(u,V,y_{1:n} ) \| )
	\\
	\leq &
	2 C_{1,Q } C_{4,Q } \psi_{Q}(y_{n+1} ) 
	(1 + \|V\| )
\end{align*}
for each $\theta \in Q$, 
$u \in {\cal P}^{N_{x} }$, 
$V \in \mathbb{R}^{d_{\theta } \times N_{x} }$
and any sequence 
$\boldsymbol y = \{y_{n} \}_{n\geq 1 }$
from ${\cal Y}$. 
Then, Assumption \ref{a4} gives 
\begin{align*}
	&
	E\left(\left. 
	\|F_{\theta}(G_{\theta }^{0:n}(u,Y_{1:n} ), H_{\theta}^{0:n}(u,V,Y_{1:n} ), Y_{n+1} ) \|
	\right|
	X_{1}=x, Y_{1}=y
	\right)
	\\
	&
	\leq 
	2 C_{1,Q } C_{4,Q } (1 + \|V\| ) 
	\max_{x'\in {\cal X} }
	\int \psi_{Q }(y') Q(dy'|x')
	< \infty
\end{align*}
for all 
$\theta \in Q$, 
$u \in {\cal P}^{N_{x} }$, $V \in \mathbb{R}^{d_{\theta } \times N_{x} }$, 
$x\in {\cal X}$, $y \in {\cal Y}$. 
Consequently, 
the dominated convergence theorem and (\ref{1*.1.21}), (\ref{1*.1.25}), 
(\ref{1*.1.27})
imply 
\begin{align} \label{l1.5.5}
	\nabla_{\theta } 
	(\Pi^{n-1} \phi )(\theta, \zeta )
	= &
	E\left(\left. 
	\nabla_{\theta }\left(\phi_{\theta }(G_{\theta }^{0:n}(u,Y_{1:n} ),Y_{n+1} ) \right)
	\right|
	X_{1}=x,Y_{1}=y
	\right)
	\nonumber\\
	= &
	E\left(\left. 
	F_{\theta }(G_{\theta }^{0:n}(u,Y_{1:n} ), H_{\theta }^{0:n}(u,0,Y_{1:n} ), Y_{n+1} ) 
	\right|
	X_{1}=x,Y_{1}=y
	\right)
	\nonumber\\
	= &
	(P^{n-1} F)(\theta, (\zeta, 0 ) )
\end{align}
for any $\theta \in Q$, 
$\zeta = (x,y,u) \in {\cal S}_{\zeta}$, $n>1$
(here, $0$ stands for $d_{\theta } \times N_{x}$ zero matrix). 
As 
$(\Pi^{n} \phi)(\theta, \zeta )$ and 
$(P^{n} F)(\theta, z )$ converge 
(respectively) to 
$\psi(\theta )$ and $g(\theta )$
uniformly in $\theta \in Q$ for each 
$z \in {\cal S}_{z}$, $\zeta \in {\cal S}_{\zeta }$
(due to (\ref{1*.1.25}), (\ref{1*.1.27}), (\ref{l1.5.1}), (\ref{l1.5.3})), 
it follows from 
(\ref{l1.5.5}) that Part (i) is true. 
Part (ii) is then a direct consequence of 
(\ref{1*.1.25}), (\ref{1*.1.27}), 
(\ref{l1.5.1}), (\ref{l1.5.3}). 
\end{IEEEproof}

\begin{lemma} \label{lemma1.6}
Suppose that Assumptions \ref{a2} -- \ref{a4} hold. 
Let $Q\subset \Theta$ be an arbitrary compact set. 
Then, there exists
a real number $C_{7,Q } \in [1,\infty )$
such that 
\begin{align} \label{l1.6.1*}
	\|(P^{n} F)(\theta', z ) - (P^{n} F)(\theta'', z ) \|
	\leq 
	C_{7,Q } \|\theta' - \theta'' \| (1 + \|V\|^{2} ) 
\end{align}
for all $\theta', \theta'' \in Q$, 
$z = (x,y,u,V) \in {\cal S}_{z}$, 
$n\geq 1$. 
\end{lemma}

\begin{IEEEproof}
Let 
Owing to Lemmas \ref{lemma1.1}, \ref{lemma1.3} and \ref{lemma1.4}, we have
\begin{align*}
	&
	\|
	F_{\theta'}(G_{\theta'}^{0:n}(u,y_{1:n} ), H_{\theta'}^{0:n}(u,V,y_{1:n} ), y_{n+1} )
	-
	F_{\theta''}(G_{\theta''}^{0:n}(u,y_{1:n} ), H_{\theta''}^{0:n}(u,V,y_{1:n} ), y_{n+1} )
	\|
	\\
	&
	\begin{aligned}[b]
		\leq &
		C_{1,Q} \psi_{Q}(y_{n+1} ) 
		(1 + \|H_{\theta'}^{0:n}(u,V,y_{1:n} ) \| + \|H_{\theta''}^{0:n}(u,V,y_{1:n} ) \| )
		\\
		& \cdot 
		(\|\theta' - \theta'' \|
		+
		\|G_{\theta'}^{0:n}(u,y_{1:n} ) - G_{\theta''}^{0:n}(u,y_{1:n} ) \| 
		+
		\|H_{\theta'}^{0:n}(u,V,y_{1:n} ) - H_{\theta''}^{0:n}(u,V,y_{1:n} ) \| )
	\end{aligned}
	\\
	&
	\leq
	9 C_{1,Q } C_{4,Q } C_{5,Q } \psi_{Q}(y_{n+1} ) 
	(1 + \|V\| )^{2} \|\theta' - \theta'' \| 
\end{align*}
for all $\theta', \theta'' \in Q$, 
$u \in {\cal P}^{N_{x} }$, 
$V \in \mathbb{R}^{d_{\theta } \times N_{x} }$, $n\geq 1$
and 
any sequence $\boldsymbol y = \{y_{n} \}_{n\geq 1}$
from ${\cal Y}$. 
Consequently, 
\begin{align*}
	\|(\Pi^{n} F)(\theta',z) - (\Pi^{n} F)(\theta'',z) \|
	\leq &
	\begin{aligned}[t]
	E\big( 
	\|
	&
	F_{\theta'}(G_{\theta'}^{0:n}(u,y_{1:n} ), H_{\theta'}^{0:n}(u,V,y_{1:n} ), y_{n+1} )
	\\
	&
	-
	F_{\theta''}(G_{\theta''}^{0:n}(u,y_{1:n} ), H_{\theta''}^{0:n}(u,V,y_{1:n} ), y_{n+1} )
	\|\;
	|X_{1}=x, Y_{1}=y
	\big)
	\end{aligned}
	\\
	\leq &
	9 C_{1,Q } C_{4,Q } C_{5,Q } 
	(1 + \|V\| )^{2} \|\theta' - \theta'' \| 
	\max_{x' \in {\cal X} } 
	\int \psi_{Q}(y') Q(dy'|x') 
\end{align*}
for each $\theta',\theta'' \in Q$, 
$z = (x,y,u,V) \in {\cal S}_{z}$. 
Then, it can be deduced from Assumption \ref{a4}
that there exists a real number 
$C_{7,Q } \in [1,\infty )$ such that 
(\ref{l1.6.1*}) holds for all 
$\theta',\theta'' \in Q$, 
$z = (x,y,u,V) \in {\cal S}_{z}$. 
\end{IEEEproof}

\begin{lemma} \label{lemma1.7}
Suppose that Assumptions \ref{a3} and \ref{a4} hold. 
Let $Q\subset \Theta$ be an arbitrary compact set. 
Then, there exist 
real numbers 
$\delta_{3,Q }, \varepsilon_{4,Q } \in (0,1)$, $C_{8,Q } \in [1, \infty )$
such that the following is true: 
\begin{enumerate}
\item
$\hat{G}_{\eta, \boldsymbol y }^{0:n}(w)$ is analytical 
in $(\eta, w )$ on 
$V_{\delta_{3,Q } }(Q ) \times V_{\delta_{3,Q } }({\cal P}^{N_{x} } )$
for each $n\geq 0$ and 
any sequence $\boldsymbol y = \{y_{n} \}_{n\geq 1}$ from ${\cal Y}$. 
\item
Inequalities   
\begin{align*}
	&
	d(\hat{G}_{\eta,\boldsymbol y}^{0:n}(w), {\cal P}^{N_{x} } )
	\leq \min\{\delta_{Q }, \delta_{1,Q }, \delta_{2,Q } \}, 
	\\
	&
	\|\hat{G}_{\eta, \boldsymbol y }^{0:n}(w') 
	-
	\hat{G}_{\eta, \boldsymbol y }^{0:n}(w'') \| 
	\leq 
	C_{8,Q } \varepsilon_{4,Q }^{n} \|w' - w'' \| 
\end{align*}
hold  
for all $\eta \in V_{\delta_{3,Q} }(Q )$, 
$w,w',w'' \in V_{\delta_{3,Q } }({\cal P}^{N_{x} } )$
and  
any sequence $\boldsymbol y = \{y_{n} \}_{n\geq 1}$ from ${\cal Y}$
($\delta_{Q}$ is specified in Assumption \ref{a4}).  
\end{enumerate} 
\end{lemma}

\begin{IEEEproof}
Let $\boldsymbol y = \{y_{n} \}_{n\geq 1}$ be an arbitrary sequence from ${\cal Y}$. 
Moreover, let 
$k_{Q} = \min\{n\geq 1: C_{3,Q } \varepsilon_{1,Q}^{n} \leq \varepsilon_{1,Q }/2 \}$,
while  
$\tilde{\delta}_{1,Q} = \min\{\delta_{Q}, \delta_{1,Q }, \delta_{2,Q } \}$, 
$\tilde{\delta}_{2,Q } = 4^{-k_{Q}} C_{2,Q }^{-k_{Q } } \tilde{\delta}_{1,Q }$. 

First, we prove by induction (in $k$) that 
\begin{align} \label{l1.7.1}
	d(\hat{G}_{\eta, \boldsymbol y}^{n:n+k}(w), {\cal P}^{N_{x} } )
	\leq 
	(2^{k+1} C_{2,Q }^{k} - 1 ) \tilde{\delta}_{2,Q }
	\leq 
	\tilde{\delta}_{1,Q }
\end{align}
for all $\eta \in V_{\tilde{\delta}_{2,Q } }(Q )$, 
$w \in V_{\tilde{\delta}_{2,Q } }({\cal P}^{N_{x} } )$, 
$n \geq 0$, $0 \leq k \leq k_{Q }$. 
Obviously, (\ref{l1.7.1}) is true when 
$k=0$, $n\geq 0$, $\eta \in V_{\tilde{\delta}_{2,Q } }(Q )$, 
$w \in V_{\tilde{\delta}_{2,Q } }({\cal P}^{N_{x} } )$. 
Suppose now that (\ref{l1.7.1}) holds for 
each 
$\eta \in V_{\tilde{\delta}_{2,Q } }(Q )$, 
$w \in V_{\tilde{\delta}_{2,Q } }({\cal P}^{N_{x} } )$, 
$n \geq 0$
and some $0\leq k < k_{Q }$. 
Then, Lemma \ref{lemma1.2} implies 
\begin{align*}
	\|\hat{G}_{\eta, \boldsymbol y }^{n:n+k+1}(w) - G_{\theta}(u,y_{n+k+1} ) \|
	= &
	\|\hat{G}_{\eta}(\hat{G}_{\eta, \boldsymbol y }^{n:n+k}(w), y_{n+k+1} )
	-
	\hat{G}_{\theta}(u,y_{n+k+1} ) \|
	\\
	\leq &
	C_{2,Q } (\|\eta - \theta \| + \|\hat{G}_{\eta, \boldsymbol y }^{n:n+k}(w) - u \| )
\end{align*}
for any  
$\theta \in Q$, 
$\eta \in V_{\tilde{\delta}_{2,Q } }(Q )$, 
$u \in {\cal P}^{N_{x} }$, 
$w \in V_{\tilde{\delta}_{2,Q } }({\cal P}^{N_{x} } )$, 
$n \geq 0$. Therefore, 
\begin{align*}
	d(\hat{G}_{\eta, \boldsymbol y }^{n:n+k+1}(w), {\cal P}^{N_{x} } )
	\leq &
	C_{2,Q } 
	\left(d(\eta, Q ) 
	+ 
	d(\hat{G}_{\eta, \boldsymbol y }^{n:n+k}(w), {\cal P}^{N_{x} } ) 
	\right)
	\\ 
	\leq &
	2^{k+1} C_{2,Q }^{k+1} \tilde{\delta}_{2,Q }
	\\
	\leq &
	(2^{k+2} C_{2,Q }^{k+1} - 1 ) \tilde{\delta}_{2,Q }
	\leq 
	\tilde{\delta}_{1,Q }
\end{align*}
for any 
$\eta \in V_{\tilde{\delta}_{2,Q } }(Q )$, 
$w \in V_{\tilde{\delta}_{2,Q } }({\cal P}^{N_{x} } )$, 
$n \geq 0$. Hence, (\ref{l1.7.1}) is satisfied for 
all $\eta \in V_{\tilde{\delta}_{2,Q } }(Q )$, 
$w \in V_{\tilde{\delta}_{2,Q } }({\cal P}^{N_{x} } )$, 
$n \geq 0$, $0 \leq k \leq k_{Q }$. 

Let $\tilde{\delta}_{3,Q} = \tilde{\delta}_{2,Q }/2$. 
Since 
$\hat{G}_{\eta, \boldsymbol y }^{n:n}(w)=w$
and 
$\hat{G}_{\eta, \boldsymbol y }^{n:n+k+1}(w) =
\hat{G}_{\eta}(\hat{G}_{\eta, \boldsymbol y }^{n:n+k}(w), y_{n+k+1} )$, 
it can be deduced from Assumption \ref{a4} and (\ref{l1.7.1}) 
that 
$\hat{G}_{\eta, \boldsymbol y}^{n:n+k}(w)$ is analytic 
in $(\eta,w)$ on 
$V_{\tilde{\delta}_{2,Q} }(Q ) \times 
V_{\tilde{\delta}_{2,Q} }({\cal P}^{N_{x} } )$
for each $n\geq 0$, $0 \leq k \leq k_{Q }$
(notice that a composition of two analytic functions is analytic, too). 
Due to Assumption \ref{a4} and (\ref{l1.7.1}), we also have 
\begin{align} \label{l1.7.3}
	\|\hat{G}_{\eta, \boldsymbol y }^{n:n+k+1}(w) \|
	=
	\|\hat{G}_{\eta}(\hat{G}_{\eta, \boldsymbol y }^{n:n+k}(w), y_{n+k+1} ) \|
	\leq 
	K_{Q }
\end{align}
for all 
$\eta \in V_{\tilde{\delta}_{2,Q } }(Q )$, 
$w \in V_{\tilde{\delta}_{2,Q } }({\cal P}^{N_{x} } )$, 
$n\geq 0$, $0\leq k \leq k_{Q}$
($K_{Q}$ is defined in Assumption \ref{a4}). 
As a consequence of Cauchy inequality for analytic functions and (\ref{l1.7.3}), 
there exists a real number 
$\tilde{C}_{1,Q } \in [1, \infty )$
depending exclusively on 
$K_{Q }$, $d_{\theta }$, $N_{x}$
($\tilde{C}_{1,Q }$ can be selected as 
$\tilde{C}_{1,Q } = 
4 (d_{\theta } + N_{x} ) K_{Q }/ \tilde{\delta }_{2,Q }^{2}$) 
such that 
\begin{align*}
	\max\{
	\|\nabla_{(\eta,w)} \hat{G}_{l,\eta, \boldsymbol y}^{n:n+k}(w) \|, 
	\|\nabla_{(\eta,w)}^{2} \hat{G}_{l,\eta, \boldsymbol y}^{n:n+k}(w) \|
	\}
	\leq 
	\tilde{C}_{1,Q }
\end{align*}
for any $\eta \in V_{\tilde{\delta}_{2,Q } }(Q )$, 
$w \in V_{\tilde{\delta}_{2,Q } }({\cal P}^{N_{x} } )$, 
$n\geq 0$, $0\leq k \leq k_{Q}$, $1 \leq l \leq N_{x}$
($\hat{G}_{l,\eta, \boldsymbol y}^{n:n+k}(w)$ denote the 
$l$-th component of 
$\hat{G}_{\eta, \boldsymbol y}^{n:n+k}(w)$). 
Consequently, there exists another real number
$\tilde{C}_{2,Q } \in [1 \infty )$
depending exclusively on 
$K_{Q }$, $d_{\theta }$, $N_{x}$
such that 
\begin{align}
	&\label{l1.7.5}
	\max\{
	\|\hat{G}_{\eta', \boldsymbol y }^{n:n+k}(w') 
	-
	\hat{G}_{\eta'', \boldsymbol y }^{n:n+k}(w'') \|, 
	\|\nabla_{w} \hat{G}_{\eta', \boldsymbol y }^{n:n+k}(w') 
	-
	\nabla_{w} \hat{G}_{\eta'', \boldsymbol y }^{n:n+k}(w'') \|
	\}
	\nonumber \\
	&
	\leq 
	\tilde{C}_{2,Q }
	(\|\eta' - \eta'' \| + \|w' - w'' \| ) 
\end{align}
for each $\eta', \eta'' \in V_{\tilde{\delta}_{3,Q } }(Q )$, 
$w', w'' \in V_{\tilde{\delta}_{3,Q } }({\cal P}^{N_{x} } )$, 
$n\geq 0$, $0\leq k \leq k_{Q}$. 

Let $\tilde{\delta}_{4,Q } = 
\min\{\tilde{\delta}_{3,Q }, 4^{-1} \tilde{C}_{2,Q }^{-1} \varepsilon_{1,Q } \}$. 
Owing to Lemma \ref{lemma1.3} (Part (i)), we have 
\begin{align*}
	\|G_{\theta, \boldsymbol y }^{n:n+k_{Q } }(u') 
	- 
	G_{\theta, \boldsymbol y }^{n:n+k_{Q } }(u'') \|
	\leq 
	C_{3,Q} \varepsilon_{1,Q}^{k_{Q} } \|u' - u'' \| 
	\leq 
	(\varepsilon_{1,Q }/2 ) \|u' - u'' \|
\end{align*}
for all $\theta \in Q$, 
$u',u'' \in [0,\infty )^{N_{x} }\setminus \{0\}$, $n\geq 0$. 
Therefore,
$\|\nabla_{u} G_{\theta, \boldsymbol y }^{n:n+k_{Q } }(w) \| \leq \varepsilon_{1,Q }/2$
for each $\theta \in Q$, 
$u \in [0,\infty )^{N_{x} }\setminus \{0\}$, $n\geq 0$, which, together with (\ref{l1.7.5}) yields 
\begin{align*}
	\|\nabla_{w} \hat{G}_{\eta, \boldsymbol y }^{n: n+k_{Q } }(w) \|
	\leq &
	\|\nabla_{u} G_{\theta, \boldsymbol y }^{n: n+k_{Q } }(u) \|
	+ 
	\|\nabla_{w} \hat{G}_{\eta, \boldsymbol y }^{n: n+k_{Q } }(w) 
	-
	\nabla_{w} \hat{G}_{\theta, \boldsymbol y }^{n: n+k_{Q } }(u) \|
	\\
	\leq &
	\varepsilon_{1,Q }/2 
	+
	\tilde{C}_{2,Q } (\|\theta - \eta \| + \|u - w \| )
\end{align*}
for any 
$\theta \in Q$, 
$\eta \in V_{\tilde{\delta}_{3,Q } }(Q )$, 
$u \in {\cal P}^{N_{x} }$, 
$w \in V_{\tilde{\delta}_{3,Q } }({\cal P}^{N_{x} } )$, 
$n\geq 0$. 
Consequently, 
\begin{align*}
	\|\nabla_{w} \hat{G}_{\eta, \boldsymbol y }^{n: n+k_{Q } }(w) \|
	\leq 
	\varepsilon_{1,Q }/2 
	+
	\tilde{C}_{2,Q } 
	(d(\eta, Q ) + d(w,{\cal P}^{N_{x} } ) )
	\leq 
	\varepsilon_{1,Q }
\end{align*}
for each 
$\eta \in V_{\tilde{\delta}_{4,Q } }(Q )$, 
$w \in V_{\tilde{\delta}_{4,Q } }({\cal P}^{N_{x} } )$, 
$n\geq 0$. 
Thus, 
\begin{align} \label{l1.7.7}
	\|\hat{G}_{\eta, \boldsymbol y }^{n: n+k_{Q } }(w') 
	- 
	\hat{G}_{\eta, \boldsymbol y }^{n: n+k_{Q } }(w'') \|
	\leq 
	\int_{0}^{1} 
	\|\nabla_{w} \hat{G}_{\eta, \boldsymbol y }^{n: n+k_{Q } }(tw' + (1-t)w'' ) \|
	\|w' - w'' \| dt 
	\leq 
	\varepsilon_{1,Q } \|w' - w'' \|
\end{align}
for all  
$\eta \in V_{\tilde{\delta}_{4,Q } }(Q )$, 
$w \in V_{\tilde{\delta}_{4,Q } }({\cal P}^{N_{x} } )$, 
$n\geq 0$. 

Let $\tilde{\delta}_{5,Q } = 
(1 - \varepsilon_{1,Q } ) \tilde{\delta}_{4,Q } \tilde{C}_{2,Q }^{-1}$. 
Now, we prove by induction (in $i$) that 
\begin{align}\label{l1.7.9}
	d(\hat{G}_{\eta, \boldsymbol y }^{0:ik_{Q } }(w), {\cal P}^{N_{x} } )
	\leq 
	\tilde{\delta}_{4,Q }
\end{align}
for each 
$\eta \in V_{\tilde{\delta}_{5,Q } }(Q )$, 
$w \in V_{\tilde{\delta}_{4,Q } }({\cal P}^{N_{x} } )$, 
$i\geq 0$. 
Obviously, (\ref{l1.7.9}) is true when 
$i=0$, 
$\eta \in V_{\tilde{\delta}_{5,Q } }(Q )$, 
$w \in V_{\tilde{\delta}_{4,Q } }({\cal P}^{N_{x} } )$. 
Suppose that (\ref{l1.7.9}) holds 
for all 
$\eta \in V_{\tilde{\delta}_{5,Q } }(Q )$, 
$w \in V_{\tilde{\delta}_{4,Q } }({\cal P}^{N_{x} } )$
and some $i\geq 0$. 
Then, (\ref{l1.7.5}), (\ref{l1.7.7}) imply 
\begin{align*}
	\|\hat{G}_{\eta, \boldsymbol y }^{0:(i+1)k_{Q } }(w)
	-
	G_{\theta, \boldsymbol y }^{ik_{Q }:(i+1)k_{Q} }(u) \|
	\leq &
	\|\hat{G}_{\eta, \boldsymbol y }^{ik_{Q }:(i+1)k_{Q} }
	(\hat{G}_{\eta, \boldsymbol y }^{0:ik_{Q } }(w) ) 
	-
	\hat{G}_{\eta, \boldsymbol y }^{ik_{Q }:(i+1)k_{Q} }(u) \|
	\\
	&
	+
	\|\hat{G}_{\eta, \boldsymbol y }^{ik_{Q }:(i+1)k_{Q} }(u) 
	-
	\hat{G}_{\theta, \boldsymbol y }^{ik_{Q }:(i+1)k_{Q} }(u) \|
	\\
	\leq &
	\varepsilon_{1,Q } 
	\|\hat{G}_{\eta, \boldsymbol y }^{0:ik_{Q } }(w) - u \|
	+
	\tilde{C}_{2,Q } \|\theta - \eta \|
\end{align*}
for any 
$\theta \in Q$, 
$\eta \in V_{\tilde{\delta}_{5,Q } }(Q )$, 
$u \in {\cal P}^{N_{x} }$,
$w \in V_{\tilde{\delta}_{4,Q } }({\cal P}^{N_{x} } )$. 
Therefore, 
\begin{align*}
	d(\hat{G}_{\eta, \boldsymbol y }^{0:(i+1)k_{Q } }(w), {\cal P}^{N_{x} } )
	\leq 
	\varepsilon_{1,Q } 
	d(\hat{G}_{\eta, \boldsymbol y }^{0:ik_{Q } }(w), {\cal P}^{N_{x} } )
	+
	\tilde{C}_{2,Q } d(\eta, Q )
	\leq 
	\varepsilon_{1,Q } \tilde{\delta}_{4,Q } 
	+ 
	\tilde{C}_{2,Q } \tilde{\delta}_{5,Q }
	=
	\tilde{\delta}_{4,Q }
\end{align*}
for each 
$\eta \in V_{\tilde{\delta}_{5,Q } }(Q )$, 
$w \in V_{\tilde{\delta}_{4,Q } }({\cal P}^{N_{x} } )$. 
Hence, (\ref{l1.7.9}) holds 
for all  
$\eta \in V_{\tilde{\delta}_{5,Q } }(Q )$, 
$w \in V_{\tilde{\delta}_{4,Q } }({\cal P}^{N_{x} } )$, 
$i\geq 0$. 

Let $\delta_{3,Q } = \min\{\tilde{\delta}_{4,Q }, \tilde{\delta}_{5,Q } \}$. 
As 
$\hat{G}_{\eta, \boldsymbol y}^{0:0}(w) = w$
and 
$\hat{G}_{\eta, \boldsymbol y }^{0:(i+1)k_{Q } }(w) = 
\hat{G}_{\eta, \boldsymbol y}^{ik_{Q }: (i+1)k_{Q } }
(\hat{G}_{\eta, \boldsymbol y }^{0:ik_{Q } }(w) )$, 
it can be deduced from (\ref{l1.7.9}) that 
$\hat{G}_{\eta, \boldsymbol y }^{0:ik_{Q } }(w)$ is analytical 
in $(\eta, w )$ on 
$V_{\tilde{\delta}_{5,Q } }(Q ) 
\times V_{\tilde{\delta}_{4,Q } }({\cal P}^{N_{x} } )$
for each $i\geq 0$ 
(notice that 
$\hat{G}_{\eta, \boldsymbol y}^{ik_{Q }: (i+1)k_{Q } }(w)$ is analytic in 
$(\eta, w )$ on 
$V_{\tilde{\delta}_{5,Q } }(Q ) 
\times V_{\tilde{\delta}_{4,Q } }({\cal P}^{N_{x} } )$
for any $i\geq 0$). 
Since 
$\hat{G}_{\eta, \boldsymbol y}^{0:n}(w) =
\hat{G}_{\eta, \boldsymbol y}^{ik_{Q }:n}
(\hat{G}_{\eta, \boldsymbol y}^{0:ik_{Q } }(w) )$
for 
$i=\lfloor n/k_{Q } \rfloor$, 
we conclude from (\ref{l1.7.9}) that 
$\hat{G}_{\eta, \boldsymbol y}^{0:n}(w)$ is analytical 
in $(\eta, w )$ on 
$V_{\tilde{\delta}_{5,Q } }(Q) 
\times V_{\tilde{\delta}_{4,Q } }({\cal P}^{N_{x} } )
\supseteq 
V_{\delta_{3,Q } }(Q ) \times V_{\delta_{3,Q } }({\cal P}^{N_{x} } )$
for all $n\geq 0$ 
(notice that 
$\hat{G}_{\eta, \boldsymbol y}^{ik_{Q }:ik_{Q }+j}(w)$ is analytical 
in $(\eta, w )$ on 
$V_{\tilde{\delta}_{5,Q } }(Q ) 
\times V_{\tilde{\delta}_{4,Q } }({\cal P}^{N_{x} } )$
for any $i\geq 0$, $0\leq j \leq k_{Q }$). 
On the other side, (\ref{l1.7.1}), (\ref{l1.7.9}) yield 
\begin{align} \label{l1.7.11}
	d(\hat{G}_{\eta, \boldsymbol y}^{0:n}(w), {\cal P}^{N_{x} } )
	=
	d(\hat{G}_{\eta, \boldsymbol y}^{ik_{Q }:n}
	(\hat{G}_{\eta, \boldsymbol y}^{0:ik_{Q } }(w) ), {\cal P}^{N_{x} } )
	\leq 
	\tilde{\delta}_{1,Q } 
	=
	\min\{\delta_{Q }, \delta_{1,Q }, \delta_{2,Q } \}
\end{align}
for all 
$\eta \in V_{\tilde{\delta}_{5,Q } }(Q ) \supseteq
V_{\delta_{3,Q } }(Q )$, 
$w \in V_{\tilde{\delta}_{5,Q } }({\cal P}^{N_{x} } ) \supseteq
V_{\delta_{3,Q } }({\cal P}^{N_{x} } )$, $n\geq 0$
and 
$i=\lfloor n/k_{Q } \rfloor$. 

Let $\varepsilon_{4,Q } = \varepsilon_{1,Q }^{1/k_{Q } }$, 
$C_{8,Q } = \tilde{C}_{2,Q } \varepsilon_{1,Q }^{-1}$. 
Owing to (\ref{l1.7.7}), (\ref{l1.7.9}), we have 
\begin{align*}
	\|\hat{G}_{\eta,\boldsymbol y }^{0:(i+1)k_{Q } }(w') 
	- 
	\hat{G}_{\eta,\boldsymbol y }^{0:(i+1)k_{Q } }(w'')  \|
	= &
	\|\hat{G}_{\eta,\boldsymbol y }^{ik_{Q }:(i+1)k_{Q } }
	(\hat{G}_{\eta,\boldsymbol y }^{0:ik_{Q } }(w') ) 
	- 
	\hat{G}_{\eta,\boldsymbol y }^{ik_{Q }:(i+1)k_{Q } }
	(\hat{G}_{\eta,\boldsymbol y }^{0:ik_{Q } }(w'') ) \|
	\\
	\leq & 
	\varepsilon_{1,Q } 
	\|\hat{G}_{\eta,\boldsymbol y }^{0:ik_{Q } }(w') 
	- 
	\hat{G}_{\eta,\boldsymbol y }^{0:ik_{Q } }(w'')  \|
\end{align*}
for any $\eta \in V_{\tilde{\delta}_{5,Q } }(Q )$, 
$w',w'' \in V_{\tilde{\delta}_{4,Q } }({\cal P}^{N_{x} } )$, $i\geq 0$. 
Therefore, 
\begin{align*}
	\|\hat{G}_{\eta,\boldsymbol y }^{0:ik_{Q } }(w') 
	- 
	\hat{G}_{\eta,\boldsymbol y }^{0:ik_{Q } }(w'')  \|
	\leq 
	\varepsilon_{1,Q }^{i} \|w' - w'' \|
\end{align*}
for each $\eta \in V_{\tilde{\delta}_{5,Q } }(Q )$, 
$w',w'' \in V_{\tilde{\delta}_{4,Q } }({\cal P}^{N_{x} } )$, $i\geq 0$. 
Consequently, (\ref{l1.7.5}), (\ref{l1.7.9}) yield
\begin{align*}
	\|\hat{G}_{\eta,\boldsymbol y }^{0:n}(w') 
	- 
	\hat{G}_{\eta,\boldsymbol y }^{0:n}(w'') \|
	= &
	\|\hat{G}_{\eta,\boldsymbol y }^{ik_{Q }:n }
	(\hat{G}_{\eta,\boldsymbol y }^{0:ik_{Q } }(w') ) 
	- 
	\hat{G}_{\eta,\boldsymbol y }^{ik_{Q }:n }
	(\hat{G}_{\eta,\boldsymbol y }^{0:ik_{Q } }(w'') ) \|
	\\
	\leq &
	\tilde{C}_{2,Q } 
	\|\hat{G}_{\eta,\boldsymbol y }^{0:ik_{Q } }(w') 
	- 
	\hat{G}_{\eta,\boldsymbol y }^{0:ik_{Q } }(w'')  \|
	\\
	\leq &
	\tilde{C}_{2,Q } \varepsilon_{1,Q }^{i} \|w' - w'' \|
	\\
	\leq &
	C_{8,Q } \varepsilon_{4,Q }^{n} \|w' - w'' \|
\end{align*}
for each 
$\eta \in V_{\tilde{\delta}_{5,Q } }(Q ) \supseteq
V_{\delta_{3,Q } }(Q )$, 
$w',w'' \in V_{\tilde{\delta}_{4,Q } }({\cal P}^{N_{x} } ) \supseteq
V_{\delta_{3,Q } }({\cal P}^{N_{x} } )$, $n\geq 0$, 
$i = \lfloor n/k_{Q } \rfloor$
(notice that 
$\tilde{C}_{2,Q } \varepsilon_{1,Q }^{i} 
= \tilde{C}_{2,Q } \varepsilon_{4,Q}^{-(n-ik_{Q } ) } \varepsilon_{4,Q }^{n}
\leq C_{8,Q } \varepsilon_{4,Q }^{n}$). 
Then, it is clear that $\delta_{3,Q }$, $\varepsilon_{4,Q }$, 
$C_{8,Q }$ meet the requirements of the lemma. 
\end{IEEEproof}

\subsection{Analyticity}\label{subsection1.2*} 

In this subsection, using the results of the Subsection \ref{subsection1.1*}
(Lemma \ref{lemma1.7}), 
the analyticity of the objective function $f(\cdot )$ 
is shown and Theorem \ref{theorem1} is proved. 
The proof is based on the analytic continuation techniques 
and the methods developed in 
\cite{han&marcus1}. 

\begin{IEEEproof}[Proof of Theorem \ref{theorem1}]
Let 
\begin{align*} 
	\hat{\psi}_{\eta}^{n}(w,x)
	=
	E\left(\left.
	\hat{\phi}_{\eta}(\hat{G}_{\eta}^{0:n}(w,Y_{1:n} ), Y_{n+1} )
	\right|X_{1} = x
	\right)
\end{align*}
for $\eta\in \mathbb{C}^{d_{\theta } }$, $w\in \mathbb{C}^{N_{x} }$, $x\in {\cal X}$, $n\geq 1$.
Then, using (\ref{1*.1.7}), it is straightforward to verify 
\begin{align} \label{t1.1}
	\hat{\psi}_{\eta}^{n+1}(w,x) 
	= &
	E\left(\left.
	E\left(\left.
	\hat{\phi}_{\eta}
	\big(\hat{G}_{\eta}^{0:n}(\hat{G}_{\eta}(w,Y_{1} ), Y_{2:n+1} ), Y_{n+2} \big)
	\right|
	X_{1},X_{2},Y_{1}
	\right)
	\right|
	X_{1}=x
	\right)
	\nonumber \\
	= &
	E(
	\hat{\psi}_{\eta}^{n}(\hat{G}_{\eta}(w,Y_{1}), X_{2} )
	|X_{1}=x)
\end{align}
for each 
$\eta \in \mathbb{C}^{d_{\theta } }$, $w \in \mathbb{C}^{N_{x} }$, $x\in {\cal X}$, 
$n \geq 0$. 
It is also easy to show
\begin{align}\label{t1.3}
	&
	\hat{\psi}_{\eta}^{n}(w',x') - \hat{\psi}_{\eta}^{n}(w'',x'') 
	\nonumber\\
	&
	\begin{aligned}[b]
		= &
		E\left(\left.
		\hat{\phi}_{\eta}(\hat{G}_{\eta}^{0:n}(w',Y_{1:n} ), Y_{n+1} )
		-
		\hat{\phi}_{\eta}(\hat{G}_{\eta}^{0:n}(e_{0},Y_{1:n} ), Y_{n+1} )
		\right|
		X_{1}=x' 
		\right)
		\\
		&
		-
		E\left(\left.
		\hat{\phi}_{\eta}(\hat{G}_{\eta}^{0:n}(w'',Y_{1:n} ), Y_{n+1} )
		-
		\hat{\phi}_{\eta}(\hat{G}_{\eta}^{0:n}(e_{0},Y_{1:n} ), Y_{n+1} )
		\right|
		X_{1}=x'' 
		\right)
		\\
		&
		+
		\begin{aligned}[t]
		\sum_{k=1}^{n-1} 
		\sum_{x\in {\cal X} }
		&
		E\left(\left.
		\hat{\phi}_{\eta}(\hat{G}_{\eta}^{0:n-k+1}(e_{0},Y_{k:n} ), Y_{n+1} )
		-
		\hat{\phi}_{\eta}(\hat{G}_{\eta}^{0:n-k}(e_{0},Y_{k+1:n} ), Y_{n+1} )
		\right|
		X_{k}=x  
		\right)
		\\
		&
		\cdot
		(p^{k-1}(x|x') - \pi(x) )
		\end{aligned}
		\\
		&
		-
		\begin{aligned}[t]
		\sum_{k=1}^{n-1} 
		\sum_{x\in {\cal X} }
		&
		E\left(\left.
		\hat{\phi}_{\eta}(\hat{G}_{\eta}^{0:n-k+1}(e_{0},Y_{k:n} ), Y_{n+1} )
		-
		\hat{\phi}_{\eta}(\hat{G}_{\eta}^{0:n-k}(e_{0},Y_{k+1:n} ), Y_{n+1} )
		\right|
		X_{k}=x  
		\right)
		\\
		&
		\cdot
		(p^{k-1}(x|x'') - \pi(x) )
		\end{aligned}
		\\
		&
		+
		\sum_{x\in {\cal X} }
		E(\hat{\phi}_{\eta}(\hat{G}_{\eta}(e_{0},Y_{n} ),Y_{n+1} )|X_{n} = x )
		(p^{n-1}(x|x') - \pi(x) )
		\\
		&
		-
		\sum_{x\in {\cal X} }
		E(\hat{\phi}_{\eta}(\hat{G}_{\eta}(e_{0},Y_{n} ),Y_{n+1} )|X_{n} = x )
		(p^{n-1}(x|x') - \pi(x) )
	\end{aligned}
\end{align}
for all 
$\eta \in \mathbb{C}^{d_{\theta } }$, 
$w',w'' \in \mathbb{C}^{N_{x} }$, 
$x',x'' \in {\cal X}$, 
$n\geq 1$, 
where 
$e_{0} = [1 \cdots 1]^{T}/N_{x} \in \mathbb{R}^{N_{x} }$
and
$p^{k-1}(x'|x) = P(X_{k} = x'|X_{1} = x)$,  
$\pi(x) = \lim_{k\rightarrow \infty } P(X_{k}=x)$. 
On the other side, 
Assumption \ref{a2} implies that $\pi(\cdot )$ is well-defined 
and that 
there exist real numbers 
$\tilde{\varepsilon} \in (0,1)$, $\tilde{C} \in [1,\infty )$  
such that 
\begin{align}\label{t1.5}
	|p^{n}(x'|x) - \pi(x') |
	\leq 
	\tilde{C} \tilde{\varepsilon}^{n}
\end{align}
for each $x,x'\in {\cal X}$, $n\geq 0$. 

Let $Q\subset \Theta$ be an arbitrary compact set, while 
$\tilde{\delta}_{1,Q } = 
\min\{\delta_{Q}, \delta_{1,Q }, \delta_{2,Q }, \delta_{3,Q } \}$, 
$\tilde{\delta}_{2,Q} = \tilde{\delta}_{1,Q }/2$. 
Owing to Assumption \ref{a4} and Lemma \ref{lemma1.7}, 
$\hat{\phi}_{\eta}(\hat{G}_{\eta}^{0:n}(w,y_{1:n} ), y_{n+1} )$ 
is analytic in 
$(\eta,w)$ on 
$V_{\tilde{\delta}_{1,Q } }(Q ) \times 
V_{\tilde{\delta}_{1,Q } }({\cal P}^{N_{x} } )$
for each $n\geq 0$ and any sequence 
$\boldsymbol y = \{y_{n} \}_{n\geq 1}$ from ${\cal Y}$. 
Due to Assumption \ref{a4} and Lemma \ref{lemma1.7}, 
we also have 
\begin{align*}
	|\hat{\phi}_{\eta}(\hat{G}_{\eta}^{0:n}(w,y_{1:n} ), y_{n+1} ) |
	\leq 
	\psi_{Q}(y_{n+1} )
\end{align*}
for all 
$\eta \in V_{\tilde{\delta}_{1,Q } }(Q )$, 
$w \in V_{\tilde{\delta}_{1,Q } }({\cal P}^{N_{x} } )$, $n\geq 0$
and any sequence 
$\boldsymbol y = \{y_{n} \}_{n\geq 1}$ from ${\cal Y}$. 
Consequently, Cauchy inequality for analytic functions implies that there exists a 
real number $\tilde{C}_{1,Q } \in [1, \infty )$ such that 
\begin{align} \label{t1.7}
	\|\nabla_{\eta } 
	\hat{\phi}_{\eta}(\hat{G}_{\eta}^{0:n}(w,y_{1:n} ), y_{n+1} ) \|
	\leq 
	\tilde{C}_{1,Q } 
	\psi_{Q}(y_{n+1} )
\end{align}
for each  
$\eta \in V_{\tilde{\delta}_{2,Q } }(Q )$, 
$w \in V_{\tilde{\delta}_{2,Q } }({\cal P}^{N_{x} } )$, $n\geq 0$
and any sequence 
$\boldsymbol y = \{y_{n} \}_{n\geq 1}$ from ${\cal Y}$. 
Since 
\begin{align} \label{t1.9}
	E(\psi_{Q }(Y_{n+1} )|X_{1}=x )
	\leq 
	\max_{x'\in {\cal X} }
	\int \psi_{Q }(y') Q(dy'|x') 
	<
	\infty
\end{align}
for all $x\in {\cal X}$, $n \geq 0$, 
it follows from the dominated convergence theorem 
and (\ref{t1.7}) that 
$\hat{\psi}_{\eta }^{n}(w,x)$ is 
differentiable
(and thus, analytic) in $\eta$
on 
$V_{\tilde{\delta}_{2,Q } }(Q )$ 
for any $w \in V_{\tilde{\delta}_{2,Q } }({\cal P}^{N_{x} } )$, 
$n\geq 0$. 

Let $\tilde{\varepsilon}_{Q } = \max\{\varepsilon_{4,Q }, \tilde{\varepsilon} \}$. 
Due to Lemmas \ref{lemma1.1} and \ref{lemma1.7}, we have 
\begin{align}
	&\label{t1.21}
	|\hat{\phi}_{\eta}(\hat{G}_{\eta}^{0:n}(w',y_{1:n} ), y_{n+1} ) 
	-
	\hat{\phi}_{\eta}(\hat{G}_{\eta}^{0:n}(w'',y_{1:n} ), y_{n+1} ) |
	\leq 
	C_{1,Q } C_{8,Q } \varepsilon_{4,Q }^{n} \psi_{Q}(y_{n+1} )
	\|w' - w'' \|, 
	\\
	&\label{t1.23}
	|\hat{\phi}_{\eta}(\hat{G}_{\eta}^{0:n-k+1}(w,y_{k:n} ), y_{n+1} ) 
	-
	\hat{\phi}_{\eta}(\hat{G}_{\eta}^{0:n-k}(w,y_{k+1:n} ), y_{n+1} ) |
	\nonumber\\
	& \;\;\; 
	\leq 
	C_{1,Q }
	\psi_{Q}(y_{n+1} ) 
	\|\hat{G}_{\eta }^{0:n-k}(\hat{G}_{\eta}(w, y_{k} ), y_{k+1:n} )
	-
	\hat{G}_{\eta}^{0:n-k}(w,y_{k+1:n} ) 	\|
	\nonumber\\
	& \;\;\; 
	\leq 
	C_{1,Q } C_{8,Q } \varepsilon_{4,Q }^{n-k} \psi_{Q}(y_{n+1} ) 
	\|\hat{G}_{\eta}(w, y_{k} ) - w \|
\end{align}
for each 
$\eta \in V_{\tilde{\delta}_{2,Q } }(Q )$, 
$w,w',w'' \in V_{\tilde{\delta}_{2,Q } }({\cal P}^{N_{x} } )$, 
$n\geq 1$, $0 < k \leq n$ 
and any sequence 
$\boldsymbol y = \{y_{n} \}_{n\geq 1}$ from 
${\cal Y}$. 
Using (\ref{t1.5}), (\ref{t1.9}) -- (\ref{t1.23}), 
we deduce that there exists a real number 
$\tilde{C}_{2,Q } \in [1, \infty )$ such that 
the absolute value of the each term on right-hand side of 
(\ref{t1.3}) is bounded by 
$\tilde{C}_{2,Q } \tilde{\varepsilon}_{Q }^{n}$
for any 
$\eta \in V_{\tilde{\delta}_{2,Q } }(Q )$, 
$w',w'' \in V_{\tilde{\delta}_{2,Q } }({\cal P}^{N_{x} } )$, 
$x,x' \in {\cal X}$, 
$n\geq 1$. 
Therefore, 
\begin{align}\label{t1.25}
	|\hat{\psi}_{\eta}^{n}(w',x') - \hat{\psi}_{\eta}^{n}(w'',x'') |
	\leq 
	2\tilde{C}_{2,Q} \tilde{\varepsilon}_{Q}^{n} (n+1)
\end{align}
for all 
$\eta \in V_{\tilde{\delta}_{2,Q } }(Q )$, 
$w',w'' \in V_{\tilde{\delta}_{2,Q } }({\cal P}^{N_{x} } )$, 
$x',x'' \in {\cal X}$, 
$n\geq 1$. 
Consequently, (\ref{t1.1}) yields 
\begin{align} \label{t1.27}
	|\hat{\psi}_{\eta}^{n+1}(w,x) - \hat{\psi}_{\eta}^{n}(w,x) |
	\leq 
	E\left(\left. 
	|\hat{\psi}_{\eta}^{n}(\hat{G}_{\eta}(w,Y_{1}),X_{2} ) - \hat{\psi}_{\eta}^{n}(w,x) |
	\right|
	X_{1} = x 
	\right)
	\leq 
	2\tilde{C}_{2,Q} \tilde{\varepsilon}_{Q}^{n} (n+1)
\end{align}
for each $\eta \in V_{\tilde{\delta}_{2,Q } }(Q )$, 
$w \in V_{\tilde{\delta}_{2,Q } }({\cal P}^{N_{x} } )$, $x \in {\cal X}$, 
$n\geq 1$. 
Owing to (\ref{t1.25}), (\ref{t1.27}), 
there exists a function 
$\hat{\psi}: \mathbb{C}^{d_{\theta } } \rightarrow \mathbb{C}$ 
such that 
$\hat{\psi}_{\eta}^{n}(w,x)$ converges to 
$\hat{\psi}(\eta )$ uniformly in 
$(\eta, w, x ) \in 
V_{\tilde{\delta}_{2,Q } }(Q ) \times 
V_{\tilde{\delta}_{2,Q } }({\cal P}^{N_{x} } ) \times {\cal X}$. 
As the uniform limit of analytic functions is 
also an analytic function 
(see \cite[Theorem 2.4.1]{taylor}), 
$\hat{\psi}(\cdot )$ is 
analytic on $V_{\tilde{\delta}_{2,Q} }(Q)$. 
On the other side, since 
\begin{align*}
	\hat{\phi}_{\theta }^{n}(u,x) 
	=
	E\left(
	\left.\phi_{\theta}(G_{\theta}^{0:n}(u,Y_{1:n} ), Y_{n+1} )
	\right| X_{1} = x 
	\right)
	=
	E\left(\left.
	(\Pi^{n-1} \phi )(\theta, (x,Y_{1},u) ) 
	\right|
	X_{1} = x 
	\right)
\end{align*}
for all $\theta \in \Theta$, $u \in {\cal P}^{N_{x} }$, 
$x \in {\cal X}$, $n \geq 1$, 
Lemma \ref{lemma1.5} implies 
$f(\theta ) = \hat{\psi}(\theta )$ for any $\theta \in Q$. 
Then, it is clear that Part (i) is true, 
while  
Part (ii) follows from the Lojasiewicz inequality 
(see e.g., \cite{kurdyka}, \cite{lojasiewicz1}, \cite{lojasiewicz2})
and the analyticity of $f(\cdot )$. 
\end{IEEEproof}

As a direct consequence of 
\cite[Theorem \L I, Page 775]{kurdyka}
and
Theorem \ref{theorem1}, we have the following corollary: 

\begin{corollary}\label{corollary1}
Let Assumptions \ref{a2} -- \ref{a4} hold. 
Then, for any compact set 
$Q \subset \Theta$ and real number 
$a \in f(Q)$, 
there exist real numbers 
$\delta_{Q,a} \in (0,1)$, 
$\mu_{Q,a} \in (1,2]$, 
$M_{Q,a} \in [1,\infty )$ such that 
\begin{align*}
	|f(\theta ) - a |
	\leq 
	M_{Q,a} \|\nabla f(\theta ) \|^{\mu_{Q,a} }
\end{align*}
for all $\theta\in Q$
satisfying $|f(\theta ) - a |\leq \delta_{Q,a}$.
\end{corollary}

\begin{remark}
Obviously, 
if $Q\subseteq \{\theta'\in \mathbb{R}^{d_{\theta } }: 
\|\theta' - \theta\| \leq \delta_{\theta } \}$
and $a=f(\theta )$ for some $\theta \in \mathbb{R}^{d_{\theta } }$,
then $\mu_{Q,a}$ and $M_{Q,a}$
can be selected as $\mu_{Q,a} = \mu_{\theta }$ and 
$M_{Q,a} = M_{\theta }$
($\delta_{\theta }$, $\mu_{\theta }$, 
$M_{\theta }$ are specified in the statement of Theorem \ref{theorem1}).
\end{remark}

\subsection{Decomposition of Algorithm (\ref{1.1}) -- (\ref{1.5})}\label{subsection1.3*} 

Relying on the results 
of Subsection \ref{subsection1.1*} (Lemmas \ref{lemma1.1} -- \ref{lemma1.6}), 
equivalent representations of recursion (\ref{1.1}) -- (\ref{1.5}) 
and their asymptotic properties are analyzed in this subsection. 
The analysis is based on the techniques developed in 
\cite[Part II]{benveniste&metivier&priouret}. 
The results of this subsection are a crucial prerequisite 
for the analysis carried out in the next subsection. 

In this subsection, the following notation is used. 
For $n\geq 0$, let 
$Z_{n+1} = (X_{n+1}, Y_{n+1}, U_{n}, V_{n} )$, 
while 
\begin{align*}
	&
	\xi_{n}
	=
	F(\theta_{n}, Z_{n+1} )
	-
	\nabla f(\theta_{n} ), 
	\\
	&
	\phi'_{n}
	=
	\alpha_{n} 
	(\nabla f(\theta_{n} ) )^{T} \xi_{n}, 
	\\
	&
	\phi''_{n} 
	=
	\int_{0}^{1}
	(\nabla f(\theta_{n} + t (\theta_{n+1} - \theta_{n} ) ) 
	- 
	\nabla f(\theta_{n} ) )^{T}
	(\theta_{n+1} - \theta_{n} )
	dt
\end{align*}
and 
$\phi_{n} = \phi'_{n} + \phi''_{n}$
($F(\cdot, \cdot )$ is defined in the beginning of the previous subsection). 
Then, 
algorithm (\ref{1.1}) -- (\ref{1.5}) admits the following representations: 
\begin{align*}
	\theta_{n+1} 
	= &
	\theta_{n} + \alpha_{n} F(\theta_{n}, Z_{n+1} )
	\\
	= &
	\theta_{n} + \alpha_{n} (\nabla f(\theta_{n} ) + \xi_{n} ), 
	\;\;\; n \geq 0. 
\end{align*}
Moreover, we have 
\begin{align*}
	f(\theta_{n+1} )
	=
	f(\theta_{n} ) 
	+
	\alpha_{n} \|\nabla f(\theta_{n} ) \|^{2} 
	+
	\phi_{n}
\end{align*}
for $n\geq 0$. 
We also conclude 
\begin{align*}
	P(Z_{n+1} \in B|\theta_{0},Z_{0},\dots,\theta_{n},Z_{n} )
	=
	P_{\theta_{n} }(Z_{n}, B )
\end{align*}
w.p.1 for $n\geq 0$ and any Borel-measurable set $B \subseteq {\cal S}_{z}$ 
($P_{\theta }(\cdot, \cdot )$ is also introduced in the beginning of the previous subsection).

\begin{lemma} \label{lemma2.1}
Suppose that Assumptions \ref{a2} -- \ref{a4} hold. 
Then, there exists a Borel-measurable function 
$\Phi: \Theta \times {\cal S}_{z} \rightarrow \mathbb{R}^{d_{\theta } }$ 
with the following properties: 
\begin{enumerate}
\item
$\Phi(\theta, \cdot )$ is integrable with respect to 
$P_{\theta}(z,\cdot )$ and 
\begin{align} \label{l2.1.1*}
	F(\theta, z ) 
	-
	\nabla f(\theta )
	=
	\Phi(\theta, z ) 
	-
	(P \Phi )(\theta, z )
\end{align}
for all $\theta \in \Theta$, $z \in {\cal S}_{z}$. 
\item
For any compact set $Q\subset\Theta$ and a real number $s \in (0,1)$, 
there exists 
a Borel-measurable function 
$\varphi_{Q,s }: {\cal S}_{z} \rightarrow [1, \infty )$
such that 
\begin{align}
	& \label{l2.1.3*}
	\max\{
	\|F(\theta, z ) \|, 
	\|\Phi(\theta, z ) \|, 
	\|(P\Phi )(\theta, z ) \|
	\}
	\leq 
	\varphi_{Q,s}(z), 
	\\
	& \label{l2.1.5*}
	\|(P\Phi )(\theta', z ) - (P\Phi)(\theta'', z ) \|
	\leq 
	\varphi_{Q,s}(z) \|\theta' - \theta'' \|^{s}, 
	\\
	& \label{l2.1.7*}
	\sup_{n\geq 0}
	E\left(\left. 
	\varphi_{Q,s}^{2}(Z_{n} ) 
	I_{ \{\tau_{Q} \geq n \} }
	\right|
	Z_{0} = z
	\right)
	< \infty
\end{align}
for all $\theta, \theta', \theta'' \in Q$, 
$z \in {\cal S}_{z}$, where 
\begin{align*}
	\tau_{Q}
	=
	\inf\{n\geq 0: \theta_{n} \not\in Q \}. 
\end{align*}
\end{enumerate}
\end{lemma}

\begin{IEEEproof}
Let $Q\subseteq\Theta$ be an arbitrary compact set. 
Owing to Lemmas \ref{lemma1.1} and \ref{lemma1.5}, 
there exists a real number 
$\tilde{C}_{1,Q} \in [1,\infty )$ such that 
\begin{align}\label{l2.1.1}
	\sum_{k=0}^{\infty } 
	\|(P^{k} F)(\theta, z ) - \nabla f(\theta ) \|
	\leq 
	\tilde{C}_{1,Q} \psi_{Q}(y) 
	(1 + \|V\|^{2} )
\end{align}
for all $\theta \in Q$, 
$z=(x,y,u,v) \in {\cal S}_{z}$
($(P^{0} F)(\theta,z)$ stands for $F(\theta, z)$). 
Consequently, 
$\sum_{k=0}^{\infty } ((P^{k} F)(\theta, z ) - \nabla f(\theta ) )$
if well-defined and finite for each 
$\theta \in Q$, 
$z \in {\cal S}_{z}$. 
We also have 
\begin{align*}
	&
	\left\|
	\sum_{k=1}^{\infty} ((P^{k} F)(\theta', z ) - \nabla f(\theta' ) )
	-
	\sum_{k=1}^{\infty} ((P^{k} F)(\theta'', z ) - \nabla f(\theta'' ) )
	\right\|
	\\
	&
	\begin{aligned}[b]
	\leq &
	\sum_{k=1}^{n} 
	\left\|(P^{k} F)(\theta', z ) - (P^{k} F)(\theta'', z ) \right\|
	+
	n
	\left\|\nabla f(\theta' ) - \nabla f(\theta'' ) \right\|
	\\
	&
	+
	\sum_{k=n+1}^{\infty} 
	\left\|(P^{k} F)(\theta', z ) - \nabla f(\theta' ) \right\|
	+
	\sum_{k=n+1}^{\infty} 
	\left\|(P^{k} F)(\theta'', z ) - \nabla f(\theta'' ) \right\|
	\end{aligned}
\end{align*}
for each $\theta', \theta'' \in \Theta$, $z \in {\cal S}_{z}$, 
$n\geq 1$. 
Then, using Lemmas \ref{lemma1.5} and \ref{lemma1.6}, 
it can be deduced 
that there exist real numbers 
$\tilde{\varepsilon}_{Q} \in (0,1)$, 
$\tilde{C}_{2,Q} \in [1, \infty )$ such that 
\begin{align}\label{l2.1.3}
	\left\|
	\sum_{k=1}^{\infty} ((P^{k} F)(\theta', z ) - \nabla f(\theta' ) )
	-
	\sum_{k=1}^{\infty} ((P^{k} F)(\theta'', z ) - \nabla f(\theta'' ) )
	\right\|
	\leq 
	\tilde{C}_{2,Q} (1 + \|V\|^{2} ) 
	(\tilde{\varepsilon}_{Q}^{n} + n\|\theta' - \theta'' \| )
\end{align}
for all $\theta', \theta'' \in Q$, 
$z = (x,y,u,V) \in {\cal S}_{z}$, $n\geq 0$ 
($(P^{0} F)(\theta, z )$ is defined as $F(\theta, z )$). 

Let $\tilde{C}_{Q} = \max\{\tilde{C}_{1,Q}, \tilde{C}_{2,Q} \}$. 
Moreover, let  
$N_{Q,s}(t) = \lceil s \log t / \log \tilde{\varepsilon}_{Q} \rceil$
for $s,t \in (0,1)$
and 
$N_{Q,s}(t) = 0$
for $s\in (0,1)$, $t \in \{0\} \cup [1,\infty )$. 
Then, it can be concluded that there exists a real number 
$\tilde{K}_{Q,s} \in [1, \infty )$ such that 
\begin{align} \label{l2.1.5}
	N_{Q,s}(t) 
	+
	\tilde{\varepsilon}_{Q}^{N_{Q,s}(t) }
	\leq
	\tilde{K}_{Q,s} t^{s}
\end{align}
for all $t\in [0,\infty )$. 

For $\theta \in \Theta$, $z = (x,y,u,V) \in {\cal S}_{z}$, let 
\begin{align*}
	&
	\Phi(\theta, z )
	=
	\sum_{k=0}^{\infty } 
	((P^{k} F)(\theta, z ) - \nabla f(\theta ) ), 
	\\
	&
	\varphi_{Q,s}(z) 
	=
	\tilde{C}_{Q} \tilde{K}_{Q,s} \psi_{Q}(y) 
	(1 + \|V\|^{2} ). 
\end{align*}
Since 
\begin{align*}
	(P \varphi_{Q,s} )(\theta, z )
	=
	\tilde{C}_{Q} \tilde{K}_{Q,s} 
	(1 + \|H_{\theta}(u,V,y) \|^{2} ) 
	E(\psi_{Q}(Y_{2} )|X_{1} = x )
	< \infty
\end{align*}
for all $\theta \in \Theta$, 
$z = (x,y,u,V) \in {\cal S}_{z}$, 
we deduce from (\ref{l2.1.1}) that 
$\Phi(\cdot, \cdot )$ is well-defined, 
integrable and satisfies 
(\ref{l2.1.1*}), (\ref{l2.1.3*})
(notice that 
$(P\Phi)(\theta, z ) = \sum_{k=1}^{\infty } ((P^{k} F)(\theta, z ) - \nabla f(\theta ) )$). 
On the other hand, 
(\ref{l2.1.3}), (\ref{l2.1.5}) imply 
\begin{align*}
	\left\|
	(P\Phi )(\theta', z ) - (P\Phi )(\theta'', z )
	\right\|
	\leq 
	\tilde{C}_{Q} \tilde{K}_{Q,s} 
	(1 + \|V\|^{2} ) \|\theta' - \theta'' \|^{s}
\end{align*}
for any $\theta', \theta'' \in Q$, 
$z = (x,y,u,V) \in {\cal S}_{z}$
(set $n = N_{Q,s}(\|\theta' - \theta'' \| )$ in (\ref{l2.1.3})). 
Thus, (\ref{l2.1.5*}) is true for each 
$\theta', \theta'' \in Q$, 
$z = (x,y,u,V) \in {\cal S}_{z}$. 

Let $\boldsymbol \theta = \{\theta_{n} \}_{n\geq 0}$
and 
$\boldsymbol Y = \{Y_{n} \}_{n\geq 1}$. 
Due to Lemma \ref{lemma1.3}, we have 
\begin{align} \label{l2.1.7}
	\varphi_{Q,s}(Z_{n+1} )
	I_{ \{\tau_{Q} > n \} }
	= &
	\tilde{C}_{Q} \tilde{K}_{Q,s} \psi_{Q}(Y_{n+1} ) 
	(1 + \|H_{\boldsymbol \theta, \boldsymbol Y}^{0:n}(U_{0},V_{0} ) \|^{2} )
	I_{ \{\tau_{Q} > n \} }
	\nonumber\\
	\leq &
	4\tilde{C}_{Q} \tilde{K}_{Q,s} C_{4,Q }^{2} \psi_{Q}(Y_{n+1} ) 
	(1 + \|V_{0} \|^{2} )
\end{align}
for each $n\geq 0$
(notice that 
$H_{\boldsymbol \theta, \boldsymbol Y}^{0:n}(U_{0},V_{0} )$
depends only on the first $n$ elements of $\boldsymbol \theta$, 
and that
$\theta_{1},\dots,\theta_{n} \in Q$ is 
sufficient for (\ref{l2.1.7}) to hold). 
Consequently, 
\begin{align*}
	E\left(\left.
	\varphi_{Q,s}^{2}(Z_{n+1} ) 
	I_{ \{\tau_{Q} > n \} }
	\right|
	Z_{1}=z
	\right)
	\leq &
	16 \tilde{C}_{Q}^{2} \tilde{K}_{Q,s}^{2} C_{4,Q}^{4} 
	(1 + \|V\| )^{4} 
	E(\psi_{Q}^{2}(Y_{n+1} )|X_{1}=x ) 
	\\
	\leq &
	16 \tilde{C}_{Q}^{2} \tilde{K}_{Q,s}^{2} C_{4,Q}^{4} 
	(1 + \|V\| )^{4} 
	\max_{x'\in {\cal X} } 
	\int \psi_{Q}^{2}(y') Q(dy'|x')
	< \infty 
\end{align*}
for all $z = (x,y,u,V) \in {\cal S}_{z}$, $n\geq 0$. 
Hence, (\ref{l2.1.7*}) is true for all $z \in {\cal S}_{z}$. 
\end{IEEEproof}

\begin{lemma} \label{lemma2.2}
Suppose that Assumption \ref{a1} holds. 
Then, there exists a real number $s \in (0,1)$ such that 
\linebreak
$\sum_{n=0}^{\infty } \alpha_{n}^{1+s} \gamma_{n}^{r} < \infty$. 
\end{lemma} 

\begin{IEEEproof}
Let $p=(2+2r)/(2+r)$, 
$q=(2+2r)/r$, 
$s=(2+r)/(2+2r)$. 
Then, using the H\"{o}lder inequality, we get
\begin{align*}
	\sum_{n=0}^{\infty } 
	\alpha_{n}^{1+s} \gamma_{n}^{r} 
	=
	\sum_{n=1}^{\infty } 
	(\alpha_{n}^{2} \gamma_{n}^{2r} )^{1/p}
	\left(
	\frac{\alpha_{n} }{\gamma_{n}^{2} } 
	\right)^{1/q} 
	\leq
	\left(
	\sum_{n=1}^{\infty }
	\alpha_{n}^{2} \gamma_{n}^{2r} 
	\right)^{1/p}
	\left(
	\sum_{n=1}^{\infty } 
	\frac{\alpha_{n} }{\gamma_{n}^{2} } 
	\right)^{1/q}. 
\end{align*}
Since 
$\gamma_{n+1}/\gamma_{n} = 1 + \alpha_{n}/\gamma_{n} = O(1)$ for 
$n\rightarrow \infty$ and 
\begin{align*}
	\sum_{n=1}^{\infty } 
	\frac{\alpha_{n} }{\gamma_{n}^{2} }
	=
	\sum_{n=1}^{\infty } 
	\frac{\gamma_{n+1} - \gamma_{n} }{\gamma_{n}^{2} }
	\leq 
	\sum_{n=1}^{\infty } 
	\left(\frac{\gamma_{n+1} }{\gamma_{n} } \right)^{2} 
	\int_{\gamma_{n} }^{\gamma_{n+1} } 
	\frac{dt}{t^{2} }
	=
	\frac{1}{\gamma_{1} }
	\max_{n\geq 0} \left(\frac{\gamma_{n+1} }{\gamma_{n} } \right)^{2}, 
\end{align*}
it is obvious that 
$\sum_{n=0}^{\infty } \alpha_{n}^{1+s} \gamma_{n}^{r}$ converges. 
\end{IEEEproof} 

\begin{lemma} \label{lemma2.3}
Suppose that Assumptions \ref{a1} -- \ref{a4} hold. 
Then, there exists an event 
$N_{0}$ such that 
$P(N_{0} ) = 0$ and 
such that 
$\sum_{n=0}^{\infty } \alpha_{n} \gamma_{n}^{r} \xi_{n}$, 
$\sum_{n=0}^{\infty } \alpha_{n} \xi_{n}$
and 
$\sum_{n=0}^{\infty } \phi_{n}$
converge on 
$\Lambda \setminus N_{0}$. 
\end{lemma}

\begin{IEEEproof}
Let 
$Q\subset\Theta$ be an arbitrary compact set, while 
$t$ is an arbitrary number from $[0,r]$. 
Moreover, let  
$\Psi:\Theta \rightarrow \mathbb{R}^{d_{\theta} \times d_{\theta } }$ be an 
arbitrary locally Lipschitz continuous function. 
Obviously, 
in order to prove the lemma, it is sufficient to 
demonstrate that 
$\sum_{n=0}^{\infty } \alpha_{n} \gamma_{n}^{t} \Psi(\theta_{n} ) \xi_{n}$
and $\sum_{n=0}^{\infty } \phi''_{n}$ converge 
w.p.1 on 
$\bigcap_{n=0}^{\infty } \{\theta_{n} \in Q \}$
(to show the convergence of 
$\sum_{n=0}^{\infty } \alpha_{n} \gamma_{n}^{r} \xi_{n}$, 
set $t=r$ and $\Psi(\theta ) = I$ for all $\theta \in \Theta$, 
where $I$ stands for $d_{\theta } \times d_{\theta }$ unit matrix; 
to demonstrate the convergence of $\sum_{n=0}^{\infty } \phi'_{n}$, 
set $t=0$ and $\Psi(\theta ) = e (\nabla f(\theta ) )^{T}$ for each $\theta \in \Theta$, 
where
$e=[1\cdots 1]^{T} \in \mathbb{R}^{d_{\theta } }$).  

Let $s\in (0,1)$ be a real number such that 
$\sum_{n=0}^{\infty} \alpha_{n}^{1+s} \gamma_{n}^{r} < \infty$, 
while 
\begin{align*}
	\tilde{C}_{Q}
	=
	\max\left\{
	\|\nabla \Psi(\theta ) \|, 
	\frac{\|\Psi(\theta' ) - \Psi(\theta'' ) \|}{\|\theta' - \theta'' \|^{s} }, 
	\frac{\|\nabla f(\theta' ) - \nabla f(\theta'' ) \|}{\|\theta' - \theta'' \| }
	: \theta, \theta', \theta'' \in Q
	\right\}. 
\end{align*}
Moreover, for $n\geq 1$, let 
\begin{align*} 
	& 
	\psi_{1,n} 
	=
	\Psi(\theta_{n} )
	(\Phi(\theta_{n}, Z_{n+1} ) - (P\Phi )(\theta_{n}, Z_{n} ) ), 
	\\
	&
	\psi_{2,n}
	=
	\Psi(\theta_{n} )
	((P\Phi )(\theta_{n}, Z_{n} ) - (P\Phi )(\theta_{n-1}, Z_{n} ) )
	+
	(\Psi(\theta_{n} ) - \Psi(\theta_{n-1} ) ) (P\Phi )(\theta_{n-1},Z_{n} ), 
	\\
	&
	\psi_{3,n}
	=
	\Psi(\theta_{n} ) 
	(P\Phi )(\theta_{n}, Z_{n+1} ).  
\end{align*}
Then, it is straightforward to verify 
\begin{align} \label{l2.3.1}
	\sum_{i=1}^{n} \alpha_{i} \gamma_{i}^{t} \Psi(\theta_{i} ) \xi_{i} 
	=
	\sum_{i=1}^{n} \alpha_{i} \gamma_{i}^{t} \psi_{1,i} 
	+
	\sum_{i=1}^{n} \alpha_{i} \gamma_{i}^{t} \psi_{2,i} 
	+
	\sum_{i=0}^{n-1} 
	(\alpha_{i+1} \gamma_{i+1}^{t} - \alpha_{i} \gamma_{i}^{t} ) \psi_{3,i} 
	-
	\alpha_{n} \gamma_{n}^{t} \psi_{3,n} 
	+
	\alpha_{0} \gamma_{0}^{t} \psi_{3,0}
\end{align}
for $n\geq 1$. 

Owing to
Assumption \ref{a1}, we have 
\begin{align*}
	&
	\alpha_{n}
	=
	\alpha_{n+1}
	(1 + \alpha_{n} (\alpha_{n+1}^{-1} - \alpha_{n}^{-1} ) )
	=
	O(\alpha_{n+1} ), 
	\\
	&
	\alpha_{n} - \alpha_{n+1} 
	=
	\alpha_{n} \alpha_{n+1} (\alpha_{n+1}^{-1} - \alpha_{n}^{-1} ) 
	=
	O(\alpha_{n+1}^{2} ), 
	\\
	&
	\gamma_{n+1}^{t} - \gamma_{n}^{t} 
	=
	\gamma_{n}^{t} 
	\left((1 + \alpha_{n}/\gamma_{n} )^{t} - 1 \right)
	=
	o(\alpha_{n} \gamma_{n}^{t} )
\end{align*}
as $n\rightarrow \infty$. 
Consequently, 
\begin{align}
	& \label{l2.3.3}
	\sum_{n=0}^{\infty } 
	\alpha_{n}^{s} \alpha_{n+1} \gamma_{n+1}^{t} 	
	=
	\sum_{n=0}^{\infty } 
	(\alpha_{n}/\alpha_{n+1} )^{s} \alpha_{n+1}^{s} \gamma_{n+1}^{t}
	< 
	\infty, 
	\\
	& \label{l2.3.5} 
	\sum_{n=0}^{\infty } 
	|\alpha_{n} \gamma_{n}^{t} - \alpha_{n+1} \gamma_{n+1}^{t} |
	\leq 
	\sum_{n=0}^{\infty } 
	\alpha_{n} |\gamma_{n}^{t} - \gamma_{n+1}^{t} |
	+
	\sum_{n=0}^{\infty } 
	|\alpha_{n} - \alpha_{n+1} | \gamma_{n+1}^{t} 
	< 
	\infty. 
\end{align}
On the other side, as a consequence of Lemma \ref{lemma2.1}, we get 
\begin{align*}
	&
	E_{\theta,z}
	\left(
	|\psi_{1,n} |^{2}
	I_{ \{\tau_{Q} > n \} }
	\right)
	\leq 
	2 \tilde{C}_{Q}^{2} 
	E_{\theta,z}
	\left(
	\varphi_{Q,s }^{2}(Z_{n+1} ) 
	I_{ \{\tau_{Q} > n \} }
	\right)
	+
	2 \tilde{C}_{Q}^{2} 
	E_{\theta,z}
	\left(
	\varphi_{Q,s }^{2}(Z_{n} ) 
	I_{ \{\tau_{Q} > n-1 \} }
	\right), 
	\\
	&
	E_{\theta,z}
	\left(
	|\psi_{2,n} | 
	I_{ \{\tau_{Q} > n \} }
	\right)
	\leq 
	2 \tilde{C}_{Q} 
	E_{\theta,z}
	\left(
	\varphi_{Q,s }(Z_{n} )
	\|\theta_{n} - \theta_{n-1} \|^{s} 
	I_{ \{\tau_{Q} > n \} }
	\right)
	\leq  
	2 \tilde{C}_{Q}
	\alpha_{n-1}^{s} 
	E_{\theta,z}
	\left(
	\varphi_{Q,s }^{2}(Z_{n} ) 
	I_{ \{\tau_{Q} > n-1 \} }
	\right)
\end{align*}
for all $\theta \in \Theta$, 
$z\in {\cal S}_{z}$, $n\geq 1$. 
Due to the same lemma, we have 
\begin{align*}	 
	&
	E_{\theta,z}
	\left(
	|\psi_{3,n} |^{2}
	I_{ \{\tau_{Q} > n \} }
	\right)
	\leq 
	\tilde{C}_{Q}^{2} 
	E_{\theta,z}
	\left(
	\varphi_{Q,s }^{2}(Z_{n+1} )
	I_{ \{\tau_{Q} > n \} }
	\right),
	\\
	&
	E_{\theta,z}\left(
	|\phi''_{n} |
	I_{ \{\tau_{Q} > n \} } 
	\right)
	\leq 
	\tilde{C}_{Q}
	E_{\theta,z}\left(
	\|\theta_{n+1} - \theta_{n} \|^{2} 
	I_{ \{\tau_{Q} > n \} } 
	\right)
	\leq 
	\tilde{C}_{Q} \alpha_{n}^{2} 
	E_{\theta,z}
	\left(
	\varphi_{Q,s }^{2}(Z_{n+1} )
	I_{ \{\tau_{Q} > n \} }
	\right)
\end{align*}
for all $\theta \in \Theta$, 
$z\in {\cal S}_{z}$, $n\geq 1$. 
Then, Lemma \ref{lemma2.1} and (\ref{l2.3.3}) yield  
\begin{align*}
	&
	E_{\theta,z}
	\left(
	\sum_{n=1}^{\infty } 
	\alpha_{n}^{2} \gamma_{n}^{2t} 
	|\psi_{1,n} |^{2}
	I_{ \{\tau_{Q} > n \} }
	\right)
	\leq 
	4 \tilde{C}_{Q}^{2} 
	\left(
	\sum_{n=1}^{\infty } 
	\alpha_{n}^{2} \gamma_{n}^{2t} 
	\right)
	\sup_{n\geq 0}
	E_{\theta,z}
	\left(
	\varphi_{Q,s }^{2}(Z_{n+1} ) 
	I_{ \{\tau_{Q} > n \} }
	\right)
	< 
	\infty, 
	\\
	&
	E_{\theta,z}
	\left(
	\sum_{n=1}^{\infty } 
	\alpha_{n} \gamma_{n}^{t} 
	|\psi_{2,n} | 
	I_{ \{\tau_{Q} > n \} }
	\right)
	\leq 
	2 \tilde{C}_{Q}
	\left(
	\sum_{n=1}^{\infty } 
	\alpha_{n-1}^{s} \alpha_{n} \gamma_{n}^{t} 
	\right)
	\sup_{n\geq 0}
	E_{\theta,z}
	\left(
	\varphi_{Q,s }^{2}(Z_{n+1} ) 
	I_{ \{\tau_{Q} > n \} }
	\right)
	<
	\infty 
\end{align*}
for any $\theta \in \Theta$, 
$z \in {\cal S}_{z}$. 
On the other side, Lemma \ref{lemma2.1} and (\ref{l2.3.5}) imply 
\begin{align*}
	&
	E_{\theta,z}
	\left(
	\sum_{n=1}^{\infty } 
	|\alpha_{n} \gamma_{n}^{t} - \alpha_{n+1} \gamma_{n+1}^{t} | \,
	|\psi_{3,n} | 
	I_{ \{\tau_{Q} > n \} }
	\right)
	\\
	& \;\;\; 
	\leq 
	\tilde{C}_{Q} 
	\left(
	\sum_{n=1}^{\infty } 
	|\alpha_{n} \gamma_{n}^{t} - \alpha_{n+1} \gamma_{n+1}^{t} |
	\right)
	\sup_{n\geq 0}
	\left(
	E_{\theta,z}
	\left(
	\varphi_{Q,s }^{2}(Z_{n+1} ) 
	I_{ \{\tau_{Q} > n \} }
	\right)
	\right)^{1/2} 
	<
	\infty, 
	\\
	&
	E_{\theta,z}
	\left(
	\sum_{n=1}^{\infty } 
	\alpha_{n+1}^{2} \gamma_{n+1}^{2t} 
	|\psi_{3,n} |^{2} 
	I_{ \{\tau_{Q} > n \} }
	\right)
%	\\
%	& \;\;\; 
	\leq 
	\tilde{C}_{Q}^{2} 
	\left(
	\sum_{n=1}^{\infty } 
	\alpha_{n+1}^{2} \gamma_{n+1}^{2t} 
	\right)
	\sup_{n\geq 0}
	E_{\theta,z}
	\left(
	\varphi_{Q,s }^{2}(Z_{n+1} ) 
	I_{ \{\tau_{Q} > n \} }
	\right)
	<
	\infty, 
	\\
	&
	E_{\theta,z}\left(
	\sum_{n=0}^{\infty } 
	|\phi''_{n} | 
	I_{ \{\tau_{Q} > n \} } 
	\right)  
	\leq 
	\tilde{C}_{Q}
	\left(
	\sum_{n=0}^{\infty } \alpha_{n}^{2} 
	\right)
	\sup_{n\geq 0} 
	E_{\theta,z}\left(
	\varphi_{Q,s}^{2}(Z_{n+1} ) 
	I_{ \{\tau_{Q} > n \} } 
	\right)
	< \infty
\end{align*}
for each $\theta \in \Theta$, 
$z \in {\cal S}_{z}$.  
Since 
\begin{align*}
	E_{\theta,z}
	\left(
	\psi_{1,n} 
	I_{ \{\tau_{Q} > n \} }
	|
	{\cal F}_{n} 
	\right)
	=
	\Psi(\theta_{n} ) 
	\left(
	E_{\theta,z}
	\left(
	\Phi(\theta_{n}, Z_{n+1} )
	|
	{\cal F}_{n} 
	\right)
	-
	(P\Phi )(\theta_{n},Z_{n} )
	\right) 
	I_{ \{\tau_{Q} > n \} }
	=
	0
\end{align*}
w.p.1 for every $\theta\in \Theta$, 
$z \in {\cal S}_{z}$, $n\geq 1$,  
it is clear that 
series 
\begin{align*}
	\sum_{n=1}^{\infty } \alpha_{n} \gamma_{n}^{t} \psi_{1,n}, 
	\;\;\; 
	\sum_{n=1}^{\infty } \alpha_{n} \gamma_{n}^{t} \psi_{2,n}, 
	\;\;\; 
	\sum_{n=1}^{\infty } (\alpha_{n} \gamma_{n}^{t} - \alpha_{n+1} \gamma_{n+1}^{t} ) \psi_{3,n}, 
	\;\;\; 
	\sum_{n=1}^{\infty } \phi''_{n}
\end{align*}
converge w.p.1 
on $\bigcap_{n=0}^{\infty } \{\theta_{n} \in Q \}$, 
and that 
$
	\lim_{n\rightarrow \infty } \alpha_{n} \gamma_{n}^{t} \psi_{3,n} = 0
$	
w.p.1 on the same event. 
Owing to this and (\ref{l2.3.1}), 
we have that 
$\sum_{n=0}^{\infty } \alpha_{n} \gamma_{n}^{t} \Psi(\theta_{n} ) \xi_{n}$
is convergent w.p.1 on $\bigcap_{n=0}^{\infty } \{\theta_{n} \in Q \}$. 
\end{IEEEproof}

\begin{lemma} \label{lemma2.4}
Suppose that Assumption \ref{a1} -- \ref{a4} hold. 
Then, on $\Lambda\setminus N_{0}$, 
$\lim_{n\rightarrow \infty } \nabla f(\theta_{n} ) = 0$
and 
$\lim_{n\rightarrow \infty } f(\theta_{n} )$
exists. 
\end{lemma}

\begin{IEEEproof}
Let $Q\subset\Theta$ be an arbitrary compact set, 
while $\omega$ is an arbitrary sample from 
$\bigcap_{n=0}^{\infty } \{\theta_{n} \in Q \} \setminus N_{0}$
(notice that all formulas which appear in the proof correspond to 
this $\omega$). 
Obviously, in order to prove the lemma, it is sufficient to 
show that $\lim_{n\rightarrow \infty } f(\theta_{n} )$ exists
and that 
$\lim_{n\rightarrow \infty } \nabla f(\theta_{n} ) = 0$. 

Since $\sum_{n=0}^{\infty } \phi_{n}$ converges and 
\begin{align*}
	\sum_{i=0}^{n-1} \alpha_{i} \|\nabla f(\theta_{i} ) \|^{2}
	=
	f(\theta_{n} ) - f(\theta_{0} ) 
	-
	\sum_{i=0}^{n-1} \phi_{i}
\end{align*}
for $n\geq 0$,  
we conclude 
$\sum_{n=0}^{\infty } \alpha_{n} \|\nabla f(\theta_{n} ) \|^{2} < \infty$
(also notice that $f(\cdot )$ is bounded on $Q$).  
As 
\begin{align*}
	f(\theta_{n} )
	=
	f(\theta_{0} )
	+
	\sum_{i=0}^{n-1} \alpha_{i} \|\nabla f(\theta_{i} ) \|^{2} 
	+
	\sum_{i=0}^{n-1} \phi_{i}
\end{align*}
for $n\geq 0$, it is clear that 
$\lim_{n\rightarrow \infty } f(\theta_{n} )$ exists. 

Let $\tilde{C}_{Q}$ be a Lipschitz constant of $\nabla f(\cdot )$
on $Q$ and an upper bound of 
$\|\nabla f(\cdot ) \|$ on the same set. 
Now, we prove 
$\lim_{n\rightarrow \infty } \nabla f(\theta_{n} ) = 0$. 
Suppose the opposite. 
Then, there exist $\varepsilon \in (0, \infty )$
and sequences 
$\{m_{k} \}_{k\geq 0}$, $\{n_{k} \}_{k\geq 0}$
(all depending on $\omega$) 
such that 
$m_{k} < n_{k} < m_{k+1}$, 
$\|\nabla f(\theta_{m_{k} } ) \| \leq \varepsilon$, 
$\|\nabla f(\theta_{n_{k} } ) \| \geq 2\varepsilon$
for $k\geq 0$, 
and such that 
$\|\nabla f(\theta_{n} ) \| \geq \varepsilon$
for $m_{k} < n \leq n_{k}$, $k\geq 0$. 
Therefore, 
\begin{align} \label{l2.4.1}
	\varepsilon
	\leq 
	\|\nabla f(\theta_{n_{k} } ) - \nabla f(\theta_{m_{k} } ) \|
	\leq 
	\tilde{C}_{Q} 
	\|\theta_{n_{k} } - \theta_{m_{k} } \|
	\leq 
	\tilde{C}_{Q}^{2} 
	\sum_{i=m_{k} }^{n_{k}-1} \alpha_{i} 
	+
	\tilde{C}_{Q} 
	\left\|
	\sum_{i=m_{k} }^{n_{k} - 1 } \alpha_{i} \xi_{i} 
	\right\|
\end{align}
for $k\geq 0$. 
We also have
\begin{align*}
	\varepsilon^{2} 
	\sum_{i=m_{k} + 1 }^{n_{k} } \alpha_{i} 
	\leq 
	\sum_{i=m_{k} + 1 }^{\infty } \alpha_{i} \|\nabla f(\theta_{i} ) \|^{2}
\end{align*}
for $k\geq 0$. 
Consequently, 
$\lim_{k\rightarrow \infty } 
\sum_{i=m_{k} }^{n_{k} - 1 } \alpha_{i} = 0$. 
However, this is not possible, since the limit process
$k\rightarrow \infty$ applied to (\ref{l2.4.1}) 
would imply 
\begin{align*}
	\varepsilon 
	\leq 
	\lim_{k\rightarrow \infty } 
	\|\nabla f(\theta_{n_{k} } ) - \nabla f(\theta_{m_{k} } ) \|
	=
	0. 
\end{align*}
Hence, $\lim_{n\rightarrow \infty } \nabla f(\theta_{n} ) = 0$. 
\end{IEEEproof}

\subsection{Convergence and Convergence Rate}\label{subsection1.4*} 

In this subsection, 
using the results of Subsections \ref{subsection1.2*}, \ref{subsection1.3*} (Corollary \ref{corollary1}, Lemmas \ref{lemma2.2}, \ref{lemma2.3}), 
the convergence and convergence rate of recursion 
(\ref{1.1}) -- (\ref{1.5}) are analyzed  
and Theorems \ref{theorem2} and \ref{theorem3} are proved. 

Throughout the subsection, we use the following notation. 
For $t\in (0,\infty )$, $n\geq 0$, let 
\begin{align*}
	a(n,t)
	=
	\max\{k\geq n: \gamma_{k} - \gamma_{n} \leq t\}. 
\end{align*}
For $0\leq n \leq k$, let 
\begin{align*}
	&
	\zeta_{n} 
	=
	\sup_{k\geq n} 
	\left\|
	\sum_{i=n}^{k} \alpha_{i} \xi_{i} 
	\right\|, 
	\\
	&
	\varepsilon'_{n,k} 
	=
	\sum_{i=n}^{k-1} \alpha_{i} \xi_{i}, 
	\\
	&
	\varepsilon''_{n,k} 
	=
	\sum_{i=n}^{k-1} \alpha_{i} (\nabla f(\theta_{i} ) - \nabla f(\theta_{n} ) ), 
	\\
	&
	\phi'_{n,k}
	=
	(\nabla f(\theta_{n} ) )^{T} (\varepsilon'_{n,k} + \varepsilon''_{n,k} ), 
	\\
	&
	\phi''_{n,k} 
	=
	\int_{0}^{1} 
	(\nabla f(\theta_{n} + t (\theta_{k} - \theta_{n} ) ) 
	- 
	\nabla f(\theta _{n} ) )^{T}
	(\theta_{k} - \theta_{n} ) dt, 
\end{align*}
while 
$\varepsilon_{n,k} = \varepsilon'_{n,k} + \varepsilon''_{n,k}$
and 
$\phi_{n,k} = \phi'_{n,k} + \phi''_{n,k}$. 
Then, it is straightforward to verify
\begin{align}
	\label{1*.3.1}
	\theta_{k}
	= &
	\theta_{n} 
	+
	\sum_{i=n}^{k-1} \alpha_{i} \nabla f(\theta_{i} ) 
	+
	\varepsilon'_{n,k}
	\nonumber\\
	=&
	\theta_{n}
	+
	(\gamma_{k} - \gamma_{n} ) \nabla f(\theta_{n} ) 
	+
	\varepsilon_{n,k}, 
	\\
	\label{1*.3.3}
	f(\theta_{k} )
	=&
	f(\theta_{n} ) 
	+
	(\gamma_{k} - \gamma_{n} ) \|\nabla f(\theta_{n} ) \|^{2} 
	+
	\phi_{n,k}
\end{align}
for $0\leq n \leq k$. 

Besides the notation introduced in the previous paragraph, 
we also rely on the following notation in this subsection. 
For a compact set $Q\subset \Theta$, 
$C_{Q} \in [1,\infty )$
denotes an upper bound of 
$\|\nabla f(\cdot ) \|$ on $Q$
and a Lipschitz constant of $\nabla f(\cdot )$ 
on the same set. 
$\hat{A}$ is the set of the accumulation points of 
$\{\theta_{n} \}_{n\geq 0}$, 
while 
\begin{align*}
	\hat{f}
	=
	\liminf_{n\rightarrow \infty } f(\theta_{n} ). 
\end{align*}
$\hat{\rho}$ and 
$\hat{B}$, $\hat{Q}$
are a random quantity 
and
random sets (respectively) defined by 
\begin{align*}
	\hat{\rho}
	=
	d(\hat{A}, \partial\Theta )/2, 
	\;\;\;\;\; 
	\hat{B}
	=
	\bigcup_{\theta \in \hat{A} } 
	\left\{
	\theta'\in\mathbb{R}^{d_{\theta } }: 
	\|\theta' - \theta \| \leq \min\{\delta_{\theta }, \hat{\rho } \} 
	\right\}, 
	\;\;\;\;\; 
	\hat{Q} = \text{cl}(\hat{B} )
\end{align*}
on $\Lambda$,
and by
\begin{align*}
	\hat{\rho}
	=
	0, 
	\;\;\;\;\; 
	\hat{B} = \hat{A}, 
	\;\;\;\;\; 
	\hat{Q} = \hat{A}
\end{align*}
otherwise. 
Overriding the definition of $\hat{\mu}$ in Theorem \ref{theorem3}, 
we specify 
random quantities 
$\hat{\delta}$, $\hat{\mu}$, $\hat{C}$, $\hat{C}$
as 
\begin{align}\label{1*.3.5}
	\hat{\delta}
	=
	\delta_{\hat{Q}, \hat{f} }, 
	\;\;\;\;\; 
	\hat{\mu}
	=
	\mu_{\hat{Q}, \hat{f} }, 
	\;\;\;\;\;
	\hat{C}
	=
	C_{\hat{Q}, \hat{f} }, 
	\;\;\;\;\; 
	\hat{M}
	=
	M_{\hat{Q}, \hat{f} }
\end{align}
on $\Lambda$, 
and as 
\begin{align*}
	\hat{\delta}
	=
	1, 
	\;\;\;\;\; 
	\hat{\mu}
	=
	2, 
	\;\;\;\;\;
	\hat{C}
	=
	1, 
	\;\;\;\;\; 
	\hat{M}
	=
	1
\end{align*}
otherwise
($\delta_{Q,a}$, $\mu_{Q,a}$, $M_{Q,a}$ are introduced in 
the statement of Corollary \ref{corollary1};
later, once Theorem \ref{theorem2} is proved, it will be clear that the 
definitions of $\hat{\mu}$ provided in Theorem \ref{theorem3} and in (\ref{1*.3.5}) 
are equivalent). 
Random quantities 
$\hat{p}$, $\hat{q}$, $\hat{r}$
are defined in the same way as in 
(\ref{t3.3*}). 
%\begin{align*}
%	\hat{r}
%	=
%	\begin{cases}
%	1/(2 - \hat{\mu} ), 
%	&\text{ if } \hat{\mu} < 2
%	\\
%	\infty, 
%	&\text{ if } \hat{\mu} = 2
%	\end{cases},
%	\;\;\;\;\; 
%	\hat{p}
%	=
%	\min\{\hat{r}, r \}, 
%	\;\;\;\;\; 
%	\hat{q}
%	=
%	\min\{(\hat{p} - 1 )/2, r-1 \}. 
%\end{align*}
Functions $u(\cdot )$ and $v(\cdot )$ are defined by 
\begin{align*}
	u(\theta )
	=
	\hat{f} - f(\theta ), 
	\;\;\;\;\;  
	v(\theta )
	=
	\begin{cases}
	(1/u(\theta ) )^{1/\hat{p} }, 
	&\text{ if } u(\theta ) > 0
	\\
	0, 
	&\text{ otherwise }
	\end{cases}
\end{align*}
for $\theta \in \Theta$. 

Obviously, on event $\Lambda$, 
$\hat{Q}$ is compact and 
satisfies 
$\hat{A} \subset \text{int}\hat{Q}$, $\hat{Q} \subset \Theta$. 
Thus, 
$\hat{\mu}$, $\hat{M}$, $\hat{p}$, $\hat{q}$, $\hat{r}$, 
$v(\cdot )$
are are well-defined on the same event
(what happens with these quantities outside 
$\Lambda$ does not affect the results provided in this subsection). 
On the other side, Corollary \ref{corollary1} implies 
\begin{align}\label{1*.3.7}
	|f(\theta ) - \hat{f} |
	\leq 
	\hat{M} \|\nabla f(\theta ) \|^{\hat{\mu} }
\end{align}
on $\Lambda$
for all $\theta \in \hat{Q}$
satisfying 
$|f(\theta ) - \hat{f} | \leq \hat{\delta}$. 

\begin{lemma} \label{lemma3.1}
Suppose that Assumptions \ref{a1} -- \ref{a4} hold. 
Then, 
$\lim_{n\rightarrow \infty } \gamma_{n}^{r} \zeta_{n} = 0$
on $\Lambda \setminus N_{0}$
($N_{0}$ is specified in the statement of Lemma \ref{lemma2.3}). 
\end{lemma}

\begin{IEEEproof}
It is straightforward to verify 
\begin{align*}
	\sum_{i=n}^{k} \gamma_{i} \xi_{i} 
	=
	\gamma_{k+1}^{-r} 
	\sum_{j=n}^{k} \alpha_{j} \gamma_{j}^{r} \xi_{j} 
	+
	\sum_{i=n}^{k} 
	(\gamma_{i}^{-r} - \gamma_{i+1}^{-r} )
	\sum_{j=n}^{i} \alpha_{j} \gamma_{j}^{r} \xi_{j}
\end{align*}
for $0\leq n \leq k$. 
Therefore, 
\begin{align*}
	\left\|\sum_{i=n}^{k} \gamma_{i} \xi_{i}  \right\|
	\leq 
	\left(
	\gamma_{k+1}^{-r} 
	+
	\sum_{i=n}^{k} (\gamma_{i}^{-r} - \gamma_{i+1}^{-r} ) 
	\right)
	\sup_{i\geq n} 
	\left\|\sum_{j=n}^{i} \alpha_{j} \gamma_{j}^{r} \xi_{j} \right\|
	=
	\gamma_{n}^{-r} 
	\sup_{i\geq n} 
	\left\|\sum_{j=n}^{i} \alpha_{j} \gamma_{j}^{r} \xi_{j} \right\|
\end{align*}
for $0\leq n \leq k$. 
Consequently, Lemma \ref{lemma2.3} implies 
\begin{align*}
	\limsup_{n\rightarrow \infty }
	\gamma_{n}^{r} \zeta_{n}
	=
	\limsup_{n\rightarrow \infty} 
	\sup_{k\geq n} 
	\left\|\sum_{i=n}^{k} \alpha_{i} \gamma_{i}^{r} \xi_{i} \right\|
	=
	0
\end{align*}
on $\Lambda \setminus N_{0}$. 
\end{IEEEproof}

\begin{lemma}\label{lemma3.2}
Suppose that Assumptions \ref{a1} -- \ref{a4} hold. 
Let $\hat{C}_{1} = (16 \hat{p} \hat{M} )^{2\hat{p} }$
(notice that $1\leq \hat{C}_{1} < \infty$ everywhere). 
Then, there exist a random quantity 
$\hat{t}$
and an integer-valued random variable 
$\sigma$
such that 
$0 < \hat{t} < 1$, 
$0\leq \sigma < \infty$
everywhere and such that 
\begin{align}
	&\label{l3.2.1*}
	\max_{n\leq k \leq a(n,\hat{t} ) }
	\|\varepsilon_{n,k} \|
	\leq 
	(\hat{t}/\hat{C}_{1} )
	(\gamma_{n}^{-r} + \|\nabla f(\theta_{n} ) \| ), 
	\\
	&\label{l3.2.3*}
	\max_{n\leq k \leq a(n,\hat{t} ) }
	|\phi_{n,k} |
	\leq 
	(\hat{t}/\hat{C}_{1} )
	(\gamma_{n}^{-2r} + \|\nabla f(\theta_{n} ) \|^{2} ), 
	\\
	&\label{l3.2.5*}
	f(\theta_{n} )
	-
	f(\theta_{a(n, \hat{t} ) } )
	+
	2^{-1} \hat{t} \|\nabla f(\theta_{n} ) \|^{2} 
	\leq 
	(\hat{t}/\hat{C}_{1} ) 
	\gamma_{n}^{-2r}, 
	\\
	&\label{l3.2.7*}
	f(\theta_{n} )
	-
	f(\theta_{a(n, \hat{t} ) } )
	+
	2^{-1} \|\nabla f(\theta_{n} ) \| \|\theta_{a(n, \hat{t} ) } - \theta_{n} \|
	\leq 
	(\hat{t}/\hat{C}_{1} ) 
	\gamma_{n}^{-2r}
\end{align}
on $\Lambda \setminus N_{0}$
for $n>\sigma$. 
\end{lemma}

\begin{IEEEproof}
Let 
$\tilde{C}_{1} = 2 \hat{C} \exp(\hat{C} )$, 
$\tilde{C}_{2} = 2 \hat{C} \tilde{C}_{1}$,
$\tilde{C}_{3} = 2 \hat{C} \tilde{C}_{2}^{2} + \tilde{C}_{2}$
and $\tilde{C}_{4} = \tilde{C}_{2} + \tilde{C}_{3}$, 
while $\hat{t} = 1/(2 \hat{C}_{1} \tilde{C}_{4} )$. 
Moreover, let 
\begin{align*}
	\tilde{\sigma}_{1}
	=& 
	\max\left(
	\left\{
	n\geq 0: \theta_{n} \not\in \hat{Q} 
	\right\} 
	\cup 
	\{0 \}
	\right), 
	\\
	\tilde{\sigma}_{2}
	=& 
	\max\left(
	\left\{
	n\geq 0: \alpha_{n} > \hat{t}/4 
	\right\} 
	\cup 
	\{0 \}
	\right), 
	\\
	\tilde{\sigma}_{3}
	=& 
	\max\left(
	\left\{
	n\geq 0: \gamma_{n}^{r} \zeta_{n} > \hat{t}/(2\hat{C}_{1}\tilde{C}_{4} ) 
	\right\} 
	\cup 
	\{0 \}
	\right)
\end{align*}
while 
$\sigma = 
\max\{\tilde{\sigma}_{1}, \tilde{\sigma}_{2}, \tilde{\sigma}_{3} \}
I_{\Lambda \setminus N_{0} }$. 
Then, it is obvious that 
$\sigma$ is well-defined, 
while Lemma \ref{lemma3.1} implies 
$0\leq \sigma < \infty$ everywhere. 
We also have 
\begin{align}
	&\label{l3.2.1}
	\max\{
	\tilde{C}_{2} \gamma_{n}^{r} \zeta_{n}, 
	\tilde{C}_{3} \gamma_{n}^{r} \zeta_{n}, 
	\tilde{C}_{3} \gamma_{n}^{2r} \zeta_{n}^{2}, 
	\tilde{C}_{4} \gamma_{n}^{r} \zeta_{n}, 
	\tilde{C}_{4} \gamma_{n}^{2r} \zeta_{n}^{2}
	\}
	\leq 
	2^{-1} \hat{C}_{1}^{-1} \hat{t}, 
	\\
	&\label{l3.2.3}
	\max\{
	\tilde{C}_{2} \hat{t}^{2}, 
	\tilde{C}_{3} \hat{t}^{2}, 
	\tilde{C}_{4} \hat{t}^{2}
	\}
	\leq 
	2^{-1} \hat{C}_{1}^{-1} \hat{t}, 
	\\
	&\label{l3.2.5}
	\hat{t} 
	\geq
	\gamma_{a(n,\hat{t} ) } - \gamma_{n}
	=
	\gamma_{a(n,\hat{t} ) + 1 } 
	-
	\gamma_{n}
	-
	\alpha_{a(n,\hat{t} ) }
	\geq
	3\hat{t}/4
\end{align}
on $\Lambda \setminus N_{0}$ for $n > \sigma$. 

Let $\omega$ be an arbitrary sample from 
$\Lambda$
(notice that all formulas which follow in the proof correspond to 
this $\omega$). 
Since $\theta_{n} \in \hat{Q}$ for $n>\sigma$,
we have
\begin{align*}
	\|\nabla f(\theta_{k} ) \|
	\leq &
	\|\nabla f(\theta_{n} ) \| 
	+ 
	\|\nabla f(\theta_{k} ) - \nabla f(\theta_{n} ) \|
	\\
	\leq &
	\|\nabla f(\theta_{n} ) \| 
	+ 
	\hat{C}
	\|\theta_{k} - \theta_{n} \|
	\\
	\leq &
	\|\nabla f(\theta_{n} ) \| 
	+ 
	\hat{C} 
	\sum_{i=n}^{k-1} 
	\alpha_{i} 
	\|\nabla f(\theta_{i} ) \|
	+
	\hat{C} 
	\|\varepsilon'_{n,k} \|
	\\
	\leq &
	\hat{C}(\zeta_{n} + 	\|\nabla f(\theta_{n} ) \| )
	+ 
	\hat{C} 
	\sum_{i=n}^{k-1} 
	\alpha_{i} 
	\|\nabla f(\theta_{i} ) \|
\end{align*}
for $\sigma < n \leq k$. 
Then, 
Bellman-Gronwall inequality yields 
\begin{align*}
	\|\nabla f(\theta_{k} ) \|
	\leq 
	\hat{C} (\zeta_{n} + \|\nabla f(\theta_{n} ) \| ) 
	\exp(\hat{C} (\gamma_{k} - \gamma_{n} ) )
	\leq 
	\hat{C} \exp(\hat{C} ) 
	(\zeta_{n} + \|\nabla f(\theta_{n} ) \| ) 
\end{align*}
for $\sigma < n \leq k \leq a(n,1)$. 
Consequently, 
\begin{align*}
	\|\theta_{k} - \theta_{n} \|
	\leq &
	\sum_{i=n}^{k-1} 
	\alpha_{i} \|\nabla f(\theta_{i} ) \|
	+
	\|\varepsilon'_{n,k} \|
	\\
	\leq &
	\zeta_{n} 
	+
	\hat{C} \exp(\hat{C} ) 
	(\zeta_{n} + \|\nabla f(\theta_{n} ) \| )
	(\gamma_{k} - \gamma_{n} ) 
	\\
	\leq &
	\tilde{C}_{1}
	(\zeta_{n} + (\gamma_{k} - \gamma_{n} ) \|\nabla f(\theta_{n} ) \| )
\end{align*}
for $\sigma < n \leq k \leq a(n,1)$. 
Therefore, 
\begin{align}\label{l3.2.7}
	\|\varepsilon_{n,k} \|
	\leq &
	\|\varepsilon'_{n,k} \| 
	+
	\hat{C}
	\sum_{i=n}^{k-1} \alpha_{i} \|\theta_{i} - \theta_{n} \| 
	\nonumber\\
	\leq &
	\zeta_{n} 
	+ 
	\hat{C} \tilde{C}_{1} 
	((\gamma_{k} - \gamma_{n} ) \zeta_{n} 	
	+ 
	(\gamma_{k} - \gamma_{n} )^{2} \|\nabla f(\theta_{n} ) \| )
	\nonumber\\
	\leq &
	\tilde{C}_{2} 
	(\zeta_{n} + (\gamma_{k} - \gamma_{n} )^{2} \|\nabla f(\theta_{n} ) \| )
\end{align}
for $\sigma < n \leq k \leq a(n,1)$
(notice that $\gamma_{k} - \gamma_{n} \leq 1$
for $n\leq k \leq a(n,1)$). 
Thus, 
\begin{align}\label{l3.2.9}
	\|\phi_{n,k} \|
	\leq &
	\|\nabla f(\theta_{n} ) \| \|\varepsilon_{n,k} \| 
	+
	\hat{C} \|\theta_{k} - \theta_{n} \|^{2} 
	\nonumber\\
	\leq &
	\tilde{C}_{2} 
	(\zeta_{n} \|\nabla f(\theta_{n} ) \| 
	+
	(\gamma_{k} - \gamma_{n} )^{2} \|\nabla f(\theta_{n} ) \|^{2} )
	+
	2 \hat{C} \tilde{C}_{1}^{2}  
	(\zeta_{n}^{2} 
	+
	(\gamma_{k} - \gamma_{n} )^{2} \|\nabla f(\theta_{n} ) \|^{2} )
	\nonumber\\
	\leq &
	\tilde{C}_{3} 
	(\zeta_{n}^{2}
	+
	\zeta_{n} \|\nabla f(\theta_{n} ) \| 
	+
	(\gamma_{k} - \gamma_{n} )^{2} \|\nabla f(\theta_{n} ) \|^{2} )
\end{align}
for $\sigma < n \leq k \leq a(n,1)$. 
On the other side, combining (\ref{1*.3.1}), (\ref{1*.3.3}), 
we get
\begin{align*}
	f(\theta_{k} )
	-
	f(\theta_{n} )
	=&
	\|\nabla f(\theta_{n} ) \| 
	\|(\gamma_{k} - \gamma_{n} ) \nabla f(\theta_{n} ) \|
	+
	\phi_{n,k}
	\\
	=&
	\|\nabla f(\theta_{n} ) \| 
	\|\theta_{k} - \theta_{n} + \varepsilon_{n,k} \|
	+
	\phi_{n,k}
	\\
	\geq&
	\|\nabla f(\theta_{n} ) \| 
	(\|\theta_{k} - \theta_{n} \| - \|\varepsilon_{n,k} \| )
	-
	|\phi_{n,k} |
\end{align*}
for $0\leq n \leq k$. 
Then, (\ref{l3.2.7}), (\ref{l3.2.9}) yield
\begin{align}\label{l3.2.21}
	f(\theta_{n} ) 
	-
	f(\theta_{k} )
	+
	\|\nabla f(\theta_{n} ) \| \|\theta_{k} - \theta_{n} \|
	\leq &
	\|\nabla f(\theta_{n} ) \| \|\varepsilon_{n,k} \| 
	+
	|\phi_{n,k} |
	\nonumber\\
	\leq &
	\tilde{C}_{3} \zeta_{n}^{2}
	+
	(\tilde{C}_{2} + \tilde{C}_{3} )
	(\zeta_{n} \|\nabla f(\theta_{n} ) \| 
	+
	(\gamma_{k} - \gamma_{n} )^{2} \|\nabla f(\theta_{n} ) \|^{2} )
	\nonumber\\
	\leq &
	\tilde{C}_{4} (\zeta_{n}^{2}
	+
	\zeta_{n} \|\nabla f(\theta_{n} ) \| 
	+
	(\gamma_{k} - \gamma_{n} )^{2} \|\nabla f(\theta_{n} ) \|^{2} )
\end{align}
for $\sigma < n \leq k \leq a(n,1)$. 

Owing to (\ref{l3.2.1}), (\ref{l3.2.3}), (\ref{l3.2.7}), (\ref{l3.2.9}), we have 
\begin{align}\label{l3.2.23}
	\|\varepsilon_{n,k} \|
	\leq &
	\tilde{C}_{2} \zeta_{n} 
	+
	\tilde{C}_{2} \hat{t}^{2} \|\nabla f(\theta_{n} ) \| 
	\nonumber\\
	\leq &
	\hat{C}_{1}^{-1} \hat{t} 
	(\gamma_{n}^{-r} + \|\nabla f(\theta_{n} ) \| ), 
	\\
	\label{l3.2.25}
	|\phi_{n,k} |
	\leq &
	\tilde{C}_{3} \zeta_{n}^{2} 
	+
	\tilde{C}_{3} \zeta_{n} \|\nabla f(\theta_{n} ) \| 
	+
	\tilde{C}_{3} \hat{t}^{2} \|\nabla f(\theta_{n} ) \|^{2} 
	\nonumber\\
	\leq &
	2^{-1} \hat{C}_{1}^{-1} \hat{t}
	(\gamma_{n}^{-2r} 
	+
	\gamma_{n}^{-r} \|\nabla f(\theta_{n} ) \| 
	+
	\|\nabla f(\theta_{n} ) \|^{2} )
	\nonumber\\
	\leq &
	\hat{C}_{1}^{-1} \hat{t} 
	(\gamma_{n}^{-2r} + \|\nabla f(\theta_{n} ) \|^{2} )
\end{align}
for $\sigma < n \leq k \leq a(n,\hat{t} )$
(notice that $\gamma_{k} -\gamma_{n} \leq \hat{t}$ for $n\leq k \leq a(n,\hat{t} )$). 
Due to (\ref{1*.3.3}), (\ref{l3.2.5}), (\ref{l3.2.25}), we have also 
\begin{align}\label{l3.2.27}
	f(\theta_{n} )
	-
	f(\theta_{a(n,\hat{t} ) } )
	\leq &
	-
	(\gamma_{a(n,\hat{t} ) } - \gamma_{n} ) \|\nabla f(\theta_{n} ) \|^{2}
	+
	|\phi_{n,a(n,\hat{t} ) } |
	\nonumber \\
	\leq &
	-
	(3\hat{t}/4) \|\nabla f(\theta_{n} ) \|^{2} 
	+
	\hat{C}_{1}^{-1} \hat{t} 
	(\gamma_{n}^{-2r} + \|\nabla f(\theta_{n} ) \|^{2} )
	\nonumber \\
	= &
	-
	(3/4 - \hat{C}_{1}^{-1} ) \hat{t} \|\nabla f(\theta_{n} ) \|^{2} 
	+
	\hat{C}_{1}^{-1} \hat{t} \gamma_{n}^{-2r}
	\nonumber \\
	\leq &
	-
	2^{-1} \hat{t} \|\nabla f(\theta_{n} ) \|^{2} 
	+
	\hat{C}_{1}^{-1} \hat{t} \gamma_{n}^{-2r}
\end{align}
for $n > \sigma$
(notice that $\hat{C}_{1} \geq 4$). 
Consequently, 
\begin{align}\label{l3.2.27'}
	\hat{C}_{1}^{-1} \hat{t} 
	\|\nabla f(\theta_{n} ) \|^{2} 
	\leq 
	2^{-1} \hat{t} \|\nabla f(\theta_{n} ) \|^{2} 
	\leq 
	\hat{C}_{1}^{-1} \hat{t} \gamma_{n}^{-2r} 
	+
	(f(\theta_{a(n,\hat{t} ) } ) - f(\theta_{n} ) )
\end{align}
for $n>\sigma$. 
On the other side, (\ref{l3.2.1}) -- (\ref{l3.2.5}), (\ref{l3.2.21}), (\ref{l3.2.27'}) 
imply 
\begin{align*}
	f(\theta_{n} )
	-
	f(\theta_{a(n,\hat{t} ) } )
	+
	\|\nabla f(\theta_{n} ) \| 
	\|\theta_{a(n,\hat{t} ) } - \theta_{n} \| 
	\leq &
	\tilde{C}_{4}
	(\zeta_{n}^{2} 
	+ 
	\zeta_{n} \|\nabla f(\theta_{n} ) \| 
	+ 
	\hat{t}^{2} \|\nabla f(\theta_{n} ) \|^{2} )
	\nonumber \\
	\leq &
	2^{-1} \hat{C}_{1}^{-1} \hat{t} 
	(\gamma_{n}^{-2r} 
	+ 
	\gamma_{n}^{-r} \|\nabla f(\theta_{n} ) \| 
	+ 
	\|\nabla f(\theta_{n} ) \|^{2} )
	\nonumber \\
	\leq &
	\hat{C}_{1}^{-1} \hat{t} 
	(\gamma_{n}^{-2r} + \|\nabla f(\theta_{n} ) \|^{2} )
	\nonumber \\
	\leq &
	2 \hat{C}_{1}^{-1} \hat{t} \gamma_{n}^{-2r} 
	+
	(f(\theta_{a(n,\hat{t} ) } ) - f(\theta_{n} ) )
\end{align*}
for $n>\sigma$. 
Therefore, 
\begin{align}\label{l3.2.29}
	2 (f(\theta_{n} ) - f(\theta_{a(n,\hat{t} ) } ) )
	+
	\|\nabla f(\theta_{n} ) \|
	\|\theta_{a(n,\hat{t} ) } - \theta_{n} \| 
	\leq 
	2 \hat{C}_{1}^{-1} \hat{t} \gamma_{n}^{-2r} 
\end{align}
for $n>\sigma$. 
Then, (\ref{l3.2.1*}) -- (\ref{l3.2.7*}) directly follow from 
(\ref{l3.2.23}), (\ref{l3.2.25}), (\ref{l3.2.27}), (\ref{l3.2.29}). 
\end{IEEEproof}

\begin{lemma} \label{lemma3.3}
Suppose that Assumptions \ref{a1} -- \ref{a4} hold. 
Let $\hat{C}_{2} = 4\hat{p} \hat{M}^{2}$
(notice that $1\leq \hat{C}_{2} < \infty$ everywhere). 
Then, there exists 
an integer-valued random variable 
$\tau$ 
such that 
$0\leq \tau < \infty$ everywhere and such that 
\begin{align}
	& \label{l3.3.1*}
	\left(
	u(\theta_{a(n,\hat{t} ) } ) - u(\theta_{n} ) 
	+
	(\hat{t}/4) \|\nabla f(\theta_{n} ) \|^{2} 
	\right)
	I_{A_{n} } 
	\leq 
	0, 
	\\
	& \label{l3.3.3*}
	\left(
	u(\theta_{a(n,\hat{t} ) } ) - u(\theta_{n} ) 
	+
	(\hat{t}/\hat{C}_{2} ) u(\theta_{n} )  
	\right)
	I_{B_{n} } 
	\leq 
	0, 
	\\
	& \label{l3.3.5*}
	\left(
	v(\theta_{a(n,\hat{t} ) } ) - v(\theta_{n} ) 
	-
	\hat{t}/\hat{C}_{2}   
	\right)
	I_{C_{n} } 
	\geq 
	0
\end{align}
on $\Lambda \setminus N_{0}$ for $n>\tau$, 
where
\begin{align*}
	A_{n}
	=&
	\{\gamma_{n}^{\hat{p} } |u(\theta_{n} ) | \geq 1 \}
	\cup
	\{\gamma_{n}^{r } \|\nabla f(\theta_{n} ) \| \geq 1 \}, 
	\\
	B_{n}
	=&
	\{\gamma_{n}^{\hat{p} } u(\theta_{n} ) \geq 1 \}
	\cap
	\{\hat{\mu} = 2 \}, 
	\\
	C_{n}
	=&
	\{\gamma_{n}^{\hat{p} } u(\theta_{n} ) \geq 1 \}
	\cap
	\{u(\theta_{a(n,\hat{t} ) } ) > 0 \}
	\cap
	\{\hat{\mu} < 2 \}
\end{align*}
($\hat{t}$ is specified in the statement of Lemma \ref{lemma3.2}). 
\end{lemma}

\begin{remark}
Inequalities (\ref{l3.3.1*}) -- (\ref{l3.3.5*}) 
can be interpreted in the following way: 
Relations 
\begin{align}
	&\label{l3.3.1*'}
	\left(
	\gamma_{n}^{\hat{p} } |u(\theta_{n} ) | \geq 1 
	\: \vee \:
	\gamma_{n}^{r } \|\nabla f(\theta_{n} ) \| \geq 1 
	\right)
	\: \wedge \: 
	n>\tau
%	\nonumber\\
%	&
	\: \Longrightarrow \: 
	u(\theta_{a(n,\hat{t} ) } ) - u(\theta_{n} ) 
	\leq
	-(\hat{t}/4) \|\nabla f(\theta_{n} ) \|^{2},  
	\\
	&\label{l3.3.3*'}
	\gamma_{n}^{\hat{p} } u(\theta_{n} ) | \geq 1 
	\: \wedge \:
	\hat{\mu} = 2
	\: \wedge \: 
	n>\tau
%	\nonumber\\
%	&
	\: \Longrightarrow \: 
	u(\theta_{a(n,\hat{t} ) } ) 
	\leq 
	(1 - \hat{t}/\hat{C}_{2} ) u(\theta_{n} ),  
	\\
	&\label{l3.3.5*'}
	\gamma_{n}^{\hat{p} } u(\theta_{n} ) \geq 1 
	\: \wedge \:
	\hat{\mu} < 2 
	\: \wedge \: 
	n>\tau
%	\nonumber\\
%	&
	\: \Longrightarrow \: 
	v(\theta_{a(n,\hat{t} ) } ) - v(\theta_{n} ) 
	\geq
	\hat{t}/\hat{C}_{2}   
\end{align}
are true on $\Lambda \setminus N_{0}$. 
\end{remark}

\begin{IEEEproof}
Let 
\begin{align*}
	&
	\tilde{\tau}_{1}
	=
	\max\left(
	\left\{n\geq 0: \theta_{n}\not\in \hat{Q} \right\}
	\cup
	\{0\}
	\right), 
	\\
	&
	\tilde{\tau}_{2}
	=
	\max\left(
	\left\{n\geq 0: |u(\theta_{n} ) | > \hat{\delta} \right\}
	\cup
	\{0\}
	\right)
\end{align*}
and 
$\tau = \max\{\sigma, \tilde{\tau}_{1}, \tilde{\tau}_{2} \} 
I_{\Lambda \setminus N_{0} }$.
Then, it is obvious that 
$\tau$ is well-defined, 
while Lemma \ref{lemma2.4} implies 
$0\leq \tau < \infty$ everywhere. 
On the other side, since $\tau \geq \sigma$
on $\Lambda \setminus N_{0}$, 
Lemma \ref{lemma3.2} (inequality (\ref{l3.2.5*})) implies 
\begin{align}\label{l3.3.1}
	u(\theta_{a(n,\hat{t} ) } )
	-
	u(\theta_{n} )
	\leq 
	-
	(\hat{t}/2) \|\nabla f(\theta_{n} ) \|^{2} 
	+
	(\hat{t}/\hat{C}_{1} ) \gamma_{n}^{-2r}
\end{align}
on $\Lambda \setminus N_{0}$
for $n>\tau$. 
As $\theta_{n} \in \hat{Q}$, $|u(\theta_{n} | \leq \hat{\delta}$
on $\Lambda \setminus N_{0}$ for $n>\tau$, 
(\ref{1*.3.7}) (i.e., Corollary \ref{corollary1}) yields 
\begin{align}\label{l3.3.3}
	|u(\theta_{n} ) |
	\leq 
	\hat{M} \|\nabla f(\theta_{n} ) \|^{\hat{\mu} }
\end{align}
on $\Lambda \setminus N_{0}$ for $n>\tau$. 

Let $\omega$ be an arbitrary sample from $\Lambda \setminus N_{0}$
(notice that all formulas which follow in the proof correspond to this $\omega$). 
First, we show (\ref{l3.3.1*}). 
We proceed by contradiction: 
Suppose that 
(\ref{l3.3.1*}) is violated for some $n>\tau$. 
Consequently, 
\begin{align}\label{l3.3.5}
	u(\theta_{a(n,\hat{t} ) } )
	-
	u(\theta_{n} )
	+
	(\hat{t}/4) \|\nabla f(\theta_{n} ) \|^{2} 
	>
	0
\end{align}
and at least one of the following two 
inequalities is true: 
\begin{align}\label{l3.3.7}
	|u(\theta_{n} ) | \geq \gamma_{n}^{-\hat{p} }, 
	\;\;\;\;\; 
	\|\nabla f(\theta_{n} ) \| \geq \gamma_{n}^{-r }.
\end{align}
If 
$|u(\theta_{n} ) | \geq \gamma_{n}^{-\hat{p} }$, 
then (\ref{l3.3.3}) implies 
\begin{align*}
	\|\nabla f(\theta_{n} ) \|^{2}
	\geq 
	\left(
	|u(\theta_{n} ) |/\hat{M} 
	\right)^{2/\hat{\mu} }
	\geq 
	(1/\hat{M} )^{2/\hat{\mu} } \gamma_{n}^{-2\hat{p}/\hat{\mu} }
	\geq 
	(4/\hat{C}_{1} ) \gamma_{n}^{-2r}
\end{align*}
(notice that $\hat{p}/\hat{\mu} \leq r$, 
$4\hat{M}^{2/\hat{\mu} } \leq 4 \hat{M}^{2} \leq \hat{C}_{1}$). 
Thus, as a result of one of (\ref{l3.3.7}), we get
\begin{align*}
	\|\nabla f(\theta_{n} ) \|^{2}
	\geq 
	(4/\hat{C}_{1} ) \gamma_{n}^{-2r}, 
\end{align*}
i.e., 
$(\hat{t}/4) \|\nabla f(\theta_{n} ) \|^{2}
\geq (\hat{t}/\hat{C}_{1} ) \gamma_{n}^{-2r}$. 
Then, (\ref{l3.3.1}) implies 
\begin{align}\label{l3.3.9}
	u(\theta_{a(n,\hat{t} ) } )
	-
	u(\theta_{n} )
	\leq
	-(\hat{t}/4) \|\nabla f(\theta_{n} ) \|^{2}, 
\end{align}
which directly contradicts (\ref{l3.3.5}). 
Hence, (\ref{l3.3.1*}) is true for $n>\tau$. 
Owing to this, (\ref{l3.3.3}) and the fact that 
$B_{n} \subset A_{n}$ for $n\geq 0$, we obtain
\begin{align*}
	\left(
	u(\theta_{a(n,\hat{t} ) } )
	-
	u(\theta_{n} )
	+
	(\hat{t}/\hat{C}_{2} ) 
	u(\theta_{n} )
	\right)
	I_{B_{n} }
	\leq &
	\left(
	u(\theta_{a(n,\hat{t} ) } )
	-
	u(\theta_{n} )
	+
	(\hat{M} \hat{t}/\hat{C}_{2} ) 
	\|\nabla f(\theta_{n} ) \|^{2}
	\right)
	I_{B_{n} }
	\\
	\leq &
	\left(
	u(\theta_{a(n,\hat{t} ) } )
	-
	u(\theta_{n} )
	+
	(\hat{t}/4 ) 
	\|\nabla f(\theta_{n} ) \|^{2}
	\right)
	I_{B_{n} }
	\leq 
	0
\end{align*}
for $n>\tau$
(notice that $u(\theta_{n} ) > 0$ on $B_{n}$; 
also notice that $4\hat{M} \leq \hat{C}_{2}$). 
Thus, (\ref{l3.3.3*}) is satisfied. 

Now, let us prove (\ref{l3.3.5*}). 
To do so, we again use contradiction: 
Suppose that (\ref{l3.3.5*}) does not hold for some 
$n>\tau$. 
Consequently, we have $\hat{\mu} < 2$, 
$u(\theta_{a(n,\hat{t} ) } ) > 0$
and 
\begin{align}
	&\label{l3.3.21}
	\gamma_{n}^{\hat{p} } u(\theta_{n} ) 
	\geq 
	1, 
	\\
	&\label{l3.3.23}
	v(\theta_{a(n,\hat{t} ) } )
	-
	v(\theta_{n} )
	<
	\hat{t}/\hat{C}_{2}. 
\end{align}
Combining (\ref{l3.3.21}) with 
(already proved) (\ref{l3.3.1*}), 
we get (\ref{l3.3.9}). 
On the other side, (\ref{l3.3.3}) yields 
\begin{align*}
	\|\nabla f(\theta_{n} ) \|^{2} 
	\geq 
	\left( 
	u(\theta_{n} ) / \hat{M}
	\right)^{2/\hat{\mu} }
	\geq 
	\hat{M}^{-2} 
	(u(\theta_{n} ) )^{1 + 1/\hat{p} }
\end{align*}
(notice that $0<u(\theta_{n} ) \leq \hat{\delta} \leq 1$, 
$2/\hat{\mu} = 1 + 1/(\hat{\mu} \hat{r} ) \leq 1 + 1/\hat{p}$). 
Therefore, (\ref{l3.3.9}) implies 
\begin{align*}
	\frac{\hat{t} }{4} 
	\leq 
	\frac{u(\theta_{n} ) - u(\theta_{a(n,\hat{t} ) } ) }
	{\|\nabla f(\theta_{n} ) \|^{2} }
	\leq &
	\hat{M}^{2} \; 
	\frac{u(\theta_{n} ) - u(\theta_{a(n,\hat{t} ) } ) }
	{(u(\theta_{n} ) )^{1 + 1/\hat{p} } }
	\\
	= &
	\hat{M}^{2} 
	\int_{u(\theta_{a(n,\hat{t} ) } ) }^{u(\theta_{n} ) }
	\frac{du}{(u(\theta_{n} ) )^{1 + 1/\hat{p} } }
	\\
	\leq &
	\hat{M}^{2} 
	\int_{u(\theta_{a(n,\hat{t} ) } ) }^{u(\theta_{n} ) }
	\frac{du}{u^{1+1/\hat{p} } }
	\\
	=&
	\frac{\hat{C}_{2} }{4}
	(v(\theta_{a(n,\hat{t} ) } ) - v(\theta_{n} ) ). 
\end{align*}
Thus, 
$
	v(\theta_{a(n,\hat{t} ) } )
	-
	v(\theta_{n} )
	\geq 
	\hat{t}/\hat{C}_{2}, 
$
which directly contradicts (\ref{l3.3.23}). 
Hence, (\ref{l3.3.3*}) is satisfied for $n>\tau$. 
\end{IEEEproof}

\begin{lemma} \label{lemma3.4}
Suppose that Assumptions \ref{a1} -- \ref{a4} hold. 
Then, 
\begin{align}	
	& \label{l3.4.1*} 
	\gamma_{n}^{\hat{p} }
	u(\theta_{n} ) 
	\geq 
	- 1, 
	\\
	& \label{l3.4.3*} 
	\|\nabla f(\theta_{n} ) \|^{2}
	\leq 
	(4/\hat{t} )
	\left(
	\varphi(u(\theta_{n} ) )
	+
	\gamma_{n}^{-\hat{p} }  
	\right)
\end{align}
on $\Lambda\setminus N_{0}$ for $n > \tau$, 
where
function $\varphi(\cdot )$
is defined by 
$\varphi(x) = x \:{\rm I}_{(0,\infty )}(x)$, $x\in \mathbb{R}$. 
\end{lemma}

\begin{IEEEproof}
Let 
$\omega$ be an arbitrary sample from $\Lambda\setminus N_{0}$
(notice that all formulas that follow in the proof correspond to 
this $\omega$). 
First, 
we prove (\ref{l3.4.1*}). 
To do so, we use contradiction:  
Assume that (\ref{l3.4.1*}) is not satisfied for 
some 
$n_{0} > \tau$, 
and  
define recursively  
$n_{k+1} = a(n_{k},\hat{t} )$ 
for $k\geq 0$.  
Now, let us show by induction that 
$\{u(\theta_{n_{k} } ) \}_{k\geq 0}$ is non-increasing: 
Suppose that 
$u(\theta_{n_{l} } ) \leq u(\theta_{n_{l-1} } )$ for 
$0\leq l \leq k$ and some $k\geq 1$. 
Consequently, 
\begin{align*}
	u(\theta_{n_{k} } ) 
	\leq 
	u(\theta_{n_{0} } ) 
	\leq 
	-
	\gamma_{n_{0} }^{-\hat{p} }
	\leq 
	-
	\gamma_{n_{k} }^{-\hat{p} }. 
\end{align*}
Then, Lemma \ref{lemma3.3} (relations (\ref{l3.3.1*}), (\ref{l3.3.1*'})) yields 
\begin{align*}
	u(\theta_{n_{k+1} } ) - u(\theta_{n_{k} } ) 
	\leq 
	-
	(\hat{t}/4) \|\nabla f(\theta_{n_{k} } ) \|^{2} 
	\leq 
	0,  
\end{align*}
i.e., $u(\theta_{n_{k+1} } ) \leq u(\theta_{n_{k} } )$. 
Thus, $\{u(\theta_{n_{k} } ) \}_{k\geq 0}$ is non-increasing. 
Therefore,
\begin{align*}
	\limsup_{n\rightarrow \infty } 
	u(\theta_{n_{k} } ) 
	\leq 
	u(\theta_{n_{0} } ) 
	< 
	0. 
\end{align*}
However, this is not possible, as 
$\lim_{n\rightarrow \infty } u(\theta_{n} ) = 0$
(due to Lemma \ref{lemma2.4}). 
Hence, (\ref{l3.4.1*}) indeed holds 
for $n> \tau$. 

Now, (\ref{l3.4.3*}) is demonstrated. 
Again, we proceed by contradiction: 
Suppose that (\ref{l3.4.3*}) is violated for some 
$n> \tau$. 
Consequently, 
\begin{align*}
	\|\nabla f(\theta_{n} ) \|^{2}
	\geq 
	(4/\hat{t} )
	\gamma_{n}^{-\hat{p} }
	\geq
	\gamma_{n}^{-2r} 
\end{align*}
(notice that $\hat{p} \leq \hat{\mu} r \leq 2r$),
which, together with Lemma \ref{lemma3.3} 
(relations (\ref{l3.3.1*}), (\ref{l3.3.1*'})), 
yields 
\begin{align*}
	u(\theta_{a(n,\hat{t} ) } ) - u(\theta_{n} )
	\leq 
	-(\hat{t}/4) \|\nabla f(\theta_{n} ) \|^{2}. 
\end{align*}
Then, (\ref{l3.4.1*}) implies 
\begin{align*}
	\|\nabla f(\theta_{n} ) \|^{2}
	\leq &
	(4/\hat{t} )
	\left(
	u(\theta_{n} ) 
	-
	u(\theta_{a(n,\hat{t} ) } ) 
	\right)
	\leq 
	(4/\hat{t} )
	\left(
	\varphi(u(\theta_{n} ) )
	+
	\gamma_{n}^{-\hat{p} }
	\right). 
\end{align*}
However, this directly contradicts our assumption 
that $n$ violates (\ref{l3.4.3*}). 
Thus, (\ref{l3.4.3*}) is satisfied for 
$n>\tau$. 
\end{IEEEproof}

\begin{lemma} \label{lemma3.5}
Suppose that Assumptions \ref{a1} -- \ref{a4} hold. 
Let $\hat{C}_{3} = 2 \hat{C}_{2}^{\hat{p} }$. 
Then, 
\begin{align} \label{l3.5.1*}
	\liminf_{n\rightarrow \infty } 
	\gamma_{n}^{\hat{p} }
	u(\theta_{n} ) 
	\leq 
	\hat{C}_{3}
\end{align}
on $\Lambda\setminus N_{0}$. 
\end{lemma}

\begin{IEEEproof}
We prove the lemma by contradiction: 
Assume that (\ref{l3.5.1*}) is violated for some 
sample $\omega$ from $\Lambda\setminus N_{0}$
(notice that the formulas which follow in the proof correspond to 
this $\omega$). 
Consequently, there exists 
$n_{0} > \tau$ 
such that 
\begin{align} \label{l3.5.1}
	\gamma_{n}^{\hat{p} } u(\theta_{n} ) 
	\geq 
	\hat{C}_{3}
\end{align}
for $n\leq n_{0}$. 

Let $\{n_{k} \}_{k\geq 0}$ be defined recursively as 
$n_{k+1} = a(n_{k},\hat{t} )$ for $k\geq 0$. 	
In what follows in the proof, we consider separately 
the cases $\hat{\mu} < 2$
and $\hat{\mu} = 2$. 

{\em Case $\hat{\mu} < 2$:}
Owing to Lemma \ref{lemma3.3} 
(relations (\ref{l3.3.5*}), (\ref{l3.3.5*'})) and (\ref{l3.5.1}), we have 
\begin{align*}
	v(\theta_{n_{k+1} } ) - v(\theta_{n_{k} } ) 
	\geq &
	\hat{t}/\hat{C}_{2}
	\geq
	(\gamma_{n_{k+1} } - \gamma_{n_{k} } ) / \hat{C}_{2}		
\end{align*}
for $k\geq 0$
(notice that $\gamma_{n}^{\hat{p} } u(\theta_{n} ) \geq 1$ due to 
(\ref{l3.5.1}); 
also notice that $\gamma_{n_{k+1} } - \gamma_{n_{k} } \leq \hat{t}$).  
Therefore, 
\begin{align*}
	v(\theta_{n_{k} } )
	\geq 
	v(\theta_{n_{0} } )
	+
	(1/\hat{C}_{2} )    
	\sum_{i=0}^{k-1} (\gamma_{n_{i+1} }  - \gamma_{n_{i} } )
	= 
	v(\theta_{n_{0} } )
	+
	(\gamma_{n_{k} }  - \gamma_{n_{0} } )/\hat{C}_{2}
\end{align*}
for $k\geq 0$. 
Then, (\ref{l3.5.1}) implies
\begin{align*}
	\left(
	v(\theta_{n_{0} } )/\gamma_{n_{k} }
	+
	(1 - \gamma_{n_{0} }/\gamma_{n_{k} } )/\hat{C}_{2}
	\right)^{-\hat{p} } 
	\geq 	
	(v(\theta_{n_{k} } ) / \gamma_{n_{k} } )^{-\hat{p} }
	= 
	\gamma_{n_{k} }^{\hat{p} }
	u(\theta_{n_{k} } )
	\geq 
	\hat{C}_{3} 
\end{align*}
for $k\geq 0$. 
However, this is impossible, since the limit process 
$k\rightarrow \infty$ (applied to the previous relation)  
yields 
$\hat{C}_{3} \leq \hat{C}_{2}^{\hat{p} }$. 
Hence, (\ref{l3.5.1*}) holds when $\hat{\mu} < 2$. 

{\em Case $\hat{\mu} = 2$:} 
Due to Lemma \ref{lemma3.3} (relations (\ref{l3.3.3*}), (\ref{l3.3.3*'})) and 
(\ref{l3.5.1}), we have 
\begin{align*}
	u(\theta_{n_{k+1} } )
	\leq
	(1 - \hat{t}/\hat{C}_{2} )
	u(\theta_{n_{k} } )
	\leq 
	\left(1 - (\gamma_{n_{k+1} } - \gamma_{n_{k} } ) / \hat{C}_{2} \right)
	u(\theta_{n_{k} } )
\end{align*}
for $k\geq 0$. 
Consequently, 
\begin{align*}
	u(\theta_{n_{k} } )
	\leq &
	u(\theta_{n_{0} } )	
	\prod_{i=0}^{k-1} 
	\left(
	1	-	(\gamma_{n_{i+1} } - \gamma_{n_{i} } )/\hat{C}_{2}
	\right)
	\\
	\leq &
	u(\theta_{n_{0} } )	
	\exp\left(
	-
	(1/\hat{C}_{2} ) 
	\sum_{i=0}^{k-1} (\gamma_{n_{i+1} } - \gamma_{n_{i} } )	
	\right)
	\\
	= &
	u(\theta_{n_{0} } )	
	\exp\left( 
	-(\gamma_{n_{k} } - \gamma_{n_{0} } )/\hat{C}_{2} 	
	\right)
\end{align*}
for $k\geq 0$. 
Then, (\ref{l3.5.1}) yields 
\begin{align*}
	u(\theta_{n_{0} } )	
	\gamma_{n_{k} }^{\hat{p} }
	\exp\left(
	-(\gamma_{n_{k} } - \gamma_{n_{0} } )/\hat{C}_{2}	
	\right)
	\geq 
	\gamma_{n_{k} }^{\hat{p} } u(\theta_{n_{k} } )
	\geq 
	\hat{C}_{3}	
\end{align*}
for $k\geq 0$. 
However, this is not possible, as the limit 
process $k\rightarrow \infty$
(applied to the previous relation) 
implies 
$\hat{C}_{3}\leq 0$. 
Thus, (\ref{l3.5.1*}) holds in the case $\hat{\mu} = 2$, too. 
\end{IEEEproof}

\begin{lemma} \label{lemma3.6}
Suppose that Assumptions \ref{a1} -- \ref{a4} hold. 
Let $\hat{C}_{4} = 6 \hat{C}_{3}$. 
Then, 
\begin{align} \label{l3.6.1*}
	\limsup_{n\rightarrow \infty } 
	\gamma_{n}^{\hat{p} }
	u(\theta_{n} ) 
	\leq 
	\hat{C}_{4} 
\end{align}
on $\Lambda\setminus N_{0}$. 
\end{lemma}

\begin{IEEEproof}
We use contradiction to prove the lemma: 
Suppose that (\ref{l3.6.1*}) is violated for some sample 
$\omega$ from $\Lambda\setminus N_{0}$
(notice that the formulas which appear in the proof correspond to 
this $\omega$). 
Since 
$\lim_{n\rightarrow \infty } (\gamma_{a(n,\hat{t} ) }/\gamma_{n} ) = 1$, 
it can be deduced from Lemma \ref{lemma3.5} that 
there exist 
$n_{0} > m_{0} > \tau$ such that 
\begin{align}
	& \label{l3.6.1}
	\gamma_{m_{0} }^{\hat{p} } u(\theta_{m_{0} } ) 
	\leq 
	2\hat{C}_{3}, 
	\\
	& \label{l3.6.3}
	\gamma_{n_{0} }^{\hat{p} } u(\theta_{n_{0} } ) 
	>
	\hat{C}_{4}, 
	\\
	& \label{l3.6.5}
	\min_{m_{0}<n\leq n_{0} }
	\gamma_{n}^{\hat{p} } u(\theta_{n} ) 
	>
	2\hat{C}_{3}, 
	\\
	&
	\max_{m_{0} \leq n < n_{0} } 
	\gamma_{n}^{\hat{p} } u(\theta_{n} ) 
	\leq 
	\hat{C}_{4}, 
\end{align}
and such that 
\begin{align}
	(\gamma_{a(m_{0},\hat{t} ) } /\gamma_{m_{0} } )^{\hat{p} }
	\leq 
	\min\{2, (1 - \hat{t}/\hat{C}_{2} )^{-1} \}. 
\end{align}
Let $l_{0} = a(m_{0},\hat{t} )$. 
As a direct consequence of Lemma \ref{lemma3.4} and 
(\ref{l3.6.1}), we get
\begin{align*}
	\|\nabla f(\theta_{m_{0} } ) \|^{2} 
	\leq 
	(4/\hat{t} )
	\left(
	\varphi(u(\theta_{m_{0} } ) ) + \gamma_{m_{0} }^{-\hat{p} } 
	\right)
	\leq 
	12(\hat{C}_{3}/\hat{t} ) \gamma_{m_{0} }^{-\hat{p} }.
\end{align*}
Consequently, Lemma \ref{lemma3.2} and (\ref{1*.3.3}) imply 
\begin{align*}
	u(\theta_{n} ) - u(\theta_{m_{0} } ) 
	\leq 
	|\phi_{m_{0},n } |
	\leq &
	(\hat{t}/\hat{C}_{1} )
	(\gamma_{m_{0} }^{-2r} + \|\nabla f(\theta_{m_{0} } ) \|^{2} )
	\\
	\leq &
	(\hat{t}/\hat{C}_{1} ) \gamma_{m_{0} }^{-2r} 
	+
	(12\hat{C}_{3}/\hat{C}_{1} ) \gamma_{m_{0} }^{-\hat{p} }
	\leq 
	\gamma_{m_{0} }^{-\hat{p} } 
\end{align*}
for $m_{0}\leq n \leq l_{0}$
(notice that $\hat{p} \leq 2r$, $\hat{t}/\hat{C}_{1} \leq 1/2$, 
$\hat{C}_{1} \geq 24\hat{C}_{3}$). 
Then, (\ref{l3.6.1}), (\ref{l3.6.5}) yield 
\begin{align}
	& \label{l3.6.7}
	u(\theta_{m_{0} } ) 
	\geq 
	u(\theta_{m_{0} + 1 } )
	-
	\gamma_{m_{0} }^{-\hat{p} }
	\geq 
	2\hat{C}_{3}
	(\gamma_{{m}_{0} }/\gamma_{m_{0}+1} )^{\hat{p} }
	\gamma_{m_{0} }^{-\hat{p} } 
	- 
	\gamma_{m_{0} }^{-\hat{p} }
	\geq
	(\hat{C}_{3} - 1 ) \gamma_{m_{0} }^{-\hat{p} }
	\geq 
	\gamma_{m_{0} }^{-\hat{p} }, 
	\\
	& \label{l3.6.9}
	u(\theta_{n} )
	\leq 
	u(\theta_{m_{0} } )
	+
	\gamma_{m_{0} }^{-\hat{p} }
	\leq 
	(2\hat{C}_{3} + 1 ) 
	(\gamma_{n}/\gamma_{m_{0} } )^{\hat{p} } \gamma_{n}^{-\hat{p} }
	\leq 
	6\hat{C}_{3} \gamma_{n}^{-\hat{p} }
	=
	\hat{C}_{4} \gamma_{n}^{-\hat{p} }
\end{align}
for $m_{0} \leq n \leq l_{0}$
(notice that 
$(\gamma_{n}/\gamma_{m_{0} } )^{\hat{p} } \leq
(\gamma_{l_{0} }/\gamma_{m_{0} } )^{\hat{p} } \leq 2$ for 
$m_{0} \leq n \leq n_{0}$). 
Using (\ref{l3.6.3}), (\ref{l3.6.9}), we conclude
$l_{0} < n_{0}$. 

In the rest of the proof, we consider separately the cases 
$\hat{\mu} < 2$ and $\hat{\mu} = 2$. 

{\em Case $\hat{\mu} < 2$:}
Owing to Lemma \ref{lemma3.3} (relations (\ref{l3.3.5*}), (\ref{l3.3.5*'})) and 
(\ref{l3.6.1}), (\ref{l3.6.7}), 
we have 
\begin{align*}
	v(\theta_{l_{0} } ) 
	\geq 
	v(\theta_{m_{0} } ) 
	+
	\hat{t}/\hat{C}_{2} 
	\geq &
	(2\hat{C}_{3} )^{-1/\hat{p} } \:
	\gamma_{m_{0} }
	+
	(\gamma_{l_{0} } - \gamma_{m_{0} } )/\hat{C}_{2}
	\\
	> &
	\min\{(2\hat{C}_{3} )^{-1/\hat{p} }, 
	\hat{C}_{2}^{-1} \}
	\gamma_{l_{0} }
	\\
	=&
	(2\hat{C}_{3} )^{-1/\hat{p} } \:
	\gamma_{l_{0} }	
\end{align*}
(notice that 
$(2\hat{C}_{3} )^{1/\hat{p} } > \hat{C}_{2}$).
Therefore, 
\begin{align*}
	u(\theta_{l_{0} } ) 
	=
	(v(\theta_{l_{0} } ) )^{-\hat{p} }
	< 
	2\hat{C}_{3} \gamma_{l_{0} }^{-\hat{p} }. 
\end{align*}
However, this directly contradicts 
(\ref{l3.6.5}) and the fact that 
$m_{0} < l_{0} < n_{0}$. 
Thus, (\ref{l3.6.1*}) holds 
when $\hat{\mu} < 2$. 

{\em Case $\hat{\mu} = 2$:}
Using Lemma \ref{lemma3.3} (relations (\ref{l3.3.3*}),(\ref{l3.3.3*'})) and (\ref{l3.6.7}), 
we get
\begin{align*}
	u(\theta_{l_{0} } )
	\leq 
	(1 - \hat{t}/\hat{C}_{2} ) u(\theta_{{m}_{0} } )
	\leq 
	2\hat{C}_{3} (1 - \hat{t}/\hat{C}_{2} )
	(\gamma_{l_{0} }/\gamma_{m_{0} } )^{\hat{p} }
	\gamma_{l_{0} }^{-\hat{p} }
	\leq 
	2\hat{C}_{3} \gamma_{l_{0} }^{-\hat{p} }. 
\end{align*}
However, this is impossible due to 
(\ref{l3.6.5}) and the fact that 
$m_{0} < l_{0} < n_{0}$. 
Hence, 
(\ref{l3.6.1*}) holds in the case $\hat{\mu} = 2$, too. 
\end{IEEEproof}

\begin{lemma} \label{lemma3.7}
Suppose that Assumptions \ref{a1} -- \ref{a4} hold. 
Then, 
\begin{align} \label{l3.7.1*}
	\|\theta_{a(n,\hat{t} ) } - \theta_{n} \|
	\leq 
	2 \gamma_{n}^{\hat{q} + 1 }
	(u(\theta_{n} ) - u(\theta_{a(n,\hat{t} } ) )
	+
	6\gamma_{n}^{-(\hat{q} + 1 ) }
\end{align}
on $\Lambda \setminus N_{0}$ for $n> \tau$. 
\end{lemma}

\begin{IEEEproof}
Let $\omega$ be an arbitrary sample from 
$\Lambda \setminus N_{0}$, 
while $n > \max\{\sigma, \tau \}$ is an arbitrary integer
(notice that all formulas which appear in the proof
correspond to these $\omega$, $n$). 
To show (\ref{l3.7.1*}), 
we consider separately the cases
$\|\nabla f(\theta_{n} ) \| \geq \gamma_{n}^{-(\hat{q} + 1 ) }$
and  
$\|\nabla f(\theta_{n} ) \| < \gamma_{n}^{-(\hat{q} + 1 ) }$. 

{\em Case $\|\nabla f(\theta_{n} ) \| \geq \gamma_{n}^{-(\hat{q} + 1 ) }$:}
Due to Lemma \ref{lemma3.2}, 
we have
\begin{align}\label{l3.7.1}
	\|\nabla f(\theta_{n} ) \|
	\|\theta_{a(n,\hat{t} ) } - \theta_{n} \|
	\leq 
	2(u(\theta_{n} ) - u(\theta_{a(n,\hat{t} ) } ) )
	+
	2	(\hat{t}/\hat{C}_{1} ) \gamma_{n}^{-2r}.   
\end{align}
On the other side, 
since 
$\|\nabla f(\theta_{n} ) \| \geq \gamma_{n}^{-(\hat{q}+1)} \geq \gamma_{n}^{-r}$
(notice that 
$\hat{q}+1 = \min\{(\hat{p}+1)/2,r\} \leq r$), 
Lemma \ref{lemma3.3} (relations (\ref{l3.3.1*}), (\ref{l3.3.1*'})) 
implies 
\begin{align*}
	u(\theta_{a(n,\hat{t} ) } )
	-
	u(\theta_{n} )
	\leq 
	- (\hat{t}/4 ) \|\nabla f(\theta_{n} ) \|^{2} 
	<
	0, 
\end{align*}
i.e., $u(\theta_{n} ) - u(\theta_{a(n,\hat{t} ) } ) > 0$. 
Then, (\ref{l3.7.1}) yields 
\begin{align*}
	\|\theta_{a(n,\hat{t} ) } - \theta_{n} \|
	\leq &
	2(u(\theta_{n} ) - u(\theta_{a(n,\hat{t} ) } ) )
	\|\nabla f(\theta_{n} ) \|^{-1} 
	+
	2 (\hat{t}/\hat{C}_{1} ) \gamma_{n}^{-2r} \|\nabla f(\theta_{n} ) \|^{-1} 
	\\
	\leq &
	2 \gamma_{n}^{\hat{q} + 1 }
	(u(\theta_{n} ) - u(\theta_{a(n,\hat{t} ) } ) )
	+
	\gamma_{n}^{-2r + (\hat{q} + 1 ) }
	\\
	\leq &
	2 \gamma_{n}^{\hat{q} + 1 }
	(u(\theta_{n} ) - u(\theta_{a(n,\hat{t} ) } ) )
	+
	\gamma_{n}^{-(\hat{q} + 1 ) }
\end{align*}
(notice that $\hat{t}/\hat{C}_{1} \leq 1/2$; 
also notice that $\hat{q}+1 \leq r$, which implies 
$2r - (\hat{q} + 1 ) \geq \hat{q} + 1$). 
Hence, (\ref{l3.7.1*}) is true when
$\|\nabla f(\theta_{n} ) \| \geq \gamma_{n}^{-(\hat{q} + 1 ) }$. 

{\em Case $\|\nabla f(\theta_{n} ) \| < \gamma_{n}^{-(\hat{q}+1) }$:}
Using Lemma \ref{lemma3.2} and (\ref{1*.3.3}), we get
\begin{align*}
	|u(\theta_{a(n,\hat{t} ) } ) - u(\theta_{n} ) |
	\leq &
	(\gamma_{a(n,\hat{t} ) } - \gamma_{n} )
	\|\nabla f(\theta_{n} ) \|^{2} 
	+ 
	|\phi_{n,a(n,\hat{t} ) } |
	\\
	\leq &
	\hat{t} \|\nabla f(\theta_{n} ) \|^{2} 
	+
	(\hat{t}/\hat{C}_{1} )(\gamma_{n}^{-2r} + \|\nabla f(\theta_{n} ) \|^{2} ) 
	\\
	\leq &
	2\gamma_{n}^{-2(\hat{q} + 1 ) } 
\end{align*}
(notice that $\hat{q} + 1 \leq r < 2r$ and $\hat{t}/\hat{C}_{1} \leq 1/2$). 
On the other side, owing to Lemma \ref{lemma3.2} and (\ref{1*.3.1}), 
we have 
\begin{align*}
	\|\theta_{a(n,\hat{t} ) } - \theta_{n} \|
	\leq &
	(\gamma_{a(n,\hat{t} ) } - \gamma_{n} )
	\|\nabla f(\theta_{n} ) \| 
	+ 
	\|\varepsilon_{n,a(n,\hat{t} ) } \|
	\\
	\leq &
	\hat{t} \|\nabla f(\theta_{n} ) \| 
	+
	(\hat{t}/\hat{C}_{1} )(\gamma_{n}^{-r} + \|\nabla f(\theta_{n} ) \| ) 
	\\
	\leq &
	2\gamma_{n}^{-(\hat{q} + 1 ) } 
\end{align*}
(notice that $\hat{q} + 1 \leq r$). 
Consequently, 
\begin{align*}
	\|\theta_{a(n,\hat{t} ) } - \theta_{n} \|
	\leq &
	2\gamma_{n}^{\hat{q} + 1 }
	(u(\theta_{n} ) - u(\theta_{a(n,\hat{t} ) } ) )
	+
	2\gamma_{n}^{\hat{q} + 1 }
	|u(\theta_{n} ) - u(\theta_{a(n,\hat{t} ) } ) |
	+
	2\gamma_{n}^{-(\hat{q} + 1 ) }
	\\
	\leq &
	2\gamma_{n}^{\hat{q} + 1 }
	(u(\theta_{n} ) - u(\theta_{a(n,\hat{t} } ) )
	+
	6\gamma_{n}^{-(\hat{q} + 1 ) }. 
\end{align*}
Thus, (\ref{l3.7.1*}) holds in the case 
$\|\nabla f(\theta_{n} ) \| < \gamma_{n}^{-(\hat{q}+1)}$. 
\end{IEEEproof}

\begin{lemma} \label{lemma3.8}
Suppose that Assumptions \ref{a1} -- \ref{a4} hold. 
Then, there exists a random quantity $\hat{C}_{5}$
such that $1\leq \hat{C}_{5} < \infty$ everywhere 
and such that 
\begin{align} \label{l3.8.1*}
	\limsup_{n\rightarrow \infty } 
	\gamma_{n}^{\hat{q} } \max_{k\geq n} \|\theta_{k} - \theta_{n} \| 
	\leq 
	\hat{C}_{5}
\end{align}
on $\Lambda \setminus N_{0}$. 
\end{lemma}

\begin{IEEEproof}
Let $\tilde{C} = 9\hat{C}_{4} (\hat{q} + 1 )$ 
and 
$\hat{C}_{5} = 20 \tilde{C} \hat{t}^{-1} (1 + 1/\hat{q} )$, 
while $\omega$ is an arbitrary sample from 
$\Lambda \setminus N_{0}$
(notice that all formulas which follow in the proof correspond to 
this $\omega$). 

As a consequence of Lemmas \ref{lemma3.4} and \ref{lemma3.6}, we get
\begin{align} 
	&\label{l3.8.1}
	\limsup_{n\rightarrow \infty } \gamma_{n}^{\hat{p} } 
	|u(\theta_{n} ) |
	\leq  
	\hat{C}_{4}, 
	\\
	&\label{l3.8.3}
	\limsup_{n\rightarrow \infty } 
	\gamma_{n}^{\hat{p} } \|\nabla f(\theta_{n} ) \|^{2} 
	\leq 
	8 \hat{C}_{4}/\hat{t}. 
\end{align}
Since $\gamma_{a(n,\hat{t} ) } - \gamma_{n} = \hat{t} + O(\alpha_{a(n,\hat{t} ) } )$
for $n\rightarrow \infty$, 
 and 
\begin{align*}
	(1 - \hat{t}/\gamma_{n} )^{\hat{q} + 1 } 
	=
	1 - \hat{t} (\hat{q} + 1 ) \gamma_{n}^{-1} + o(\gamma_{n}^{-1} )
\end{align*}
for $n\rightarrow \infty$, 
we conclude from 
(\ref{l3.8.1}), (\ref{l3.8.3}) that there exists 
$n_{0} > \max\{\sigma, \tau \}$ (depending on $\omega$)
such that 
$|u(\theta_{n} ) | \leq 2\hat{C}_{4} \gamma_{n}^{-\hat{p} }$, 
$\|\nabla f(\theta_{n} ) \| \leq (4\hat{C}_{4}/\hat{t} ) \gamma_{n}^{-\hat{p}/2 }$, 
$\gamma_{a(n,\hat{t} ) } - \gamma_{n} \geq \hat{t}/2$
and
\begin{align} \label{l3.8.3'}
	(1 - \hat{t}/\gamma_{n} )^{\hat{q} + 1}
	\geq 
	1 - (\hat{q} + 1 ) \gamma_{n}^{-1} 
\end{align}
for $n\geq n_{0}$. 
Then, (\ref{1*.3.1}) and Lemma \ref{lemma3.2} imply 
\begin{align}\label{l3.8.5}
	\|\theta_{k} - \theta_{n} \|
	\leq &
	(\gamma_{k} - \gamma_{n} ) \|\nabla f(\theta_{n} ) \|
	+
	\|\varepsilon_{n,k} \|
	\nonumber\\
	\leq &
	\hat{t} \|\nabla f(\theta_{n} ) \| 
	+
	(\hat{t}/\hat{C}_{1} )
	(\gamma_{n}^{-r} + \|\nabla f(\theta_{n} ) \| )
	\nonumber\\
	\leq &
	8 \hat{C}_{4} \gamma_{n}^{-\hat{p}/2} 
	+
	\gamma_{n}^{-r}
	\nonumber\\
	\leq &
	\tilde{C} \gamma_{n}^{-\hat{q} }
\end{align}
for $n_{0} \leq n \leq k \leq a(n,\hat{t} )$
(notice that $\hat{q} < \min\{\hat{p}/2,r\}$). 

Let $\{n_{k} \}_{k\geq 0}$ be recursively defined as 
$n_{k+1} = a(n_{k}, \hat{t} )$ for $k\geq 0$. 
Due to Lemma \ref{lemma3.7}, we have 
\begin{align} \label{l3.8.7}
	\|\theta_{n_{l} } - \theta_{n_{k} } \|
	\leq 
	\sum_{i=k}^{l-1} 
	\|\theta_{n_{i+1} } - \theta_{n_{i} } \|
	\leq &
	6 \sum_{i=k}^{l-1} \gamma_{n_{i} }^{-(\hat{q} + 1 ) }
	+ 
	2 \sum_{i=k}^{l-1} 
	\gamma_{n_{i} }^{\hat{q} + 1 } 
	(u(\theta_{n_{i} } ) - u(\theta_{n_{i+1} } ) )
	\nonumber\\
	\leq &
	6 \sum_{i=k}^{l-1} \gamma_{n_{i} }^{-(\hat{q} + 1 ) }
	+ 
	2 \sum_{i=k+1}^{l} 
	(\gamma_{n_{i} }^{\hat{q} + 1 } - \gamma_{n_{i-1} }^{\hat{q} + 1 } )
	|u(\theta_{n_{i} } ) |
	\nonumber\\
	&
	+
	2 \gamma_{n_{l} }^{\hat{q} + 1}
	|u(\theta_{n_{l} } ) |
	+
	2 \gamma_{n_{k} }^{\hat{q} + 1}
	|u(\theta_{n_{k} } ) |
\end{align}
for $l\geq k \geq 0$. 
As 
\begin{align*}
	\gamma_{n_{i} }^{\hat{q} + 1 } - \gamma_{n_{i-1} }^{\hat{q} + 1 } 
	=
	\gamma_{n_{i} }^{\hat{q} + 1 } 
	\left(
	1
	-
	\left(
	1
	-
	(\gamma_{n_{i} } - \gamma_{n_{i-1} } )/\gamma_{n_{i} }
	\right)^{\hat{q} + 1}
	\right)
	\leq 
	\gamma_{n_{i} }^{\hat{q} + 1 } 
	\left(
	1
	-
	\left(
	1
	-
	\hat{t}/\gamma_{n_{i} }
	\right)^{\hat{q} + 1}
	\right)
	\leq 
	(\hat{q} + 1 )
	\gamma_{n_{i} }^{\hat{q} } 
\end{align*}
for $i\geq 0$
(use (\ref{l3.8.3'})), we get 
\begin{align} \label{l3.8.9}
	\sum_{i=k+1}^{l} 
	(\gamma_{n_{i} }^{\hat{q}+1} -\gamma_{n_{i-1} }^{\hat{q}+1} )
	|u(\theta_{n_{i} } ) |
	\leq 
	2 \hat{C}_{4} (\hat{q} + 1 ) 
	\sum_{i=k}^{\infty } 
	\gamma_{n_{i} }^{-\hat{p}+\hat{q} }
	\leq
	\tilde{C} 
	\sum_{i=k}^{\infty } 
	\gamma_{n_{i} }^{-(\hat{q}+1 ) }
\end{align}
for $l> k \geq 0$
(notice that $\hat{p} - \hat{q} \geq (\hat{p}+1 )/2 \geq \hat{q} + 1$). 
Since 
\begin{align*}
	\gamma_{n_{l} }
	=
	\gamma_{n_{k} }
	+
	\sum_{i=k}^{l-1} 
	(\gamma_{n_{i+1} } - \gamma_{n_{i} } )
	\geq 
	\gamma_{n_{k} }
	+
	(\hat{t}/2) (l-k)
\end{align*}
for $l>k\geq 0$
(notice that $\gamma_{a(n,\hat{t} ) } - \gamma_{n} \geq \hat{t}/2)$ for $n\geq n_{0}$), 
we have
\begin{align*}
	\sum_{i=k}^{\infty } \gamma_{n_{i} }^{-(\hat{q}+1) }
	\leq &
	\sum_{i=0}^{\infty } 
	(\gamma_{n_{k} } + \hat{t} i/2 )^{-(\hat{q} + 1 ) }
	\\
	\leq &
	\gamma_{n_{k} }^{-(\hat{q}+1) }
	+
	\int_{0}^{\infty } (\gamma_{n_{k} } + \hat{t} u /2 )^{-(\hat{q}+1) } du
	\\
	= &
	\gamma_{n_{k} }^{-(\hat{q}+1) }
	+
	2 \hat{t}^{-1} \hat{q}^{-1} \gamma_{n_{k} }^{-\hat{q} }
	\\
	\leq &
	(1 + 2 \hat{t}^{-1} \hat{q}^{-1} ) \gamma_{n_{k} }^{-\hat{q} }
\end{align*}
for $k\geq 0$. 
Consequently, (\ref{l3.8.7}) and (\ref{l3.8.9}) imply 
\begin{align} \label{l3.8.11}
	\|\theta_{n_{l} } - \theta_{n_{k} } \|
	\leq 
	(6 + 2 \tilde{C} )
	\sum_{i=k}^{\infty } 
	\gamma_{n_{i} }^{-(\hat{q} + 1 ) }
	+
	4\hat{C}_{4} \gamma_{n_{k} }^{-\hat{p} + \hat{q} + 1 }
	+
	4\hat{C}_{4} \gamma_{n_{l} }^{-\hat{p} + \hat{q} + 1 }
	\leq 
	16\tilde{C} (1 + \hat{t}^{-1} \hat{q}^{-1} )
	\gamma_{n_{k} }^{-\hat{q} }
\end{align}
for $l\geq k \geq 0$
(notice that 
$\hat{p} - (\hat{q} + 1 ) \geq (\hat{p}-1)/2 \geq \hat{q}$). 
Using (\ref{l3.8.5}) and (\ref{l3.8.11}), 
we get
\begin{align*}
	\|\theta_{k} - \theta_{n} \|
	\leq &
	\|\theta_{k} - \theta_{n_{j} } \|
	+
	\|\theta_{n_{j} } - \theta_{n_{i} } \|
	+
	\|\theta_{n_{i} } - \theta_{n} \|
	\\
	\leq &
	\tilde{C} \gamma_{k}^{-\hat{q} } 
	+
	\tilde{C} \gamma_{n}^{-\hat{q} } 
	+
	16\tilde{C} (1 + \hat{t}^{-1} \hat{q}^{-1} ) 
	\gamma_{n_{i} }^{-\hat{q} }
	\\
	\leq &
	\hat{C}_{5} \gamma_{n}^{-\hat{q} } 
\end{align*}
for $k\geq n \geq n_{0}$, 
$j\geq i \geq 1$
satisfying 
$n_{i-1}\leq n < n_{i}$, $n_{j-1}\leq k < n_{j}$. 
Then, it is obvious that (\ref{l3.8.1*}) is true. 
\end{IEEEproof}

\begin{IEEEproof}[Proof of Theorems \ref{theorem2} and \ref{theorem3}]
Owing to Lemmas \ref{lemma2.4} and \ref{lemma3.8}, 
we have that on $\Lambda \setminus N_{0}$, 
$\hat{\theta } = \lim_{n\rightarrow \infty } \theta_{n}$
exists and satisfies $\nabla f(\hat{\theta } ) = 0$. 
Consequently, 
$\hat{Q} \subseteq 
\{\theta \in \mathbb{R}^{d_{\theta} }: \|\theta - \hat{\theta } \| 
\leq \delta_{\hat{\theta} } \}$
on $\Lambda \setminus N_{0}$. 
Thus, random quantities $\hat{p}$, $\hat{q}$ defined in this subsection
coincide with $\hat{p}$, $\hat{q}$ introduced in 
Theorem \ref{theorem3}
(see the remark after Corollary \ref{corollary1}). 
Then, Lemmas \ref{lemma3.4}, \ref{lemma3.6}, \ref{lemma3.8}
imply that  
(\ref{t3.1*}) is true on $\Lambda \setminus N_{0}$. 
\end{IEEEproof}

\section{Proof of Propositions \ref{proposition1} -- \ref{proposition4}} \label{section2*}

\renewcommand{\theenumi}{\alph{enumi}}

\begin{IEEEproof}[Proof of Proposition \ref{proposition1}]
Owing to Conditions (i), (ii) of the proposition, 
for any compact set $Q\subset \Theta$, there exists a real number 
$\varepsilon_{Q} \in (0,1)$ such that 
\begin{align} \label{p1.1}
	\varepsilon_{Q} 
	\leq 
	r_{\theta }(y|x',x) 
	\leq 
	\varepsilon_{Q}^{-1}
\end{align}
for all $\theta \in Q$, $x,x' \in {\cal X}$, $y \in {\cal Y}$. 
Hence, Assumption \ref{a3} is satisfied. 
On the other side, Condition (ii) implies that 
$r_{\theta}(y|x',x)$ has an (complex-valued) analytical continuation 
$\hat{r}_{\eta}(y|x',x)$
with the following properties: 
\begin{enumerate}
\item
$\hat{r}_{\eta}(y|x',x)$ maps 
$(\eta,x,x',y) \in \mathbb{C}^{d_{\theta } } \times {\cal X} \times {\cal X} \times {\cal Y}$
into $\mathbb{C}$. 
\item
$\hat{r}_{\theta }(y|x',x) = r_{\theta }(y|x',x)$ for all $\theta \in \Theta$, 
$x,x' \in {\cal X}$, $y\in {\cal Y}$. 
\item
For any compact set $Q\subset \Theta$, 
there exists a real number 
$\tilde{\delta}_{Q} \in (0,1)$ such that 
$\hat{r}_{\eta}(y|x',x)$ is analytical in $\eta$ on 
$V_{\tilde{\delta}_{Q} }(Q)$ 
for each $x,x' \in {\cal X}$, $y \in {\cal Y}$. 
\end{enumerate}
Relying on $\hat{r}_{\eta}(y|x',x)$, 
we define quantities 
$\hat{R}_{\eta}(y)$, $\hat{\phi}_{\eta}(w,y)$, $\hat{G}_{\eta}(w,y)$. 
More specifically, 
for $\eta \in \mathbb{C}^{d_{\theta } }$, $y \in {\cal Y}$,  
$\hat{R}_{\eta}(y)$ is an $N_{x}\times N_{x}$ matrix 
whose $(i,j)$ entry is 
$\hat{r}_{\eta}(y|i,j)$, while  
\begin{align}
	&\label{p1.3}
	\hat{\phi}_{\eta}(w,y) 
	=
	\begin{cases}
	\log(e^{T} \hat{R}_{\eta}(y) w ), 
	& \text{ if } e^{T} \hat{R}_{\eta}(y) w \neq 0
	\\
	0, & \text{ otherwise}
	\end{cases}
	\\
	&\label{p1.5}
	\hat{G}_{\eta}(w,y) 
	=
	\begin{cases}
	\hat{R}_{\eta}(y) w /(e^{T} \hat{R}_{\eta}(y) w ), 
	& \text{ if } e^{T} \hat{R}_{\eta}(y) w \neq 0
	\\
	0, & \text{ otherwise}
	\end{cases}
\end{align}
for $\eta\in \mathbb{C}^{d_{\theta } }$, $y\in {\cal Y}$, $w\in \mathbb{C}^{N_{x} }$. 

Let $Q\subset\Theta$ be an arbitrary compact set. 
Since
$e^{T} R_{\theta }(y) u \geq N_{x} \varepsilon_{Q}$
for all $\theta \in Q$, $y\in {\cal Y}$, $u \in {\cal P}^{N_{x} }$
(due to (\ref{p1.1})), 
we conclude that 
there exists a real number $\delta_{Q} \in (0, \tilde{\delta}_{Q} )$
such that 
$|e^{T} \hat{R}_{\eta }(y) w | \geq N_{x} \varepsilon_{Q} /2$
for all $\eta \in V_{\delta_{Q} }(Q)$, 
$w\in V_{\delta_{Q} }({\cal P}^{N_{x} } )$, $y \in {\cal Y}$. 
Therefore, 
$\hat{\phi}_{\eta}(w,y)$, $\hat{G}_{\eta}(w,y)$ are 
analytical in $(\eta,w)$ on 
$V_{\delta_{Q} }(Q) \times V_{\delta_{Q} }({\cal P}^{N_{x} } )$
for any $y\in {\cal Y}$. 
Consequently, 
$|\hat{\phi}_{\eta}(w,y)|$, $\|\hat{G}_{\eta}(w,y)\|$ are uniformly bounded in 
$(\eta,w,y)$ on 
$V_{\delta_{Q} }(Q) \times V_{\delta_{Q} }({\cal P}^{N_{x} } )
\times {\cal Y}$. Thus, Assumption \ref{a4} is satisfied, too. 
\end{IEEEproof}

\begin{IEEEproof}[Proof of Proposition \ref{proposition2}]
Conditions (i), (ii) of the proposition imply that for any compact set 
$Q\subset\Theta$, 
there exists a real number $\varepsilon_{Q}\in (0,1)$ such that 
$\varepsilon_{Q} \leq r_{\theta }(y|x',x) \leq \varepsilon_{Q}^{-1}$
for all $\theta \in Q$, $x,x' \in {\cal X}$, $y\in {\cal Y}$. 
Thus, Assumption \ref{a3} holds. 
On the other side, as a result of Condition (ii), 
$r_{\theta }(y|x',x)$ has an (complex-valued) analytical continuation 
$\hat{r}_{\eta}(z|x',x)$ with the following properties: 
\begin{enumerate}
\item
$\hat{r}_{\eta}(z|x',x)$ maps 
$(\eta,x,x',z) \in \mathbb{C}^{d_{\theta } } \times {\cal X} \times {\cal X} \times \mathbb{C}^{d_{y} }$
into $\mathbb{C}$. 
\item
$\hat{r}_{\theta}(y|x',x) = r_{\theta}(y|x',x)$ 
for all $\theta \in \Theta$, $x,x' \in {\cal X}$, $y \in {\cal Y}$. 
\item
For any compact set $Q\subset\Theta$, there exists a real number 
$\tilde{\delta}_{Q} \in (0,1)$ such that 
$\hat{r}_{\eta}(z|x',x)$ is analytical in $(\eta, z )$
on $V_{\tilde{\delta}_{Q} }(Q) \times V_{\tilde{\delta}_{Q} }({\cal Y})$
for each $x,x' \in {\cal X}$. 
\end{enumerate}
Relying on $\hat{r}_{\eta}(y|x',x)$, 
we define quantities 
$\hat{R}_{\eta}(y)$, $\hat{\phi}_{\eta}(w,y)$, $\hat{G}_{\eta}(w,y)$
in the same way as in the proof of Proposition \ref{proposition1}. 
More specifically, 
for $\eta \in \mathbb{C}^{d_{\theta } }$, $y \in {\cal Y}$, 
$\hat{R}_{\eta}(y)$ is an $N_{x} \times N_{x}$ matrix whose 
$(i,j)$ entry is $\hat{r}_{\eta}(y|i,j)$, 
while $\hat{\phi}_{\eta}(w,y)$, $\hat{G}_{\eta}(w,y)$ are defined 
by (\ref{p1.3}), (\ref{p1.5})
for $\eta\in \mathbb{C}^{d_{\theta } }$, $y \in {\cal Y}$, $w\in \mathbb{C}^{N_{x} }$. 

Let $Q\subset \Theta$ be an arbitrary compact set. 
As 
$N_{x} \varepsilon_{Q} \leq e^{T} R_{\theta}(y) u \leq N_{x} \varepsilon_{Q}^{-1}$
for any $\theta\in Q$, $y\in {\cal Y}$, 
$u\in {\cal P}^{N_{x} }$,
we have that there exists a real number 
$\delta_{Q} \in (0, \tilde{\delta}_{Q} )$ such that 
$N_{x} \varepsilon_{Q}/2 \leq 
|e^{T} \hat{R}_{\eta}(y) w | \leq 2 N_{x} \varepsilon_{Q}^{-1}$
for all $\eta \in V_{\delta_{Q} }(Q)$, 
$w \in V_{\delta_{Q} }({\cal P}^{N_{x} } )$, $y\in {\cal Y}$
(notice that $|e^{T} \hat{R}_{\eta}(y) w |$ is analytical in 
$(\eta,w,y)$ on 
$V_{\tilde{\delta}_{Q} }(Q) \times 
V_{\tilde{\delta}_{Q} }({\cal P}^{N_{x} } ) \times {\cal Y}$). 
Therefore, 
$\hat{\phi}_{\eta}(w,y)$, $\hat{G}_{\eta}(w,y)$ are analytical in 
$(\eta,w)$ on 
$V_{\delta_{Q} }(Q) \times V_{\delta_{Q} }({\cal P}^{N_{x} } )$
for any $y\in {\cal Y}$. 
Moreover, 
$|\hat{\phi}_{\eta}(w,y)|$, $\|\hat{G}_{\eta}(w,y)\|$ are uniformly bounded in 
$(\eta,w,y)$ on 
$V_{\delta_{Q} }(Q) \times V_{\delta_{Q} }({\cal P}^{N_{x} } )
\times {\cal Y}$. 
Hence, Assumption \ref{a4} holds, too. 
\end{IEEEproof}

\begin{IEEEproof}[Proof of Proposition \ref{proposition3}]
For $\alpha \in {\cal A}$, 
$\beta = [\beta_{1} \cdots \beta_{N_{\beta } } ]^{T} \in {\cal B}$, 
$x,x' \in {\cal X}$, 
let
$
	g_{\theta }^{k}(x'|x)
	=
	\beta_{x',k} p_{\alpha}(x'|x)
$. 
Then, we have
\begin{align*}
	r_{\theta }(y|x',x) 
	=
	\sum_{k=1}^{N_{\beta } } 
	f_{k}(y|x') g_{\theta}^{k}(x'|x)
\end{align*}
for all $\theta \in \Theta$, $x,x' \in {\cal X}$, $y \in {\cal Y}$. 
We also have that for any compact set $Q\subset\Theta$, 
there exists a real number $\varepsilon_{Q} \in (0,1)$ such that 
$
	\varepsilon_{Q} 
	\leq 
	g_{\theta }^{k}(x'|x) 
	\leq 
	\varepsilon_{Q}^{-1}
$
for each $\theta \in Q$, $x,x' \in {\cal X}$,
$1\leq k \leq N_{\beta }$. 
Consequently, 
\begin{align*}
	\varepsilon_{Q} 
	\sum_{k=1}^{N_{\beta} }
	f_{k}(y|x') 
	\leq 
	r_{\theta }(y|x',x) 
	\leq 
	\varepsilon_{Q}^{-1} 
	\sum_{k=1}^{N_{\beta} }
	f_{k}(y|x') 
\end{align*}
for all $\theta \in Q$, $x,x' \in {\cal X}$ 
and any compact set $Q\subset\Theta$. 
Hence, Assumption \ref{a3} holds 
(set $s_{\theta }(y|x) = \sum_{k=1}^{N_{\beta} } f_{k}(y|x)$). 
On the other side, Condition (i) implies that for each $1\leq k \leq N_{\beta}$, 
$g_{\theta }^{k}(x'|x)$ has an (complex-valued) analytical continuation 
$\hat{g}_{\eta}^{k}(x'|x)$ with the following properties: 
\begin{enumerate}
\item
$\hat{g}_{\eta}(x'|x)$ maps $(\eta,x,x') \in \mathbb{C}^{d_{\theta } }\times {\cal X} \times {\cal X}$
into $\mathbb{C}$. 
\item
$\hat{g}_{\theta }^{k}(x'|x) = g_{\theta }^{k}(x'|x)$ for all 
$\theta \in \Theta$, $x,x' \in {\cal X}$. 
\item
For any compact set $Q\subset\Theta$, 
there exists a real number $\tilde{\delta}_{Q} \in (0,1)$ such that 
$\hat{g}_{\eta}^{k}(x'|x)$ is analytical in $\eta$ on 
$V_{\tilde{\delta}_{Q} }(Q)$ for each $x,x' \in {\cal X}$. 
\end{enumerate}
Relying on $\hat{g}_{\eta}^{k}(x'|x)$, we define some new quantities. 
More specifically, 
for $\eta \in \mathbb{C}^{d_{\theta } }$, 
$w = [w_{1} \cdots w_{N_{x} } ]^{T} \in \mathbb{C}^{N_{x} }$, 
$x,x' \in {\cal X}$, $y \in {\cal Y}$, let 
\begin{align*}
	&
	\hat{r}_{\eta}(y|x',x) 
	=
	\sum_{k=1}^{N_{\beta} } 
	f_{k}(y|x') \hat{g}_{\eta}^{k}(x'|x), 
	\\
	&
	\hat{h}_{\eta,w}^{k}(x') 
	=
	\sum_{x''\in {\cal X} } \hat{g}_{\eta}^{k}(x'|x'') w_{x''}, 
\end{align*}
while $\hat{R}_{\eta}(y)$ is an $N_{x} \times N_{x}$ matrix whose 
$(i,j)$ entry is $\hat{r}_{\eta}(y|i,j)$. 
Moreover, let $\hat{\phi}_{\eta}(w,y)$, $\hat{G}_{\eta}(w,y)$ be defined 
for $\eta\in \mathbb{C}^{d_{\theta } }$, $w\in \mathbb{C}^{N_{x} }$, $y\in {\cal Y}$ in the same 
way as in (\ref{p1.3}), (\ref{p1.5}). 

Let $Q\subset\Theta$ be arbitrary compact set. 
Since 
\begin{align*}
	\varepsilon_{Q} 
	\leq 
	\sum_{x\in {\cal X} } g_{\theta }^{k}(x'|x) u_{x} 
	\leq 
	\varepsilon_{Q}^{-1}
\end{align*}
for all $\theta \in Q$, $u = [u_{1} \cdots u_{N_{x} } ]^{T} \in {\cal P}^{N_{x} }$, 
$x,x' \in {\cal X}$, $1 \leq k \leq N_{\beta}$, 
we deduce that 
there exists a real number $\delta_{Q} \in (0,\tilde{\delta}_{Q} )$ such that 
${\rm Re}\{\hat{h}_{\eta,w}^{k}(x') \} \geq \varepsilon_{Q}/2$, 
$|\hat{h}_{\eta,w}^{k}(x') | \leq 2\varepsilon_{Q}^{-1}$
for all $\eta \in V_{\delta_{Q} }(Q)$, 
$w \in V_{\delta_{Q} }({\cal P}^{N_{x} } )$, $x' \in {\cal X}$, 
$1 \leq k \leq N_{\beta }$. 
Consequently, 
\begin{align*}
	&
	|e^{T}\hat{R}_{\eta}(y) w |
	\geq 
	|{\rm Re}\{e^{T}\hat{R}_{\eta}(y) w \}|
	=
	\sum_{x'\in {\cal X} }\sum_{k=1}^{N_{\beta} } 
	f_{k}(y|x') {\rm Re}\{\hat{h}_{\eta,w}^{k}(x') \}
	\geq 
	(\varepsilon_{Q}/2 ) \psi(y)
	> 0, 
	\\
	&
	\max\{\|\hat{R}_{\eta}(y) w \|, |e^{T}\hat{R}_{\eta}(y) w | \}
	\leq 
	\sum_{x'\in {\cal X} }\sum_{k=1}^{N_{\beta} } 
	f_{k}(y|x') |\hat{h}_{\eta,w}^{k}(x') |
	\leq 
	2 \varepsilon_{Q}^{-1} \psi(y)
\end{align*}
for all 
$\eta \in V_{\delta_{Q} }(Q)$, 
$w \in V_{\delta_{Q} }({\cal P}^{N_{x} } )$, $y\in {\cal Y}$. 
Therefore, 
$\hat{\phi}_{\eta}(w,y)$, $\hat{G}_{\eta}(w,y)$ are analytical in 
$(\eta,w)$ on 
$V_{\delta_{Q} }(Q) \times V_{\delta_{Q} }({\cal P}^{N_{x} } )$
for each $y\in {\cal Y}$. 
Moreover, 
\begin{align*}
	&
	\|\hat{G}_{\eta}(w,y) \|
	\leq 
	4 \varepsilon_{Q}^{-2}, 
	\\
	&
	|\hat{\phi}_{\eta}(w,y) |
	\leq 
	|\log|e^{T} \hat{R}_{\eta}(y) w | | + 2\pi
	\leq 
	|\log\psi(y) | 
	+ 
	\log(2\varepsilon_{Q}^{-1} ) 
	+ 
	2\pi
\end{align*}
for all 
$\eta \in V_{\delta_{Q} }(Q)$, 
$w \in V_{\delta_{Q} }({\cal P}^{N_{x} } )$, $y\in {\cal Y}$. 
Then, it is clear that Assumption \ref{a4} holds, too. 
\end{IEEEproof}

\renewcommand{\theenumi}{\roman{enumi}}
\begin{lemma}\label{lemmaP}
Let the conditions of Proposition \ref{proposition4} hold. 
Then, $\phi_{\theta}(u,y)$, $G_{\theta }(u,y)$ have
(complex-valued) analytical continuations 
$\hat{\phi}_{\eta}(w,y)$, $\hat{G}_{\eta}(w,y)$ (respectively)
with the following properties: 
\begin{enumerate}
\item
$\hat{\phi}_{\eta}(w,y)$, $\hat{G}_{\eta}(w,y)$ map 
$(\eta,w,y) \in \mathbb{C}^{d_{\theta } } \times \mathbb{C}^{N_{x} } \times {\cal Y}$
into $\mathbb{C}$, $\mathbb{C}^{N_{x} }$ (respectively). 
\item
$\hat{\phi}_{\theta}(u,y) = \phi_{\theta}(u,y)$, 
$\hat{G}_{\theta}(u,y) = G_{\theta}(u,y)$
for all $\theta \in \Theta$, $u\in {\cal P}^{N_{x} }$, $y \in {\cal Y}$. 
\item
For each $\theta \in \Theta$, there exist real numbers 
$\delta_{\theta } \in (0,1)$, $K_{\theta } \in [1,\infty )$ such that 
$\hat{\phi}_{\eta}(w,y)$, $\hat{G}_{\eta}(w,y)$ are analytical in 
$(\eta, w )$ on 
$V_{\delta_{\theta } }(\theta ) \times V_{\delta_{\theta } }({\cal P}^{N_{x} } )$
for any $y\in {\cal Y}$, and such that 
\begin{align*}
	&
	|\hat{\phi}_{\eta}(w,y) |
	\leq 
	K_{\theta }(1 + y^{2} ), 
	\\
	&
	\|\hat{G}_{\eta}(w,y) \| 
	\leq 
	K_{\theta }
\end{align*}
for all $\eta \in V_{\delta_{\theta } }(\theta )$, 
$w \in V_{\delta_{\theta } }({\cal P}^{N_{x} } )$, 
$y\in {\cal Y}$. 
\end{enumerate}
\end{lemma}

\renewcommand{\theenumi}{\alph{enumi}}
\begin{IEEEproof}
Due to Condition (i) of Proposition \ref{proposition4}, 
$p_{\alpha}(x'|x)$ has an (complex-valued) analytical continuation 
$\hat{p}_{a}(x'|x)$
with the following properties 
\begin{enumerate}
\item
$\hat{p}_{a}(x'|x)$ maps $(a,x,x') \in \mathbb{C}^{d_{\alpha} } \times {\cal X} \times {\cal X}$
into $\mathbb{C}$. 
\item
$\hat{p}_{\alpha}(x'|x) = p_{\alpha}(x'|x)$ for all $\alpha \in {\cal A}$, 
$x,x' \in {\cal X}$. 
\item
For any $\alpha \in {\cal A}$, 
there exists a real number $\tilde{\delta}_{\alpha} \in (0,1)$
such that $\hat{p}_{a}(x'|x)$ is analytical in $a$ on 
$V_{\tilde{\delta}_{\alpha} }(\alpha )$ for each $x,x' \in {\cal X}$. 
\end{enumerate}
On the other side, the analytical continuation 
$\hat{q}_{b}(y|x)$ of $q_{\beta}(y|x)$ is defined by 
\begin{align*}
	\hat{q}_{b}(y|x)
	=
	\sqrt{l_{x}/\pi }
	\exp(-l_{x} (y-m_{x} )^{2} ), 
\end{align*}
for $b = [l_{1} \cdots l_{N_{x} } \: m_{1} \cdots m_{N_{x} } ]^{T} \in \mathbb{C}^{2N_{x} }$, 
$x\in {\cal X}$, $y \in {\cal Y}$.  

Let  
$\hat{r}_{\eta}(y|x',x) = \hat{q}_{b}(y|x') \hat{p}_{a}(x'|x)$ for 
$a \in \mathbb{C}^{d_{\alpha} }$, $b \in \mathbb{C}^{2N_{x} }$, 
$\eta = [a^{T} \: b^{T} ]^{T}$, $x,x' \in {\cal X}$, $y \in {\cal Y}$. 
Moreover, for $\eta \in \mathbb{C}^{d_{\theta } }$, $y \in {\cal Y}$, 
$\hat{R}_{\eta}(y)$ is an $N_{x} \times N_{x}$ matrix whose 
$(i,j)$ entry is $\hat{r}_{\eta}(y|i,j)$, 
while 
$\hat{\phi}_{\eta}(w,y)$, $\hat{G}_{\eta}(w,y)$ are defined for 
$\eta \in \mathbb{C}^{d_{\theta } }$, $w\in \mathbb{C}^{N_{x} }$, $y \in {\cal Y}$
in the same way as in (\ref{p1.3}), (\ref{p1.5}). 

Let $\alpha$, 
$\beta = 
[\lambda_{1} \cdots \lambda_{N_{x} } \: \mu_{1} \: \mu_{N_{x} } ]^{T}$
be arbitrarily vectors from ${\cal A}$, ${\cal B}$ (respectively), 
while 
$\theta = [\alpha^{T} \: \beta^{T} ]^{T}$. 
Obviously, it can be assumed without loss of generality that 
$0 < \lambda_{1} < \lambda_{2} < \cdots < \lambda_{N_{x} }$.
Since  
\begin{align*}
	\sum_{x\in {\cal X} } p_{\alpha}(x'|x) u_{x} > 0
\end{align*}
for all $x'\in {\cal X}$, $u = [u_{1} \cdots u_{N_{x} } ]^{T} \in {\cal P}^{N_{x} }$, 
there exist real numbers 
$\tilde{\delta}_{1,\theta }, \tilde{\varepsilon}_{\theta } \in (0,1)$ such that 
$\hat{R}_{\eta}(y)$ is analytical in $\eta$ on 
$V_{\tilde{\delta}_{1,\theta } }(\theta )$ for any $y \in {\cal Y}$, 
and such that 
\begin{align}
	&\label{lP.1}
	{\rm Re}\left\{
	\sum_{x\in {\cal X} } 
	\hat{p}_{a}(x'|x) w_{x} 
	\right\}
	\geq 
	\tilde{\varepsilon}_{\theta }, 
	\\
	&\label{lP.3}
	\left|
	\sum_{x\in {\cal X} } 
	\hat{p}_{a}(x'|x) w_{x} 
	\right|
	\leq 
	\tilde{\varepsilon}_{\theta}^{-1}, 
	\\
	&
	\min\{{\rm Re}\{l_{1} \}, {\rm Re}\{l_{x'} - l_{1} \} \}
	\geq 
	\tilde{\varepsilon}_{\theta }, 
	\nonumber\\
	&
	\max\{|l_{x''} |, |m_{x''} | \}
	\leq 
	\tilde{\varepsilon}_{\theta }^{-1}
	\nonumber
\end{align}
for all $a \in V_{\tilde{\delta}_{1,\theta } }(\alpha )$, 
$b = [l_{1} \cdots l_{N_{x} } \: m_{1} \cdots m_{N_{x} } ]^{T} \in 
V_{\tilde{\delta}_{1,\theta } }(\beta )$, 
$w = [w_{1} \cdots w_{N_{x} } ]^{T} \in V_{\tilde{\delta}_{1,\theta } }({\cal P}^{N_{x} } )$, 
$x' \in {\cal X}\setminus \{1\}$, $x'' \in {\cal X}$.  
Therefore, 
we have 
\begin{align*}
	|\hat{q}_{b}(y|x) |
	= &
	\sqrt{|l_{x}|/\pi}\: 
	|\exp(- {\rm Re}\{l_{x} \} y^{2} + 2 {\rm Re}\{l_{x} m_{x} \} y - {\rm Re}\{l_{x} m_{x}^{2} \} ) |
	\\
	\leq &
	\sqrt{|l_{x}|/\pi}\: 
	\exp(- {\rm Re}\{l_{x} \} y^{2} + 2 |l_{x} | |m_{x} | |y| + |l_{x} | |m_{x} |^{2} )
	\\ 
	\leq &
	(1/\sqrt{\pi\varepsilon_{\theta} } ) \: 
	\exp(-\tilde{\varepsilon}_{\theta } y^{2} 
	+ 
	2 \tilde{\varepsilon}_{\theta }^{-2} |y|
	+
	\tilde{\varepsilon}_{\theta }^{-3} )
\end{align*}
for any $b = [l_{1} \cdots l_{N_{x} } \; m_{1} \cdots m_{N_{x} } ]^{T} \in V_{\tilde{\delta}_{1,\theta } }(\beta )$, 
$x\in {\cal X}$, $y \in {\cal Y}$. 
We also have 
\begin{align*}
	\left|
	\frac{\hat{q}_{b}(y|x) }{\hat{q}_{b}(y|1) }
	\right|
	= &
	\sqrt{|l_{x}|/|l_{1} |}
	|\exp(-{\rm Re}\{l_{x} - l_{1} \} y^{2} 
	+
	2 {\rm Re}\{l_{x} m_{x} - l_{1} m_{1} \} y 
	-
	{\rm Re}\{l_{x} m_{x}^{2} - l_{1} m_{1}^{2} \} 
	) |
	\\
	\leq &
	\sqrt{|l_{x}|/|l_{1} |}
	\exp(-{\rm Re}\{l_{x} - l_{1} \} y^{2} 
	+
	2 (|l_{x} | |m_{x} | + |l_{1} | |m_{1} | ) |y| 
	+
	|l_{x} | |m_{x} |^{2}
	+
	|l_{1} | |m_{1} |^{2} ) 
	\\
	\leq &
	\varepsilon_{\theta}^{-1}
	\exp(-\tilde{\varepsilon}_{\theta } y^{2} 
	+ 
	4 \tilde{\varepsilon}_{\theta }^{-2} |y| 
	+
	2 \tilde{\varepsilon}_{\theta }^{-3} )
\end{align*}
for all 
$b = [l_{1} \cdots l_{N_{x} } \; m_{1} \cdots m_{N_{x} } ]^{T} \in V_{\tilde{\delta}_{1,\theta } }(\beta )$, 
$x\in {\cal X}\setminus \{1\}$, $y \in {\cal Y}$. 
Consequently, there exists a real number 
$\tilde{C}_{\theta } \in [1,\infty )$ such that 
\begin{align} 
	&\label{lP.5}
	\left|
	\frac{\hat{q}_{b}(y|x) }{\hat{q}_{b}(y|1) }
	\right|
	\leq 
	\tilde{C}_{\theta }, 
	\\
	&\label{lP.7}
	|\log|\hat{q}_{b}(y|x) | |
	\leq 
	\tilde{C}_{\theta } (1 + y^{2} )
\end{align}
for all $b \in V_{\tilde{\delta}_{1,\theta } }(\beta )$, 
$x\in {\cal X}$, $y \in {\cal Y}$, 
and such that 
\begin{align} \label{lP.9}
	\left|
	\frac{\hat{q}_{b}(y|x) }{\hat{q}_{b}(y|1) }
	\right|
	\leq 
	2^{-1} N_{x}^{-1} \tilde{\varepsilon}_{\theta }^{2}
\end{align}
for any
$b \in V_{\tilde{\delta}_{1,\theta } }(\beta )$, 
$x\in {\cal X}\setminus \{1\}$, $y \in [-\tilde{C}_{\theta}, \tilde{C}_{\theta } ]^{c}$
(to show that (\ref{lP.9}) holds for all sufficiently large $|y|$, 
notice that 
$\lim_{|y|\rightarrow \infty } \hat{q}_{b}(y|x)/\hat{q}_{b}(y|1) = 0$
for $x\neq 1$). 
As $\hat{q}_{b}(y|x)/q_{\beta}(y|x)$ is uniformly continuous in 
$(b,y)$ on 
$V_{\tilde{\delta}_{1,\theta } }(\beta ) \times [-\tilde{C}_{\theta}, \tilde{C}_{\theta}]$
and 
$\lim_{b\rightarrow \beta} \hat{q}_{b}(y|x)/q_{\beta}(y|x) = 1$
for any $x\in {\cal X}$, $y \in {\cal Y}$, 
there also exists a real number 
$\tilde{\delta}_{2,\theta } \in (0,1)$ such that 
\begin{align}
	&\label{lP.21}
	\left|
	\frac{\hat{q}_{b}(y|x) }{q_{\beta}(y|x) }
	-
	1
	\right|
	\leq 
	2^{-1} \tilde{\varepsilon}_{\theta }^{2}, 
	\\
	&\label{lP.23}
	\left|
	\frac{\hat{q}_{b}(y|x) }{q_{\beta}(y|x) }
	\right|
	\leq 
	2 
\end{align}
for all $b \in V_{\tilde{\delta}_{2,\theta } }(\beta )$, 
$x\in {\cal X}$, $y \in [-\tilde{C}_{\theta}, \tilde{C}_{\theta } ]$. 

Let $\delta_{\theta } = \min\{\tilde{\delta}_{1,\theta }, \tilde{\delta}_{2,\theta } \}$, 
$K_{\theta } = 8 N_{x} \tilde{C}_{\theta } \tilde{\varepsilon}_{\theta }^{-2}$. 
As a result of (\ref{lP.3}), (\ref{lP.5}), we have
\begin{align} \label{lP.31}
	\max\{\|\hat{R}_{\eta}(y) w \|, |e^{T} \hat{R}_{\eta}(y) w | \}
	\leq 
	\sum_{x'\in {\cal X} } 
	|\hat{q}_{b}(y|x') |
	\left|
	\sum_{x\in {\cal X} } \hat{p}_{a}(x'|x) w_{x} 
	\right|
	\leq 
	N_{x} \tilde{C}_{\theta } \tilde{\varepsilon}_{\theta}^{-1} 
	|\hat{q}_{b}(y|1) |
\end{align}
for all $a \in V_{\delta_{\theta } }(\alpha )$, 
$b \in V_{\delta_{\theta } }(\beta )$, 
$\eta = [a^{T} b^{T}]^{T}$, 
$y \in {\cal Y}$, 
$w = [w_{1} \cdots w_{N_{x} } ]^{T} \in V_{\delta_{\theta } }({\cal P}^{N_{x} } )$. 
Using (\ref{lP.1}), (\ref{lP.3}), (\ref{lP.9}), we get
\begin{align} \label{lP.33}
	|e^{T} \hat{R}_{\eta}(y) w |
	= &
	|\hat{q}_{b}(y|1) |
	\left|
	\sum_{x'\in {\cal X} } 
	\frac{\hat{q}_{b}(y|x') }{\hat{q}_{b}(y|1) }
	\sum_{x\in {\cal X} } 
	\hat{p}_{a}(x'|x) w_{x} 
	\right|
	\nonumber \\
	\geq &
	|\hat{q}_{b}(y|1) |
	\left(
	{\rm Re}\left\{
	\sum_{x\in {\cal X} } 
	\hat{p}_{a}(1|x) w_{x}  
	\right\}
	-
	\sum_{x' \in {\cal X}\setminus \{1\} }
	\left|\frac{\hat{q}_{b}(y|x') }{\hat{q}_{b}(y|1) } \right|
	\left|
	\sum_{x\in {\cal X} } 
	\hat{p}_{a}(x'|x) w_{x}  
	\right|
	\right)
	\nonumber \\
	\geq &
	2^{-1} \tilde{\varepsilon}_{\theta } 
	|\hat{q}_{b}(y|1) |
\end{align}
for all $a \in V_{\delta_{\theta } }(\alpha )$, 
$b \in V_{\delta_{\theta } }(\beta )$, 
$\eta = [a^{T} b^{T}]^{T}$, 
$y \in [-\tilde{C}_{\theta }, \tilde{C}_{\theta } ]^{c}$, 
$w = [w_{1} \cdots w_{N_{x} } ]^{T} \in V_{\delta_{\theta } }({\cal P}^{N_{x} } )$. 
Combining (\ref{lP.1}), (\ref{lP.3}), (\ref{lP.21}), (\ref{lP.23}), we obtain 
\begin{align} \label{lP.35}
	|e^{T} \hat{R}_{\eta}(y) w |
	\geq &
	\left|
	\sum_{x'\in {\cal X} } q_{\beta}(y|x') 
	\sum_{x \in {\cal X} } \hat{p}_{a}(x'|x) w_{x} 
	\right|
	-
	\left|
	\sum_{x'\in {\cal X} } (\hat{q}_{b}(y|x') - q_{\beta}(y|x') )
	\sum_{x \in {\cal X} } \hat{p}_{a}(x'|x) w_{x} 
	\right|
	\nonumber\\
	\geq &
	\sum_{x'\in {\cal X} } q_{\beta}(y|x') 
	{\rm Re}\left\{
	\sum_{x \in {\cal X} } \hat{p}_{a}(x'|x) w_{x} 
	\right\}
	-
	\sum_{x'\in {\cal X} } 
	q_{\beta}(y|x') 
	\left|
	\frac{\hat{q}_{b}(y|x') }{q_{\beta}(y|x') } 
	-
	1
	\right|
	\left|
	\sum_{x \in {\cal X} } \hat{p}_{a}(x'|x) w_{x} 
	\right|
	\nonumber\\
	\geq &
	2^{-1} \tilde{\varepsilon}_{\theta } 
	\sum_{x' \in {\cal X} } q_{\beta}(y|x') 
	\nonumber\\
	\geq &
	2^{-1} \tilde{\varepsilon}_{\theta } 
	q_{\beta}(y|1) 
	\nonumber\\
	\geq &
	4^{-1} \tilde{\varepsilon}_{\theta } 
	|\hat{q}_{b}(y|1) |
\end{align}
for any $a \in V_{\delta_{\theta } }(\alpha )$, 
$b \in V_{\delta_{\theta } }(\beta )$, 
$\eta = [a^{T} b^{T}]^{T}$, 
$y \in [-\tilde{C}_{\theta }, \tilde{C}_{\theta } ]$, 
$w = [w_{1} \cdots w_{N_{x} } ]^{T} \in V_{\delta_{\theta } }({\cal P}^{N_{x} } )$. 
Then, it can concluded from (\ref{lP.33}), (\ref{lP.35}) that 
$\hat{\phi}_{\eta}(w,y)$, $\hat{G}_{\eta}(w,y)$ are analytical in 
$(\eta, w )$ on 
$V_{\delta_{\theta } }(\theta ) \times V_{\delta_{\theta } }({\cal P}^{N_{x} } )$
for each $y \in {\cal Y}$. 
On the other side, (\ref{lP.7}), (\ref{lP.31}) -- (\ref{lP.35}) imply 
\begin{align*}
	&
	|\hat{\phi}_{\eta}(w,y) |
	\leq 
	|\log|e^{T} \hat{R}_{\eta}(y) w | |
	+
	2\pi 
	\leq 
	\tilde{C}_{\theta }(1 + y^{2} )
	+
	\log(N_{x} \tilde{C}_{\theta } \tilde{\varepsilon}_{\theta }^{-1} )
	+
	2\pi 
	\leq 
	K_{\theta } (1 + y^{2} ), 
	\\
	&
	\|\hat{G}_{\eta}(w,y) \|
	\leq 
	4 N_{x} \tilde{C}_{\theta } \tilde{\varepsilon}_{\theta }^{-2} 
	\leq 
	K_{\theta }
\end{align*}
for any $\eta \in V_{\delta_{\theta } }(\theta )$, 
$w \in V_{\delta_{\theta } }({\cal P}^{N_{x} } )$, 
$y \in {\cal Y}$. 
Hence, the lemma's assertion holds. 
\end{IEEEproof}

\begin{IEEEproof}[Proof of Proposition \ref{proposition4}]
Let $Q\subset\Theta$ be an arbitrary compact set. 
Then, owing to Conditions (i), (ii) of the proposition, 
there exists a real number 
$\varepsilon_{Q} \in (0,1)$ such that 
$\varepsilon_{Q} \leq p_{\alpha }(x'|x) \leq \varepsilon_{Q}^{-1}$
for all $\alpha \in {\cal A}$, $x,x' \in {\cal X}$
satisfying $[\alpha^{T} \beta^{T} ]^{T} \in Q$ for some 
$\beta \in {\cal B}$. 
Therefore, 
\begin{align*}
	\varepsilon_{Q} q_{\beta}(y|x') 
	\leq 
	r_{\theta }(y|x',x) 
	\leq 
	\varepsilon_{Q}^{-1} q_{\beta}(y|x')
\end{align*}
for all $\alpha \in {\cal A}$, $\beta \in {\cal B}$, 
$\theta = [\alpha^{T} \beta^{T} ]^{T}$, 
$x,x' \in {\cal X}$, $y \in {\cal Y}$ 
satisfying $\theta \in Q$. 
Thus, Assumption \ref{a3} is true. 

Since the collection of sets 
$\{V_{\delta_{\theta }/2 }(\theta ) \}_{\theta \in Q }$
covers 
$Q$ and 
since $Q$ is compact, 
there exists a finite subset $\tilde{Q}$ of 
$Q$ such that $Q$ is covered by 
$\{V_{\delta_{\theta }/2 }(\theta ) \}_{\theta \in \tilde{Q} }$. 
Let $\delta_{Q} = \min\{\delta_{\theta}/2: \theta \in \tilde{Q} \}$, 
$K_{Q} = \max\{K_{\theta}: \theta \in \tilde{Q} \}$
($\delta_{\theta }$, $K_{\theta }$ are defined in the statement of 
Lemma \ref{lemmaP}). 
Obviously, $\delta_{Q} \in (0,1)$, $K_{Q} \in [1,\infty )$. 
It can also be deduced that for each $\theta \in Q$, 
$V_{\delta_{Q} }(\theta ) \times V_{\delta_{Q} }({\cal P}^{N_{x} } )$
is contained in one of the sets from the collection 
$\{V_{\delta_{\theta } }(\theta ) \}_{\theta \in \tilde{Q} }$. 
Thus, 
$V_{\delta_{Q} }(Q) \times V_{\delta_{Q} }({\cal P}^{N_{x} } )
\subseteq 
\bigcup_{\theta \in \tilde{Q} } 
V_{\delta_{\theta } }(\theta ) \times V_{\delta_{\theta } }({\cal P}^{N_{x} } )$. 
Then, as an immediate consequence of Lemma \ref{lemmaP}, 
we have that Assumption \ref{a4} holds. 
\end{IEEEproof}

\section{Conclusion} 

We have studied the asymptotic properties of recursive maximum likelihood estimation
in hidden Markov models. 
We have analyzed the asymptotic behavior of the asymptotic log-likelihood function
and the convergence and convergence rate of 
the recursive maximum likelihood algorithm. 
Using the principle of analytical continuation, 
we have shown the analyticity of the asymptotic log-likelihood for analytically 
parameterized hidden Markov models. 
Relying on this result and Lojasiewicz inequality, 
we have demonstrated the point-convergence of the recursive maximum likelihood 
algorithm, 
and we have derived relatively tight bounds on the convergence rate. 
The obtained results cover a relatively broad class of hidden Markov models
with finite state space and continuous observations. 
They can also be extended to batch (i.e., non-recursive) 
maximum likelihood estimators such as those studied in 
\cite{bickel&ritov&ryden}, 
\cite{douc&moulines&ryden}, 
\cite{leroux}, 
\cite{ryden1}. 
In the future work, 
attention will be given to the possibility of extending 
the result of this paper to hidden Markov models with continuous state space.
The possibility of obtaining similar rate of convergence results for 
non-analytically parameterized hidden Markov models will be explored, too.

\end{document}